\documentclass[a4paper,12pt]{article}

\usepackage{amsmath}
\usepackage{amsfonts}
\usepackage{amssymb}
\usepackage{theorem}
\usepackage{graphicx}

\theorembodyfont{\it}
\newtheorem{theorem}{Theorem}[section]
\newtheorem{corollary}[theorem]{Corollary}

\newtheorem{lemma}[theorem]{Lemma}
\newtheorem{proposition}[theorem]{Proposition}
\newtheorem{star-theorem}{Theorem}

\theorembodyfont{\rm}
\newtheorem{definition}[theorem]{Definition}
\newtheorem{notation}[theorem]{Notation}

\renewenvironment{enumerate}
{
\begin{list}{\makebox[1.0cm]{\emph{(\roman{enumi})}}}
{
\usecounter{enumi}
\setlength{\topsep}{0pt}
\setlength{\partopsep}{0pt}
\setlength{\parsep}{0pt}
\setlength{\itemsep}{0pt}
\setlength{\labelsep}{0in}
\setlength{\labelwidth}{1.0cm}
\setlength{\leftmargin}{1.0cm}
}
}
{
\end{list}
}

\begin{document}

\title{On the Geometry of Goursat Structures \thanks{Submitted to:
ESAIM Control, Optimisation, and Calculus of Variations. Available as an e-print
at the Mathematics Archive front end: {\tt http://front.math.ucdavis.edu}.}}

\author{William {\sc Pasillas-L\'epine} and Witold {\sc Respondek} \bigskip \\
Institut national des sciences appliqu\'ees de Rouen\\
D\'epartement g\'enie math\'ematique\\
Place \'Emile Blondel --- 76 131 Mont Saint Aignan Cedex\\
{\tt wresp@lmi.insa-rouen.fr}}

\date{October 1999}

\maketitle

\begin{abstract}
A Goursat structure on a manifold of dimension$~n$ is a rank two distribution
$\mathcal{D}$ such that dim $\mathcal{D}^{(i)}=i+2$, for $i=0,...,n-2$, where
$\mathcal{D}^{(i)}$ denotes the derived flag of$~\mathcal{D}$, which is
defined by$~\mathcal{D}^{(0)}=\mathcal{D}$ and $\mathcal{D}^{(i+1)}%
=\mathcal{D}^{(i)}+[\mathcal{D}^{(i)},\mathcal{D}^{(i)}]$. Goursat structures
appeared first in the work of E. von Weber and E. Cartan, who have shown that
on an open and dense subset they can be converted into the so-called Goursat
normal form. Later, Goursat structures have been studied by Kumpera and Ruiz.
Contact structures on three manifolds and Engel structures on four manifolds
are examples of Goursat structures. In the paper, we introduce a new invariant
for Goursat structures, called the singularity type, and prove that the growth
vector and the abnormal curves of all elements of the derived flag are
determined by this invariant. Then we show, using a generalized version of
Backlund's theorem, that abnormal curves of all elements of the derived flag
do not determine the local equivalence class of a Goursat structure if $n>8$.
We also propose a new proof of a classical theorem of Kumpera and Ruiz. All
results are illustrated by the $n$-trailer system, which, as we show, turns
out to be a universal model for all local Goursat structures.
\end{abstract}

\newpage

\tableofcontents

\newpage

\section*{Introduction}

Let $\mathcal{D}$ be a smooth rank$~k$ distribution on a smooth manifold$~M$,
that is a map that assigns smoothly to each point$~p$ in$~M$ a linear subspace
$\mathcal{D}(p)\subset T_{p}M$ of dimension$~k$. The derived flag
of$~\mathcal{D}$ is the sequence defined by $\mathcal{D}^{(0)}=\mathcal{D}$
and $\mathcal{D}^{(i+1)}=\mathcal{D}^{(i)}+[\mathcal{D}^{(i)},\mathcal{D}%
^{(i)}]$, for $i\geq1$. A Goursat structure on a manifold $M$ of dimension
$n\geq3$ is a rank two distribution $\mathcal{D}$ such that, for $0\leq i\leq
n-2$, the elements of its derived flag satisfy $\dim\mathcal{D}^{(i)}(p)=i+2$,
for each point $p$ in $M$. Goursat structures were introduced, using the dual
language of Pfaffian systems by E. von Weber in 1898. The first period of
interest in this special class of distributions culminated in the work of
Cartan and Goursat. A new period was initiated by Giaro, Kumpera, and Ruiz at
the end of the seventies. A renewal of interest in Goursat structures has been
observed from the mid of nineties.

There are at least three reasons explaining those one century long studies.
The first reason is that any Goursat structure on $\mathbb{R}^{n}$ can be
locally converted (on an open and dense subset, as it was observed only later
by Giaro, Kumpera, and Ruiz~\cite{giaro-kumpera-ruiz}) into the so-called
Goursat normal form, also known as chained form:
\[
\left(
\begin{array}
[c]{c}%
\tfrac{\partial}{\partial x_{n}}%
\end{array}
,
\begin{array}
[c]{c}%
x_{n}\tfrac{\partial}{\partial x_{n-1}}+\cdots+x_{3}\tfrac{\partial}{\partial
x_{2}}+\tfrac{\partial}{\partial x_{1}}%
\end{array}
\right)  .
\]
It seems that Weber~\cite{weber-article} was the first to exhibit this
property and, indeed, Goursat~\cite{goursat} attributes to him this result. In
fact, the starting point of Weber's studies was the following question: ``When
is a given distribution equivalent to Goursat normal form?'', which had led
him to discover the concept of derived flag. This question is very natural
because Goursat normal form has a clear geometric interpretation. Indeed, let
us consider the space $J^{k}(\mathbb{R},\mathbb{R})$ of $k$-jets of maps from
$\mathbb{R}$ to $\mathbb{R}$. On the one hand, a necessary condition for a
curve in $J^{k}(\mathbb{R},\mathbb{R})$ to be a prolongation of a graph of a
function from $\mathbb{R}$ to $\mathbb{R}$ is that it is an integral curve of
a distribution which, in the canonical coordinates of $J^{k}(\mathbb{R}%
,\mathbb{R})$, is spanned by the Goursat normal form on $\mathbb{R}^{k+2}$. On
the other hand, a necessary and sufficient condition for a diffeomorphism of
$J^{k}(\mathbb{R},\mathbb{R})$ to map prolongations of graphs of functions
into prolongations of graphs of functions is to preserve the distribution
spanned by the Goursat normal form on$~\mathbb{R}^{k+2}$. Such diffeomorphisms
are called contact transformations~\cite{olver-equivalence} of order$~k$ and
have been intensively studied by B\"{a}cklund~\cite{backlund}, and by Lie and
Scheffers~\cite{lie-scheffers}.

The second reason of interest in Goursat structures has been the classical
problem, first considered by Monge, of integrating underdetermined
differential equations (equivalently, Pfaffian systems) without integration.
To be more precise, let $\mathcal{D}$ be a rank $k$ distribution on $M$. The
problem (see e.g.~\cite{goursat-monge} and~\cite{zervos}) is to find~$k$
smooth functions $\varphi_{1},\ldots,\varphi_{k}$ such
that any integral curve $\gamma(t)$ of $\mathcal{D}$ can be expressed as a
smooth function of $\varphi_{1},\ldots,\varphi_{k}$ and their time-derivatives
along $\gamma(t)$. The most important achievement of the
first period of studies on Goursat structures was a result of E.
Cartan~\cite{cartan-equivalence-absolue}, who showed that a rank two
distribution posses the above described property if and only if it is
transformable into Goursat normal form.

The third reason of importance of Goursat structures is that they describe the
nonholonomic constraints of many mechanical systems. For example, the
kinematical constraints of a passenger car are described by a Goursat
structure on $\mathbb{R}^{2}\times(S^{1})^{2}$; those of a truck by a Goursat
structure on $\mathbb{R}^{2}\times(S^{1})^{3}$. Moreover, for Goursat
structures the nonholonomic motion planning problem can be solved explicitly;
either by transforming them into Goursat normal form, as suggested by Murray
and Sastry (see e.g.~\cite{murray-nilpotent} and~\cite{murray-sastry}), or by
using the concept of flatness, introduced in control theory by Fliess,
L\'{e}vine, Martin, and Rouchon, which is the above described property of
calculating the trajectories without integration (see
e.g.~\cite{martin-rouchon-driftless} and~\cite{fliess-intro-flat}).

As we said, the second period of studies on Goursat structures began with a
work of Giaro, Kumpera, and Ruiz~\cite{giaro-kumpera-ruiz}, who observed that
there are Goursat structures which are not locally equivalent to Goursat
normal form. This observation raised the problem of classification of Goursat
structures and that of finding their invariants, and has led Kumpera and Ruiz
to write their important paper~\cite{kumpera-ruiz}, where they gave a complete
classification, up to dimension~$7$, together with a set of general results on
Goursat structures.

In the nineties, research on Goursat structures was concentrated around two
main topics: the classification problem and the nonholonomic motion planning
problem for mechanical systems described by Goursat structures. Among results
concerning the classification problem, Murray~\cite{murray-nilpotent}
obtained, using the concept of growth vector\footnote{The Lie flag of a
distribution $\mathcal{D}$ is the sequence defined by $\mathcal{D}%
_{0}=\mathcal{D}$ and $\mathcal{D}_{i+1}=\mathcal{D}_{i}+[\mathcal{D}%
_{0},\mathcal{D}_{i}]$, for $i\geq1$. The sequence $(\dim\mathcal{D}%
_{i}(p))_{i\geq0}$ is called the growth vector of $\mathcal{D}$ at $p$.}, an
easily checkable necessary and sufficient condition for a Goursat structure to
be equivalent to Goursat normal form (his condition simplifies those of
Libermann~\cite{libermann} and Kumpera and Ruiz~\cite{kumpera-ruiz}). Cheaito
and Mormul~\cite{cheaito-mormul} corrected the classification in dimension~$7$
(see also \cite{gaspar}) and obtained a complete classification in
dimension~$8$. Mormul~\cite{mormul-R9} obtained a complete classification in
dimension$~9$. It turns out that this dimension is the highest one in which
there is a finite number of non-equivalent Goursat structures.\ Indeed,
Cheaito, Mormul, and the authors~\cite{cheaito-mormul-pasillas-respondek}
showed that in higher dimensions there are real continuous parameters in the
classification. Note, however, that in each dimension all Goursat structures
are finitely determined, which implies that there are no functional parameters
in the classification.

Most of the work concerning mechanical control systems described by Goursat
structures has been motivated by the study of the $n$-trailer system. It would
be impossible to give here a complete set of references on this subject. We
have thus chosen to cite two books \cite{laumond-book,li-canny-book}, and to
give a few references concerning standard control theory problems for the
$n$-trailer system and chained systems.

The controllability of the $n$-trailer system has been proved by Laumond both
for regular~\cite{laumond-trailer} and singular~\cite{laumond-singularities}
configurations. Improved bounds for the nonholonomy degree of the $n$-trailer
at singular configurations have been obtained by S\o
rdalen~\cite{sordalen-bound}, Luca and Risler~\cite{luca-risler}, and
Jean~\cite{jean-trailer}. For regular configurations, an explicit conversion
of the $n$-trailer system into chained form has been obtained by S\o rdalen~
\cite{sordalen-trailer}; for singular configurations, an explicit conversion
of the $n$-trailer system into Kumpera-Ruiz normal form has been obtained by
the authors~\cite{pasillas-respondek-cdc} (see also Section~\ref{sec-trailer}).

Open loop motion planning has been investigated for general nonholonomic
systems by Brockett~\cite{brockett-singular}, Lafferriere and
Sussmann~\cite{laferriere-sussmann}, and Liu~\cite{liu-approximation}. For
chained systems, these results have been considerably simplified by Murray and
Sastry~\cite{murray-sastry} (using the special properties of chained form) and
by Fliess \emph{et al.}~\cite{fliess-intro-flat} (using the concept of
flatness). Combined with the conversion of the $n$-trailer into chained form
obtained by S\o rdalen~\cite{sordalen-trailer}, they have led to a solution of
the nonholonomic motion planning problem for the $n$-trailer system
(see~e.g.~\cite{laumond-jacobs-taix-murray}, \cite{fliess-trailer},
and~\cite{tilbury-murray-sastry}).

Path tracking of non-abnormal open loop trajectories has been studied by
Fliess \emph{et al.}~\cite{fliess-intro-flat}, Jiang and
Nijmeijer~\cite{jiang-nijmeijer}, and Walsh~\emph{et~al.}%
~\cite{walsh-tilbury-sastry-murray-laumond}. Since for chained systems
constant trajectories (points) are abnormal, the proposed path tracking
strategies cannot be applied to achieve pointwise stabilization. Indeed, the
linearization of a chained system around such trajectories is not
controllable. The first who observed the difficulties of pointwise
stabilization for control systems without drift was
Brockett~\cite{brockett-stabilization}. General algorithms for pointwise
stabilization of nonholonomic systems can be found in the work of
Coron~\cite{coron-driftless}, Pomet~\cite{pomet-driftless}, McCloskey and
Murray~\cite{closkey-murray}, and Morin and Samson~\cite{morin-samson-cocv}.
For chained systems, improved results have been obtained by
Samson~\cite{samson-chained}, S\o rdalen and Egeland~\cite{sordalen-egeland},
and Teel~\emph{et al.}~\cite{teel-murray-walsh}. These methods have been
successfully applied to the $n$-trailer system (see e.g.
Samson~\cite{samson-chained}, S\o rdalen and Wichlund~\cite{sordalen-wichlund}%
, and the references given there).

\newpage 

\noindent Our paper is devoted to a study of the geometry of Goursat
structures.\ We will introduce a new local invariant for Goursat structures,
called the singularity type, whose definition is based on the following
observation.\ If $\mathcal{D}$ is a Goursat structure then each element
$\mathcal{D}^{(i)}$ of its derived flag contains an involutive subdistribution
$\mathcal{C}_{i}\subset\mathcal{D}^{(i)}$ that has constant corank one in
$\mathcal{D}^{(i)}$ and is characteristic for $\mathcal{D}^{(i+1)}$. The
singularity type reflects the geometry of incidence between the distributions
$\mathcal{D}^{(i)}$ and the distributions$~\mathcal{C}_{i}$. Although, as we
prove, the singularity type keeps the same information about a Goursat
structure as the growth vector, that information is encoded in the singularity
type in a much more systematic and, what is extremely important, in a much
more geometric way.\ In particular, the geometric information contained in the
singularity type enables us to describe completely all abnormal curves of all
elements of the derived flag.\ This can be summarized in the following
Theorem, which is a combination of
Theorem~\ref{thm-singularity-type-growth-vector}
and~Theorem~\ref{thm-singularity-type-abnormal}, and gives one of the main
results of the paper.

\begin{star-theorem}
Let $\mathcal{D}$\ and $\tilde{\mathcal{D}}$\ be two Goursat structures
defined on two manifolds~$M$\ and~$\tilde{M}$, respectively, of dimension
$n\geq3$. Fix two points $p$\ and $\widetilde{p}$\ of~$M$\ and~$\tilde{M}$,
respectively. The three following conditions are equivalent:

\begin{enumerate}
\item  The singularity type of $\mathcal{D}$\ at $p$\ equals the singularity
type of $\tilde{\mathcal{D}}$\ at $\tilde{p}$.

\item  The growth vector of $\mathcal{D}$\ at $p$\ equals the growth vector of
$\tilde{\mathcal{D}}$\ at $\tilde{p}$.

\item  There exists a diffeomorphism $\varphi$, with $\tilde{p}=\varphi(p)$,
between two small enough neighborhoods of $p$\ and $\tilde{p}$\ that
transforms the abnormal curves of $\mathcal{D}^{(i)}$\ into the abnormal
curves of $\tilde{\mathcal{D}}^{(i)}$, for each~$i\geq0$.
\end{enumerate}
\end{star-theorem}

An important example that we will use to illustrate our results on Goursat
structures will be the $n$-trailer system, that is a mobile robot (unicycle)
towing$~n$ trailers. In the paper we will calculate rigid curves of the
$n$-trailer and give their natural mechanical interpretation: they correspond
to motions that fix the positions of the centers of at least two trailers. We
will also show how to transform locally the $n$-trailer system into a
Kumpera-Ruiz normal form, and we will prove a surprising result stating that
any Goursat structure is locally equivalent to the $n$-trailer system around a
well chosen point of its configuration space. This result will enable us to
use for any Goursat structure a deep result of Jean~\cite{jean-trailer}
devoted to singular configurations of the $n$-trailer system, in particular we
will extend to all Goursat structures Jean's formula for the growth vector of
the $n$-trailer system. In our work, the singularity type will replace the
angles of the $n$-trailer system that appear in Jean's theorem.

In the paper, we will propose an inductive procedure of constructing
Kumpera-Ruiz normal forms of Goursat structures based on two types of
prolongations: regular and singular. This construction provides a systematic
and unifying approach to many results of the paper. In particular, it will be
used to show that any Goursat structure can be brought to a Kumpera-Ruiz
normal form; to study generalized contact transformations, that is
transformations which preserve Goursat structures; and to define the above
mentioned transformations that transform locally the $n$-trailer system into a
Kumpera-Ruiz normal form, and, conversely, that convert locally an arbitrary
Goursat structure into the $n$-trailer system around a well chosen point of
its configuration space.

Recent studies (see \cite{jakubczyk-complex-abnormals} and
\cite{montgomery-survey}) show that most distributions are determined by their
abnormal curves. Our complete description of abnormal curves of Goursat
structures enables us to conclude that this is not the case of Goursat
structures.\ Indeed, combining our study with the main theorem
of~\cite{cheaito-mormul-pasillas-respondek} leads us to the following result:
There exist two Goursat structures $\mathcal{D}$\ and $\tilde{\mathcal{D}}%
$\ defined at $p$\ and $\tilde{p}$, respectively, that are not locally
equivalent but for which there exists a diffeomorphism $\varphi$, with
$\tilde{p}=\varphi(p)$, between two small enough neighborhoods of $p$\ and
$\tilde{p}$\ that transforms the abnormal curves of $\mathcal{D}^{(i)}$\ into
the abnormal curves of$~\tilde{\mathcal{D}}^{(i)}$, for each $i\geq0$ (see
Proposition~\ref{prop-R9} and Proposition~\ref{prop-R11}).

The paper is organized as follows. In the first Section we introduce Goursat
structures, we give some examples in small dimension, and we define Goursat
normal form. In the second Section, we provide an inductive definition of
Kumpera-Ruiz normal form. The proposed concept of prolongations enable us to
give a new proof of the Kumpera-Ruiz theorem, which states that any Goursat
structure can be converted locally into a Kumpera-Ruiz normal form.\ In the
third Section, we introduce the $n$-trailer system and we construct
transformations that bring locally the $n$-trailer system into a Kumpera-Ruiz
normal form and, conversely, that bring an arbitrary Goursat structure into
the $n$-trailer system. In Section four, we introduce our main invariant of
Goursat structures, namely, the singularity type. We also compute the
singularity type for Kumpera-Ruiz normal forms and for the $n$-trailer
system.\ As we have said, the singularity type keeps the same information
about Goursat structures as the growth vector although in both invariants that
information is encoded in a different way. Section five is devoted to study
relations between these two invariants.\ In particular, we give a formula to
compute the growth vector of an arbitrary Goursat structure and another to
compute the singularity type using the growth vector. In Section six we study
abnormal curves of Goursat structures.\ We give a complete description of
absolutely continuous abnormal curves for all elements of the derived flag of
any Goursat structure.\ We prove that the whole information about all abnormal
curves is given by the singularity type.\ In Section seven we study
generalized contact transformations, that is transformations which preserve
Goursat structures (also called symmetries) and we give formulas to calculate
them from first order contract transformations.\ Those formulas are used to
analyze examples of Goursat structures that are non-equivalent but that have
diffeomorphic abnormal curves.\ The paper ends with three Appendices.\ The
first is devoted to a class of distributions that, although of rank greater
than two, are very close to Goursat structures.\ This class was also studied
by E.\ von Weber. In the second Appendix we construct a normal form which we
use in our study of rigidity of integral curves of Goursat
structures.\ Finally, in the third Appendix, we illustrate through a set of
figures different configurations of the $n$-trailer system corresponding to
various Kumpera-Ruiz normal forms in dimensions$~3$, $4$, $5$, and$~6$.

\newpage

\section{Goursat Structures}

\label{sec-goursat}

\subsection{Derived Flag and Goursat Structures}

A rank\emph{\ }$k$ \emph{distribution }$\mathcal{D}$ on a smooth manifold $M$
is a map that assigns smoothly to each point~$p$ in$~M$ a linear subspace
$\mathcal{D}(p)\subset T_{p}M$ of dimension~$k$. In other words, a rank~$k$
distribution is a smooth rank $k$ subbundle of the tangent bundle $TM$. Such a
field of tangent $k$-planes is spanned locally by $k$ pointwise linearly
independent smooth vector fields $f_{1},\ldots,f_{k}$ on $M$, which will be
denoted by $\mathcal{D}=(f_{1},\ldots,f_{k})$.

Two distributions $\mathcal{D}$ and $\tilde{\mathcal{D}}$ defined on two
manifolds $M$ and $\tilde M$, respectively, are \emph{equivalent} if there
exists a smooth diffeomorphism $\varphi$ between $M$ and $\tilde M$ such that
\[
(\varphi_{*}\mathcal{D})(\tilde p)=\tilde{\mathcal{D}}(\tilde p),
\]
for each point$~\tilde p$ in $\tilde M$. They are \emph{locally equivalent} at
two points $p$ and $\tilde p$ that belong to $M$ and $\tilde M$, respectively,
if there exists two small enough neighborhoods $U$ and$~\tilde U$ of $p$ and
$\tilde p$, respectively, and a diffeomorphism$~\varphi$ between $U$ and
$\tilde U$ such that $\varphi(p)=\tilde p$ and $(\varphi_{*}\mathcal{D}%
)(\tilde p)=\tilde{\mathcal{D}}(\tilde p)$, for each point$~\tilde p$ in
$\tilde U$.

The \emph{derived flag} of a distribution $\mathcal{D}$ is the sequence
$\mathcal{D}^{(0)}\subset\mathcal{D}^{(1)}\subset\cdots$ defined inductively
by
\begin{equation}
\mathcal{D}^{(0)}=\mathcal{D}\text{\quad and\quad}\mathcal{D}^{(i+1)}%
=\mathcal{D}^{(i)}+[\mathcal{D}^{(i)},\mathcal{D}^{(i)}]\text{, \quad for
}i\geq0\text{.}\label{derived-flag}%
\end{equation}
This sequence should not be confused with the Lie flag (\ref{lie-flag}), which
will be introduced in Section~\ref{sec-growth}.

\begin{definition}
\label{def-goursat-structure}A \emph{Goursat structure} on a manifold $M$ of
dimension $n\geq3$ is a rank two distribution $\mathcal{D}$ such that, for
$0\leq i\leq n-2$, the elements of its derived flag satisfy $\dim
\mathcal{D}^{(i)}(p)=i+2$, for each point $p$ in $M$.
\end{definition}

\subsection{Examples of Goursat Structures}

A Goursat structure on a three-manifold is a \emph{contact structure,} and is
locally equivalent to the distribution spanned by
\[
\left(
\begin{array}
[c]{c}%
\tfrac{\partial}{\partial x_{3}}%
\end{array}
,
\begin{array}
[c]{c}%
x_{3}\tfrac{\partial}{\partial x_{2}}+\tfrac{\partial}{\partial x_{1}}%
\end{array}
\right)  ,
\]
which is called \emph{Pfaff-Darboux normal form}. A Goursat structure on a
four-manifold is an \emph{Engel structure,} and is locally equivalent to the
distribution spanned by
\[
\left(
\begin{array}
[c]{c}%
\tfrac{\partial}{\partial x_{4}}%
\end{array}
,
\begin{array}
[c]{c}%
x_{4}\tfrac{\partial}{\partial x_{3}}+x_{3}\tfrac{\partial}{\partial x_{2}%
}+\tfrac{\partial}{\partial x_{1}}%
\end{array}
\right)  ,
\]
which is called \emph{Engel normal form}. Observe that, for a generic field of
planes$~\mathcal{D}$ on$~\mathbb{R}^{3}$, we have $\dim\mathcal{D}^{(1)}%
(p)=3$, for any point$~p$ on an open and dense subset of$~\mathbb{R}^{3}$; for
a generic field of planes$~\mathcal{D}$ on$~\mathbb{R}^{4}$, we have
$\dim\mathcal{D}^{(1)}(p)=3$ and $\dim\mathcal{D}^{(2)}(p)=4$, for any
point$~p$ on an open and dense subset of$~\mathbb{R}^{4}$. Therefore, in a
small enough neighborhood of a generic point, a generic field of planes on a
manifold of dimension three or four is a Goursat structure. Note, however,
that starting from dimension five the class of Goursat structures is of
infinite codimension within the class of all rank two distributions. Indeed,
for a generic field of planes$~\mathcal{D}$ on$~\mathbb{R}^{n}$, for $n\geq5$,
we have $\dim\mathcal{D}^{(1)}(p)=3$ and $\dim\mathcal{D}^{(2)}(p)=5$, for any
point$~p$ on an open and dense subset of$~\mathbb{R}^{n}$.

We give now a mechanical example. For $n\geq0$, the distribution spanned by
the following pair of vector fields:
\begin{equation}
\left(
\begin{array}
[c]{c}%
\tfrac\partial{\partial\theta_{n}}%
\end{array}
,
\begin{array}
[c]{c}%
\cos(\theta_{0})\pi_{0}\tfrac\partial{\partial\xi_{1}}+\sin(\theta_{0})\pi
_{0}\tfrac\partial{\partial\xi_{2}}+%
{\textstyle\sum\limits_{i=0}^{n-1}}
\sin(\theta_{i+1})\pi_{i+1}\tfrac\partial{\partial\theta_{i}}%
\end{array}
\right)  ,\label{nonholonomic-manipulator}%
\end{equation}
where $\pi_{i}=%
{\textstyle\prod\nolimits_{j=i+1}^{n}}
\cos(\theta_{j})$ and $\pi_{n}=1$, is a Goursat structure on $\mathbb{R}%
^{2}\times(S^{1})^{n+1}$ equipped with coordinates $(\xi_{1},\xi_{2}%
,\theta_{0},\ldots,\theta_{n})$. This distribution is the kinematical model
for the ``nonholonomic manipulator'' of S\o rdalen, Nakamura, and
Chung~\cite{sordalen-nakamura-chung}. Another example is the $n $-trailer
system (see Section~\ref{sec-trailer}), which will play a fundamental role in
this article.

\subsection{Goursat Normal Form}

The concepts of derived flag and Goursat structure were introduced, using the
dual language of Pfaffian systems, by E. von Weber \cite{weber-article} in
order to characterize the class of Pfaffian systems that can be converted into
a particular normal form, also introduced by him, which is now known as
Goursat normal form (see (\ref{goursat-normal-form}) below; see also
Appendix~\ref{sec-weber}).

Although it is not clear who was the first to prove that Goursat structures
can be converted locally into Goursat normal form, at least on an open and
dense subset \cite[Theorem VI]{weber-article} (compare
\cite{bryant-chern-gardner-goldschmidt-griffiths},
\cite{cartan-equivalence-absolue}, \cite{goursat}, \cite{kumpera-ruiz},
\cite{murray-nilpotent}, and \cite{weber-article}); it is clear that the
foundations of this result were set by Engel and Weber (see
\cite{cartan-weber}, \cite{engel}, \cite{weber-article}, and the references
given there; see also Appendix \ref{sec-weber}). The importance of this result
was, however, fully understood only later, by E.~Cartan when he solved a long
standing problem of that time: the characterization of explicitly integrable
Monge equations~\cite{cartan-equivalence-absolue} (see also
\cite{bryant-chern-gardner-goldschmidt-griffiths}, \cite{giaro-kumpera-ruiz},
\cite{goursat}, \cite{hilbert}, \cite{martin-rouchon-driftless}, and
\cite{zervos}).

\begin{theorem}
[Weber-Cartan-Goursat]\label{thm-weber-cartan-goursat}Let $\mathcal{D}$ be a
Goursat structure defined on a manifold $M$ of dimension $n\geq3$. There
exists an open and dense subset $U\subset M$ such that, for any point $p$ in
$U$, the distribution $\mathcal{D}$ is locally equivalent at $p$ to the
distribution spanned by the following pair of vector fields:
\begin{equation}
\left(
\begin{array}
[c]{c}%
\tfrac{\partial}{\partial x_{n}}%
\end{array}
,
\begin{array}
[c]{c}%
x_{n}\tfrac{\partial}{\partial x_{n-1}}+x_{n-1}\tfrac{\partial}{\partial
x_{n-2}}+\cdots+x_{3}\tfrac{\partial}{\partial x_{2}}+\tfrac{\partial
}{\partial x_{1}}%
\end{array}
\right)  ,\label{goursat-normal-form}%
\end{equation}
considered on a small enough neighborhood $V\subset\mathbb{R}^{n}$ of zero.
Moreover, if $n=3$ or$~4$ then $U=M$.
\end{theorem}

\noindent\noindent In control theory, the normal form
(\ref{goursat-normal-form}) is usually called \emph{chained form}. In order to
keep the classical name, we will call it \emph{Goursat normal form}. An
elegant characterization, using the growth vector, of the above mentioned open
and dense set\noindent$~U$ was obtained by Murray~\cite{murray-nilpotent} (see
Theorem~\ref{thm-murray}). Observe that, in most of the above mentioned
references, Goursat structures are not defined by distributions but by their
duals, that is by Pfaffian systems. Note also that many other names have been
given to Goursat structures: ``systeme vom Charakter eins und vom Rang
zwei''~\cite{weber-article}, ``syst\`{e}mes de classe z\'{e}%
ro''~\cite{cartan-equivalence-absolue}, ``syst\`{e}mes sp\'{e}%
ciaux''~\cite{goursat}, ``syst\`{e}mes en drapeaux''~\cite{kumpera-ruiz},
``systems of Goursat type''~ \cite{bryant-hsu}, and ``systems that satisfy the
Goursat condition''~\cite{mormul-R9}.

\section{Kumpera-Ruiz's Theorem}

\label{sec-kumpera-ruiz}

\subsection{Kumpera-Ruiz Normal Forms}

If at a given point a Goursat structure can be converted into Goursat normal
form then this point is called \emph{regular}; otherwise, it is called
\emph{singular}. The first who observed the existence of singular points were
Giaro, Kumpera and Ruiz \cite{giaro-kumpera-ruiz}. This initial observation
has led Kumpera and Ruiz to write their pioneering paper \cite{kumpera-ruiz},
where they introduced the normal forms that we will consider in this section.

We start with the Pfaff-Darboux and Engel normal forms, given respectively on
$\mathbb{R}^{3}$, equipped with coordinates $(x_{1},x_{2},x_{3})$, and
$\mathbb{R}^{4}$, equipped with coordinates $(x_{1},x_{2},x_{3},x_{4})$, by
the pairs of vector fields $\kappa^{3}=(\kappa_{1}^{3},\kappa_{2}^{3})$ and
$\kappa^{4}=(\kappa_{1}^{4},\kappa_{2}^{4})$, defined by
\[%
\begin{array}
[c]{l}%
\kappa_{1}^{3}=\tfrac\partial{\partial x_{3}}\\
\kappa_{2}^{3}=x_{3}\tfrac\partial{\partial x_{2}}+\tfrac\partial{\partial
x_{1}}%
\end{array}
\qquad\qquad
\]
and
\[%
\begin{array}
[c]{l}%
\kappa_{1}^{4}=\tfrac\partial{\partial x_{4}}\\
\kappa_{2}^{4}=x_{4}\tfrac\partial{\partial x_{3}}+x_{3}\tfrac\partial
{\partial x_{2}}+\tfrac\partial{\partial x_{1}}.
\end{array}
\]
Loosely speaking, we can write
\[%
\begin{array}
[c]{l}%
\kappa_{1}^{4}=\tfrac\partial{\partial x_{4}}\\
\kappa_{2}^{4}=x_{4}\kappa_{1}^{3}+\kappa_{2}^{3}.
\end{array}
\qquad\qquad\qquad
\]
In order to make this precise we will adopt the following natural notation.
Consider a vector field
\[
f^{n-1}=%
{\textstyle\sum_{i=1}^{n-1}}
f_{i}^{n-1}(x_{1},\ldots,x_{n-1})\tfrac\partial{\partial x_{i}}
\]
on $\mathbb{R}^{n-1}$ equipped with coordinates $(x_{1},\ldots,x_{n-1})$. We
can lift $f^{n-1}$ to a vector field, denoted also by $f^{n-1}$, on
$\mathbb{R}^{n}$ equipped with coordinates $(x_{1},\ldots,x_{n-1},x_{n})$ by
taking
\[
f^{n-1}=%
{\textstyle\sum_{i=1}^{n-1}}
f_{i}^{n-1}(x_{1},\ldots,x_{n-1})\tfrac\partial{\partial x_{i}}+0\cdot
\tfrac\partial{\partial x_{n}}.
\]
That is, we lift $f^{n-1}$ by translating it along the $x_{n}$-direction.

\begin{notation}
\label{not-lift}From now on, in any expression of the form $\kappa_{2}%
^{n}=\alpha(x)\kappa_{1}^{n-1}+\beta(x)\kappa_{2}^{n-1}$, the vector fields
$\kappa_{1}^{n-1}$ and $\kappa_{2}^{n-1}$ should be considered as the above
defined lifts of $\kappa_{1}^{n-1}$ and $\kappa_{2}^{n-1}$, respectively.
\end{notation}

Let $\kappa^{n-1}=(\kappa_{1}^{n-1},\kappa_{2}^{n-1})$ denote a pair of vector
fields on $\mathbb{R}^{n-1}$. A \emph{regular prolongation},\ with parameter
$c_{n}$, of $\kappa^{n-1}$, denoted by $\kappa^{n}=R_{c_{n}}(\kappa^{n-1})$,
is a pair of vector fields $\kappa^{n}=(\kappa_{1}^{n},\kappa_{2}^{n})$
defined on $\mathbb{R}^{n}$ by
\begin{equation}%
\begin{array}
[c]{l}%
\kappa_{1}^{n}=\tfrac\partial{\partial x_{n}}\\
\kappa_{2}^{n}=(x_{n}+c_{n})\kappa_{1}^{n-1}+\kappa_{2}^{n-1},
\end{array}
\label{regular-prolongation}%
\end{equation}
where $c_{n}$ belongs to $\mathbb{R}$. The \emph{singular prolongation} of
$\kappa^{n-1}$, denoted by $\kappa^{n}=S(\kappa^{n-1})$, is the pair of vector
fields $\kappa^{n}=(\kappa_{1}^{n},\kappa_{2}^{n})$ defined on $\mathbb{R}%
^{n}$ by
\begin{equation}%
\begin{array}
[c]{l}%
\kappa_{1}^{n}=\frac\partial{\partial x_{n}}\\
\kappa_{2}^{n}=\kappa_{1}^{n-1}+x_{n}\kappa_{2}^{n-1}.
\end{array}
\qquad\quad\label{singular-prolongation}%
\end{equation}

\begin{definition}
\label{def-kumpera-ruiz}A pair of vector fields $\kappa^{n}$ on $\mathbb{R}%
^{n} $, for $n\geq3$, is called a \emph{Kumpera-Ruiz normal form} if
$\kappa^{n}=\sigma_{n-3}\circ\cdots\circ\sigma_{1}(\kappa^{3})$, where each
$\sigma_{i}$, for $1\leq i\leq n-3$, equals either $S$ or $R_{c_{i}}$, for
some real constants~$c_{i}$.
\end{definition}

We will also call a Kumpera-Ruiz normal form the restriction of a Kumpera-Ruiz
normal form to any open subset of $\mathbb{R}^{n}$. Most often, the
coordinates $x_{1},\ldots,x_{n}$ will be the elements of a coordinate chart
$x:M\rightarrow\mathbb{R}^{n}$, defined in a neighborhood of a given point $p$
in $M$. If we have $x(p)=0$ then we will say that the Kumpera-Ruiz normal form
is \emph{centered} at $p$. For example, the two models considered in
\cite{giaro-kumpera-ruiz}:
\begin{align*}
& \left(
\begin{array}
[c]{c}%
\tfrac\partial{\partial x_{5}}%
\end{array}
,
\begin{array}
[c]{c}%
x_{5}\tfrac\partial{\partial x_{4}}+x_{4}\tfrac\partial{\partial x_{3}}%
+x_{3}\tfrac\partial{\partial x_{2}}+\tfrac\partial{\partial x_{1}}%
\end{array}
\right) \\
& \left(
\begin{array}
[c]{c}%
\tfrac\partial{\partial x_{5}}%
\end{array}
,
\begin{array}
[c]{c}%
\tfrac\partial{\partial x_{4}}+x_{5}\left(  x_{4}\tfrac\partial{\partial
x_{3}}+x_{3}\tfrac\partial{\partial x_{2}}+\tfrac\partial{\partial x_{1}%
}\right)
\end{array}
\right)  ,
\end{align*}
defined by $R_{0}(\kappa^{4})$ and $S(\kappa^{4})$, respectively, are
Kumpera-Ruiz normal forms on $\mathbb{R}^{5}$, equipped with coordinates
$(x_{1},\ldots,x_{5})$, centered at zero.

\subsection{Kumpera-Ruiz's Theorem}

The following theorem of Kumpera and Ruiz shows clearly the importance of
their normal forms. We will give a proof of this theorem at the end of this
Section because many of our results are based on a construction that also
appears in our proof. Moreover, we would like to emphasize two features of our
proof. Firstly, it is quite close to the original ideas of E. von Weber.
Indeed, though we use distributions instead of Pfaffian systems, the two
proofs share the same fundamental Lemma (compare \cite[Theorem V]%
{weber-article} and Proposition~\ref{prop-extended-engel}; see also
Appendix~\ref{sec-weber}). Secondly, it is to our knowledge the only one that
does not mix the language of vector fields and differential forms (everything
is done in terms of vector fields). For alternative proofs we refer the reader
to \cite{cheaito-mormul} and to the original work of Kumpera and
Ruiz~\cite{kumpera-ruiz}.

\begin{theorem}
[Kumpera-Ruiz]\label{thm-kumpera-ruiz}Any Goursat structure on a manifold $M$
of dimension~$n$ is locally equivalent, at any point $p$ in $M$, to a
distribution spanned by a Kumpera-Ruiz normal form centered at $p$ and defined
on a suitably chosen neighborhood $U\subset\mathbb{R}^{n}$ of zero.
\end{theorem}

This theorem is the cornerstone to understand the geometry of Goursat
structures. On the one hand, this result implies that locally, even at
singular points, Goursat structures do not have \emph{functional} invariants;
this property makes them precious but rare and distinguishes them from generic
rank two distributions on $n$-manifolds, which do have local functional
invariants when $n\geq5$ (see \cite{cartan-cinq-variables},
\cite{jakubczyk-przytycki}, \cite{vershik-gershkovich},
\cite{zhitomirskii-survey}). But on the other hand, the real constants that
appear in Kumpera-Ruiz normal forms are unavoidable; this fact has been
observed only recently and implies that Goursat structures do have \emph{real}
invariants (see \cite{cheaito-mormul-pasillas-respondek}, \cite{mormul-R9},
and Section \ref{sec-contact}).

Though our definition of Kumpera-Ruiz normal forms was inductive, it is
possible to give the following equivalent explicit
definition~(\ref{cras-normal-form}), which will also be used in the paper.
Observe that in the normal form~(\ref{cras-normal-form}), we use a double
indexation$~x_{j}^{i}$ of coordinates, for $0\leq i\leq m+1$, where the
integer $m$ gives the \emph{number of singularities} of the normal form, that
is the number of singular prolongations (provided that $\sigma_{1}$ is
regular, which can always be assumed without lose of generality).

\begin{corollary}
\label{cor-cras-indexation} Any Goursat structure defined on a manifold $M$ of
dimension $n\geq4$ is locally equivalent, at any point $p$ in $M$, to a
distribution spanned in a small neighborhood of zero by a pair of vector
fields that has the following form:
\begin{equation}
(\;\tfrac{\partial}{\partial x_{1}^{0}}\;,\;%
{\textstyle\sum\limits_{i=0}^{m}}
\,(%
{\textstyle\prod\limits_{j=0}^{i-1}}
x_{k_{j}}^{j})(%
{\textstyle\sum\limits_{j=1}^{k_{i}-1}}
(x_{j}^{i}+c_{j}^{i})\tfrac{\partial}{\partial x_{j+1}^{i}}+\tfrac{\partial
}{\partial x_{1}^{i+1}})\;),\label{cras-normal-form}%
\end{equation}
where the coordinates $x_{j}^{i}$, for $0\leq i\leq m+1$ and $1\leq j\leq
k_{i}$, are centered at $p$; the integer $m$ is such that $0\leq m\leq n-4$;
and the integers $k_{i}$, for $0\leq i\leq m+1$, satisfy $k_{0}\geq
1,\ldots,k_{m-1}\geq1$, $k_{m}\geq3$, $k_{m+1}=1$ and $\sum_{i=0}^{m+1}%
k_{i}=n$. The constants $c_{j}^{i}$, for $1\leq j\leq k_{i}-1$, are real constants.
\end{corollary}

\subsection{Low Dimensional Examples}

Let us recall the complete local classification of Goursat structures on
manifolds of dimension five and six, obtained by Giaro, Kumpera and Ruiz (see
\cite{giaro-kumpera-ruiz} and \cite{kumpera-ruiz}).

(i) Any Goursat structure on a five-manifold is locally equivalent to one of
the two following Kumpera-Ruiz normal forms%

\begin{align}
& \left(
\begin{array}
[c]{c}%
\tfrac\partial{\partial x_{5}}%
\end{array}
,
\begin{array}
[c]{c}%
x_{5}\tfrac\partial{\partial x_{4}}+x_{4}\tfrac\partial{\partial x_{3}}%
+x_{3}\tfrac\partial{\partial x_{2}}+\tfrac\partial{\partial x_{1}}%
\end{array}
\right) \label{KR-normal-forms-R5a}\\
& \left(
\begin{array}
[c]{c}%
\tfrac\partial{\partial x_{5}}%
\end{array}
,
\begin{array}
[c]{c}%
\tfrac\partial{\partial x_{4}}+x_{5}\left(  x_{4}\tfrac\partial{\partial
x_{3}}+x_{3}\tfrac\partial{\partial x_{2}}+\tfrac\partial{\partial x_{1}%
}\right)
\end{array}
\right)  ,\label{KR-normal-forms-R5b}%
\end{align}
which are not locally equivalent at zero.

(ii) Any Goursat structure on a six-manifold is locally equivalent to one of
the five following Kumpera-Ruiz normal forms
\begin{align}
& \ \left(
\begin{array}
[c]{c}%
\tfrac\partial{\partial x_{6}}%
\end{array}
,
\begin{array}
[c]{c}%
x_{6}\tfrac\partial{\partial x_{5}}+x_{5}\tfrac\partial{\partial x_{4}}%
+x_{4}\tfrac\partial{\partial x_{3}}+x_{3}\tfrac\partial{\partial x_{2}%
}+\tfrac\partial{\partial x_{1}}%
\end{array}
\right) \label{KR-normal-forms-R6a}\\
& \ \left(
\begin{array}
[c]{c}%
\tfrac\partial{\partial x_{6}}%
\end{array}
,
\begin{array}
[c]{c}%
\tfrac\partial{\partial x_{5}}+x_{6}\left(  x_{5}\tfrac\partial{\partial
x_{4}}+x_{4}\tfrac\partial{\partial x_{3}}+x_{3}\tfrac\partial{\partial x_{2}%
}+\tfrac\partial{\partial x_{1}}\right)
\end{array}
\right) \label{KR-normal-forms-R6b}\\
& \ \left(
\begin{array}
[c]{c}%
\tfrac\partial{\partial x_{6}}%
\end{array}
,
\begin{array}
[c]{c}%
x_{6}\tfrac\partial{\partial x_{5}}+\tfrac\partial{\partial x_{4}}%
+x_{5}\left(  x_{4}\tfrac\partial{\partial x_{3}}+x_{3}\tfrac\partial{\partial
x_{2}}+\tfrac\partial{\partial x_{1}}\right)
\end{array}
\right) \label{KR-normal-forms-R6c}\\
& \ \left(
\begin{array}
[c]{c}%
\tfrac\partial{\partial x_{6}}%
\end{array}
,
\begin{array}
[c]{c}%
(x_{6}+1)\tfrac\partial{\partial x_{5}}+\tfrac\partial{\partial x_{4}}%
+x_{5}\left(  x_{4}\tfrac\partial{\partial x_{3}}+x_{3}\tfrac\partial{\partial
x_{2}}+\tfrac\partial{\partial x_{1}}\right)
\end{array}
\right) \label{KR-normal-forms-R6d}\\
& \ \left(
\begin{array}
[c]{c}%
\tfrac\partial{\partial x_{6}}%
\end{array}
,
\begin{array}
[c]{c}%
\tfrac\partial{\partial x_{5}}+x_{6}\left(  \tfrac\partial{\partial x_{4}%
}+x_{5}\left(  x_{4}\tfrac\partial{\partial x_{3}}+x_{3}\tfrac\partial
{\partial x_{2}}+\tfrac\partial{\partial x_{1}}\right)  \right)
\end{array}
\right)  ,\label{KR-normal-forms-R6e}%
\end{align}
which are pairwise locally non-equivalent at zero. Observe that these two
results are not implied by Theorem \ref{thm-kumpera-ruiz}. Indeed, they show
that in dimension five and six the constants that appear in Kumpera-Ruiz's
Theorem can be ``normalized'' to either $0$ or $1$.

For $n=7$, $8$ and $9$ the complete classification is more delicate, but there
is still a finite number of models (see \cite{cheaito-mormul}, \cite{gaspar},
\cite{kumpera-ruiz}, and \cite{mormul-R9}). For $n\geq10$, the number of local
models is infinite (see \cite{cheaito-mormul-pasillas-respondek},
\cite{mormul-R9}, and Section \ref{sec-contact}) and the complete
classification remains an open problem (see recent results in
\cite{mormul-dijon}).

\subsection{Proof of Kumpera-Ruiz's Theorem}

Our proof of Theorem \ref{thm-kumpera-ruiz} will be based on the following
classical result, which was originally formulated in the dual language of
Pfaffian systems \cite[Theorem V]{weber-article} (see also
\cite{cartan-equivalence-absolue}, \cite{goursat}, \cite{kumpera-ruiz}, and
Appendix~\ref{sec-weber}).

\begin{proposition}
[E. von Weber]\label{prop-extended-engel}Let $\mathcal{D}$ be a Goursat
structure on a manifold $M$ of dimension $n\geq4$. In a small enough
neighborhood of any point $p$ in $M$, the distribution $\mathcal{D}$ is
equivalent to a distribution spanned on $\mathbb{R}^{n}$ by a pair of vector
fields that has the following form:
\begin{equation}
\left(
\begin{array}
[c]{c}%
\tfrac{\partial}{\partial y_{n}}%
\end{array}
,
\begin{array}
[c]{c}%
y_{n}\zeta_{1}^{n-1}+\zeta_{2}^{n-1}%
\end{array}
\right) \label{extended-engel-normal-form}%
\end{equation}
where $\zeta_{1}^{n-1}$ and $\zeta_{2}^{n-1}$ are the lifts of a pair of
vector fields that span a Goursat structure on$~\mathbb{R}^{n-1}$ and the
coordinates $y_{1},\ldots,y_{n}$ are centered at~$p$.\ 
\end{proposition}

\noindent\textbf{Proof of Proposition \ref{prop-extended-engel}} It is well
known (see e.g. \cite{bryant-hsu}, \cite{sussmann-liu}, and
\cite{zhitomirskii-nice}) that any Goursat structure $\mathcal{D}$ on a
manifold of dimension $n\geq4$ admits a canonical line field $\mathcal{L}%
\subset\mathcal{D}$ uniquely defined by $[\mathcal{L},\mathcal{D}%
^{(1)}]\subset\mathcal{D}^{(1)}$. Observe that in the preliminary normal form
(\ref{extended-engel-normal-form}) of Proposition~\ref{prop-extended-engel}
this line field is given by $\mathcal{L}=(\tfrac\partial{\partial y_{n}})$.

It is clear that, applying around $p$ the flow-box theorem to a vector field
that spans $\mathcal{L}$, we can chose local coordinates $(z_{1},\ldots
,z_{n})$, centered at $p$, such that $\mathcal{D}$ is locally equivalent to a
distribution spanned on $\mathbb{R}^{n}$ by a pair of vector fields that has
the following form:
\[
\left(
\begin{array}
[c]{c}%
\tfrac\partial{\partial z_{n}}%
\end{array}
,
\begin{array}
[c]{c}%
{\textstyle\sum_{i=2}^{n-1}}
\alpha_{i}(z)\tfrac\partial{\partial z_{i}}+\tfrac\partial{\partial z_{1}}%
\end{array}
\right)  ,
\]
where $\mathcal{L}=(\tfrac\partial{\partial z_{n}})$. Since $\dim
\mathcal{D}^{(1)}(p)=3$ there exists an integer $i$ such that $\frac
{\partial\alpha_{i}}{\partial z_{n}}(p)\neq0$. We can assume that $i=n-1$ and,
moreover, that $\alpha_{n-1}(0)=0$. Otherwise, replace the coordinate
$z_{n-1}$ by $z_{n-1}-z_{1}\alpha_{n-1}(0)$. Now, if we define $y_{n}%
=\alpha_{n-1}(z)$ and $y_{i}=z_{i}$, for $1\leq i\leq n-1$, we get that
$\mathcal{D}$ is locally equivalent to a distribution spanned on
$\mathbb{R}^{n}$ by a pair of vector fields that has the following form:
\[
\left(
\begin{array}
[c]{c}%
\tfrac\partial{\partial y_{n}}%
\end{array}
,
\begin{array}
[c]{c}%
\;y_{n}\tfrac\partial{\partial y_{n-1}}+%
{\textstyle\sum_{i=2}^{n-2}}
\beta_{i}(y)\tfrac\partial{\partial y_{i}}+\tfrac\partial{\partial y_{1}}%
\end{array}
\right)  .
\]
But the inclusion $[\mathcal{L},\mathcal{D}^{(1)}]\subset\mathcal{D}^{(1)}$
clearly implies $\frac{\partial^{2}\beta_{i}}{\partial y_{n}^{2}}\equiv0$ for
$2\leq i\leq n-2$. That is $\beta_{i}(y)=a_{i}(\overline{y}_{n-1})y_{n}%
+b_{i}(\overline{y}_{n-1})$, where $\overline{y}_{n-1}=(y_{1},\ldots,y_{n-1}%
)$. Define
\[
\zeta_{1}^{n-1}=\tfrac\partial{\partial y_{n-1}}+%
{\textstyle\sum_{i=2}^{n-2}}
a_{i}(\overline{y}_{n-1})\tfrac\partial{\partial y_{i}}\text{\quad and\quad
}\zeta_{2}^{n-2}=%
{\textstyle\sum_{i=2}^{n-2}}
b_{i}(\overline{y}_{n-1})\tfrac\partial{\partial y_{i}}+\tfrac\partial
{\partial y_{1}}.
\]
We conclude that $\mathcal{D}$ is equivalent to $(\;\tfrac\partial{\partial
y_{n}}\;,\;y_{n}\zeta_{1}^{n-1}+\zeta_{2}^{n-1}\;)$, where both $\zeta
_{1}^{n-1}$ and $\zeta_{2}^{n-1}$ are lifts (see Notation~\ref{not-lift}) of
vector fields defined on $\mathbb{R}^{n-1}$. Put $\mathcal{F}=(\zeta_{1}%
^{n-1},\zeta_{2}^{n-1})$. Clearly $\dim\mathcal{D}^{(i+1)}=\dim\mathcal{F}%
^{(i)}+1$, for $0\leq i\leq n-3$. It follows that the
distribution$~\mathcal{F}$ is a Goursat structure on$~\mathbb{R}^{n-1}$%
.\hfill$\square$

\medskip\ 

\noindent\textbf{Proof of Theorem \ref{thm-kumpera-ruiz}} On three-manifolds,
Theorem~\ref{thm-kumpera-ruiz} is obviously true. Indeed, it is the solution
of the Pfaff problem (see \cite{darboux} and \cite{frobenius}; see
also~\cite{bryant-chern-gardner-goldschmidt-griffiths} for a modern approach).
We will proceed by induction on $n\geq4$, showing that if any Goursat
structure on an $(n-1)$-manifold can be converted locally into a Kumpera-Ruiz
normal form then the same is true for any Goursat structure on a manifold of
dimension~$n$.

Let $\mathcal{D}$ be a Goursat structure on a manifold $M$ of dimension
$n\geq4$ and let $p$ be an arbitrary point in$~M$. It follows from Proposition
\ref{prop-extended-engel} that $\mathcal{D}$ is equivalent, in a small enough
neighborhood of $p$, to a distribution spanned on $\mathbb{R}^{n}$ by a pair
of vector fields~$(\zeta_{1}^{n},\zeta_{2}^{n})$ that has the following form:
\[%
\begin{array}
[c]{l}%
\zeta_{1}^{n}=\tfrac\partial{\partial y_{n}}\\
\zeta_{2}^{n}=y_{n}\zeta_{1}^{n-1}+\zeta_{2}^{n-1}.
\end{array}
\]
In the rest of the proof we will assume that $\mathcal{D}=(\zeta_{1}^{n}%
,\zeta_{2}^{n})$. Note that the $y$-coordinates are centered at zero.

Our aim is to build a local change of coordinates
\[
(x_{1},\ldots,x_{n})=\phi^{n}(y_{1},\ldots,y_{n}),
\]
a Kumpera-Ruiz normal form $(\kappa_{1}^{n},\kappa_{2}^{n})$ on $\mathbb{R}%
^{n}$, and three smooth functions $\mu_{n}$, $\nu_{n}$, and $\eta_{n}$ such
that
\begin{equation}%
\begin{array}
[c]{l}%
\phi_{*}^{n}(\zeta_{1}^{n})=(\nu_{n}\circ\psi^{n})\kappa_{1}^{n}\\
\phi_{*}^{n}(\zeta_{2}^{n})=(\eta_{n}\circ\psi^{n})\kappa_{1}^{n}+(\mu
_{n}\circ\psi^{n})\kappa_{2}^{n},
\end{array}
\label{proof-kumpera-ruiz-01}%
\end{equation}
where $\psi^{n}=(\phi^{n})^{-1}$ denotes the inverse of the local
diffeomorphism $\phi^{n}$ and both $\mu_{n}(0)\neq0$ and $\nu_{n}(0)\neq0$.
Moreover, we will impose the $x$-coordinates to be centered at zero. That is
$\phi^{n}(0)=0$. Observe that the triangular form in
(\ref{proof-kumpera-ruiz-01}) appears because both $\zeta_{1}^{n}$ and
$\kappa_{1}^{n}$ span the canonical line fields of the distributions spanned
by $(\zeta_{1}^{n},\zeta_{2}^{n})$ and $(\kappa_{1}^{n},\kappa_{2}^{n})$, respectively.

By Proposition \ref{prop-extended-engel}, the distribution spanned by
$(\zeta_{1}^{n},\zeta_{2}^{n})$ is defined by the lifts of a pair of vector
fields $(\zeta_{1}^{n-1},\zeta_{2}^{n-1})$ that span a Goursat structure on
$\mathbb{R}^{n-1}$. Since the Theorem is assumed to be true on $\mathbb{R}%
^{n-1}$, the distribution spanned by $(\zeta_{1}^{n-1},\zeta_{2}^{n-1})$ is
locally equivalent to a Kumpera-Ruiz normal form $(\kappa_{1}^{n-1},\kappa
_{2}^{n-1})$ defined on $\mathbb{R}^{n-1}$ and centered at zero. It follows
that there exists a local diffeomorphism $(x_{1},\ldots,x_{n-1})=\phi
^{n-1}(y_{1},\ldots,y_{n-1})$ and four smooth functions $\nu_{n-1}$,
$\lambda_{n-1}$, $\eta_{n-1}$, and $\mu_{n-1}$ such that:
\begin{equation}%
\begin{array}
[c]{l}%
\phi_{*}^{n-1}(\zeta_{1}^{n-1})=(\nu_{n-1}\circ\psi^{n-1})\kappa_{1}%
^{n-1}+(\lambda_{n-1}\circ\psi^{n-1})\kappa_{2}^{n-1}\\
\phi_{*}^{n-1}(\zeta_{2}^{n-1})=(\eta_{n-1}\circ\psi^{n-1})\kappa_{1}%
^{n-1}+(\mu_{n-1}\circ\psi^{n-1})\kappa_{2}^{n-1},
\end{array}
\label{proof-kumpera-ruiz-02}%
\end{equation}
where $\psi^{n-1}=(\phi^{n-1})^{-1}$ denotes the inverse of the local
diffeomorphism $\phi^{n-1}$ and $(\nu_{n-1}\mu_{n-1}-\lambda_{n-1}\eta
_{n-1})(0)\neq0$.

Let $\phi^{n}=(\phi^{n-1},\phi_{n})^{T}$ be a diffeomorphism of $\mathbb{R}%
^{n}$ such that $\phi^{n-1}$ depends on the first $n-1$ coordinates only.
Moreover, let $f$ be a vector field on $\mathbb{R}^{n}$ of the form $f=\alpha
f^{n-1}+f_{n}$, where $\alpha$ is a smooth function on $\mathbb{R}^{n}$, the
vector field $f^{n-1}$ is the lift of a vector field on $\mathbb{R}^{n-1}$
(see Notation~\ref{not-lift}), and the only non-zero component of $f_{n}$ is
the last one. A direct computation shows that:
\begin{equation}
\phi_{*}^{n}(f)=(\alpha\circ\psi^{n})\phi_{*}^{n-1}(f^{n-1})+\left(
(\mathrm{L}_{f}\phi_{n})\circ\psi^{n}\right)  \tfrac\partial{\partial x_{n}%
}.\label{triangular-trick}%
\end{equation}
Note that the vector field $\phi_{*}^{n-1}(f^{n-1})$ is lifted along the
$x_{n}$-coordinate, which is defined by~$\phi_{n}$.

\medskip\ 

\noindent\emph{Regular case:} If $\mu_{n-1}(0)\neq0$ then we can complete
$\phi^{n-1}$ to a zero-preserving local diffeomorphism of~$\mathbb{R}^{n}$ by
taking $\phi^{n}=(\phi^{n-1},\phi_{n})^{T}$, where
\[
\phi_{n}(y)=\dfrac{\nu_{n-1}y_{n}+\eta_{n-1}}{\lambda_{n-1}y_{n}+\mu_{n-1}%
}-\frac{\eta_{n-1}(0)}{\mu_{n-1}(0)}.
\]
In this case, we define $c_{n}=(\eta_{n-1}/\mu_{n-1})(0)$ and
\[
\nu_{n}=\mathrm{L}_{\zeta_{1}^{n}}\phi_{n},\text{\quad}\eta_{n}=\mathrm{L}%
_{\zeta_{2}^{n}}\phi_{n},\text{\quad and\quad}\mu_{n}=\lambda_{n-1}y_{n}%
+\mu_{n-1}.
\]
Observe that $\nu_{n}(0)=\mathrm{L}_{\zeta_{1}^{n}}\phi_{n}(0)=(\nu_{n-1}%
\mu_{n-1}-\lambda_{n-1}\eta_{n-1})(0)\neq0$ and that $\mu_{n}(0)=\mu
_{n-1}(0)\neq0$. Thus the right hand side of~(\ref{proof-kumpera-ruiz-01})
defines a locally invertible transformation. Moreover, the Kumpera-Ruiz normal
form $(\kappa_{1}^{n},\kappa_{2}^{n})$ is defined to be the regular
prolongation, with parameter $c_{n}$, of $(\kappa_{1}^{n-1},\kappa_{2}^{n-1})$.

Let us check that, in this case, relation (\ref{proof-kumpera-ruiz-01}) holds.
Together, relations~(\ref{proof-kumpera-ruiz-02}) and~(\ref{triangular-trick})
give:
\begin{align*}
\phi_{*}^{n}(\zeta_{2}^{n})  & =(y_{n}\circ\psi^{n})\phi_{*}^{n-1}(\zeta
_{1}^{n-1})+\phi_{*}^{n-1}(\zeta_{2}^{n-1})+\left(  (\mathrm{L}_{\zeta_{2}%
^{n}}\phi_{n})\circ\psi^{n}\right)  \tfrac\partial{\partial x_{n}}\\
\  & =\left(  \left(  \nu_{n-1}y_{n}+\eta_{n-1}\right)  \circ\psi^{n}\right)
\kappa_{1}^{n-1}+\left(  \left(  \lambda_{n-1}y_{n}+\mu_{n-1}\right)
\circ\psi^{n}\right)  \kappa_{2}^{n-1}+(\eta_{n}\circ\psi^{n})\kappa_{1}^{n}\\
& =\left(  \left(  \lambda_{n-1}y_{n}+\mu_{n-1}\right)  \circ\psi^{n}\right)
\left(  \left(  \tfrac{\nu_{n-1}y_{n}+\eta_{n-1}}{\lambda_{n-1}y_{n}+\mu
_{n-1}}\circ\psi^{n}\right)  \kappa_{1}^{n-1}+\kappa_{2}^{n-1}\right)
+(\eta_{n}\circ\psi^{n})\kappa_{1}^{n}\\
\  & =\left(  \mu_{n}\circ\psi^{n}\right)  \left(  \left(  x_{n}+c_{n}\right)
\kappa_{1}^{n-1}+\kappa_{2}^{n-1}\right)  +(\eta_{n}\circ\psi^{n})\kappa
_{1}^{n}\\
\  & =(\eta_{n}\circ\psi^{n})\kappa_{1}^{n}+(\mu_{n}\circ\psi^{n})\kappa
_{2}^{n}.
\end{align*}
Moreover, we have
\[
\phi_{*}^{n}(\zeta_{1}^{n})=\left(  (\mathrm{L}_{\zeta_{1}^{n}}\phi_{n}%
)\circ\psi^{n}\right)  \tfrac\partial{\partial x_{n}}=(\nu_{n}\circ\psi
^{n})\kappa_{1}^{n+3}.
\]
It follows that, in the regular case, relation (\ref{proof-kumpera-ruiz-01}) holds.

\medskip\ 

\noindent\emph{Singular case:} If $\mu_{n-1}(0)=0$ then we can complete
$\phi^{n-1}$ to a zero-preserving local diffeomorphism of $\mathbb{R}^{n}$ by
taking $\phi^{n}=(\phi^{n-1},\phi_{n})^{T}$, where
\[
\phi_{n}(y)=\dfrac{\lambda_{n-1}y_{n}+\mu_{n-1}}{\nu_{n-1}y_{n}+\eta_{n-1}}.
\]
Observe that $\mu_{n-1}(0)=0$ implies $\phi_{n}(0)=0$. Additionally, since
$\mu_{n-1}(0)=0$ and $(\nu_{n-1}\mu_{n-1}-\lambda_{n-1}\eta_{n-1})(0)\neq0$,
we have $\lambda_{n-1}(0)\neq0$ and $\eta_{n-1}(0)\neq0$. In this case, we
define
\[
\nu_{n}=\mathrm{L}_{\zeta_{1}^{n}}\phi_{n},\text{\quad}\eta_{n}=\mathrm{L}%
_{\zeta_{2}^{n}}\phi_{n},\text{\quad and\quad}\mu_{n}=\nu_{n-1}y_{n}%
+\eta_{n-1}.
\]
Observe that $\nu_{n}(0)=\mathrm{L}_{\zeta_{1}^{n}}\phi_{n}(0)=(\nu_{n-1}%
\mu_{n-1}-\lambda_{n-1}\eta_{n-1})(0)\neq0$ and that $\mu_{n}(0)=\eta
_{n-1}(0)\neq0$. Thus the right hand side of~(\ref{proof-kumpera-ruiz-01})
defines a locally invertible transformation. Moreover, the Kumpera-Ruiz normal
form $(\kappa_{1}^{n},\kappa_{2}^{n})$ is defined to be the singular
prolongation of $(\kappa_{1}^{n-1},\kappa_{2}^{n-1})$.

Let us check that, again, relation (\ref{proof-kumpera-ruiz-01}) holds.
Together, relations~(\ref{proof-kumpera-ruiz-02}) and~(\ref{triangular-trick})
give:
\begin{align*}
\phi_{*}^{n}(\zeta_{2}^{n})  & =(y_{n}\circ\psi^{n})\phi_{*}^{n-1}(\zeta
_{1}^{n-1})+\phi_{*}^{n-1}(\zeta_{2}^{n-1})\ +\left(  (\mathrm{L}_{\zeta
_{2}^{n}}\phi_{n})\circ\psi^{n}\right)  \tfrac\partial{\partial x_{n}}\\
\  & =\left(  \left(  \nu_{n-1}y_{n}+\eta_{n-1}\right)  \circ\psi^{n}\right)
\kappa_{1}^{n-1}+\left(  \left(  \lambda_{n-1}y_{n}+\mu_{n-1}\right)
\circ\psi^{n}\right)  \kappa_{2}^{n-1}\ +(\eta_{n}\circ\psi^{n})\kappa_{1}%
^{n}\\
\  & =\left(  \left(  \nu_{n-1}y_{n}+\eta_{n-1}\right)  \circ\psi^{n}\right)
\left(  \kappa_{1}^{n-1}+\left(  \tfrac{\lambda_{n-1}y_{n}+\mu_{n-1}}%
{\nu_{n-1}y_{n}+\eta_{n-1}}\circ\psi^{n}\right)  \kappa_{2}^{n-1}\right)
+(\eta_{n}\circ\psi^{n})\kappa_{1}^{n}\\
\  & =\left(  \mu_{n}\circ\psi^{n}\right)  \left(  \kappa_{1}^{n-1}%
+x_{n}\kappa_{2}^{n-1}\right)  +(\eta_{n}\circ\psi^{n})\kappa_{1}^{n}\\
\  & =(\eta_{n}\circ\psi^{n})\kappa_{1}^{n}+(\mu_{n}\circ\psi^{n})\kappa
_{2}^{n}.
\end{align*}
Like in the previous case, we have
\[
\phi_{*}^{n}(\zeta_{1}^{n})=\left(  (\mathrm{L}_{\zeta_{1}^{n}}\phi_{n}%
)\circ\psi^{n}\right)  \tfrac\partial{\partial x_{n}}=(\nu_{n}\circ\psi
^{n})\kappa_{1}^{n+3}.
\]
It follows that relation (\ref{proof-kumpera-ruiz-01}) holds in both
cases.\hfill$\square$

\section{The N-Trailer System}

\label{sec-trailer}

\subsection{Definition of the N-Trailer System}

The kinematical model for a unicycle-like mobile robot towing $n$ trailers
such that the tow hook of each trailer is located at the center of its unique
axle is usually called, in control theory, the $n$-trailer system (see the
books~\cite{laumond-book} and~\cite{li-canny-book}; the papers
\cite{fliess-intro-flat}, \cite{jakubczyk-trailer}, \cite{jean-trailer},
\cite{laumond-trailer}, \cite{fliess-trailer}, \cite{samson-chained},
\cite{sordalen-trailer}, \cite{teel-murray-walsh}, and
\cite{tilbury-murray-sastry}; and the references given there). Figures
representing this system are given in Appendix~\ref{sec-figures}. For
simplicity, we will assume that the distances between any two consecutive
trailers are equal.

\begin{definition}
\label{def-trailer}The $n$\emph{-trailer system} is the distribution defined
on $\mathbb{R}^{2}\times(S^{1})^{n+1}$, for $n\geq0$, by the following pair of
vector fields:
\begin{equation}
\left(
\begin{array}
[c]{c}%
\begin{array}
[c]{c}%
\tfrac{\partial}{\partial\theta_{n}}%
\end{array}
,
\begin{array}
[c]{c}%
\pi_{0}\cos(\theta_{0})\tfrac{\partial}{\partial\xi_{1}}+\pi_{0}\sin
(\theta_{0})\tfrac{\partial}{\partial\xi_{2}}+%
{\textstyle\sum\limits_{i=0}^{n-1}}
\pi_{i+1}\sin(\theta_{i+1}-\theta_{i})\tfrac{\partial}{\partial\theta_{i}}%
\end{array}
\end{array}
\right)  ,\label{trailer-system}%
\end{equation}
where $\pi_{i}=%
{\textstyle\prod\nolimits_{j=i+1}^{n}}
\cos(\theta_{j}-\theta_{j-1})$ and $\pi_{n}=1$.
\end{definition}

In the above definition, the functions $\xi_{1}$, $\xi_{2}$, $\theta_{0}$,...,
$\theta_{n}$ are coordinates on the manifold $\mathbb{R}^{2}\times
(S^{1})^{n+1}$. The coordinates $\xi_{1}$ and $\xi_{2}$ represent the position
of the last trailer, while the coordinates $\theta_{0},\ldots,\theta_{n}$
represent, starting from the last trailer, the angles between each trailer's
axle and the $\xi_{1}$-axis. It is easy to check that the $n$-trailer system
is a Goursat structure (see e.g.~\cite{laumond-singularities}).

We give now an equivalent inductive definition of the $n$-trailer. This
definition already appears in~\cite{jean-trailer} and reminds the one given in
the previous section for Kumpera-Ruiz normal forms. To start with, consider
the pair of vector fields $(\tau_{1}^{0},\tau_{2}^{0})$ on $\mathbb{R}%
^{2}\times S^{1}$ that describe the kinematics of the unicycle-like mobile
robot towing no trailers:
\[%
\begin{array}
[c]{l}%
\tau_{1}^{0}=\tfrac{\partial}{\partial\theta_{0}}\\
\tau_{2}^{0}=\cos(\theta_{0})\tfrac{\partial}{\partial\xi_{1}}+\sin(\theta
_{0})\tfrac{\partial}{\partial\xi_{2}}.
\end{array}
\]
The $n$-trailer system can be defined by applying successively a sequence of
prolongations to this mobile robot. In order to do this, suppose that a pair
of vector fields $\tau^{n-1}=(\tau_{1}^{n-1},\tau_{2}^{n-1})$ on
$\mathbb{R}^{2}\times(S^{1})^{n}$, associated to the mobile robot towing $n-1$
trailers, has been defined. The pair of vector fields $\tau^{n}=(\tau_{1}%
^{n},\tau_{2}^{n})$ on $\mathbb{R}^{2}\times(S^{1})^{n+1}$ corresponding to
the $n$-trailer system is then given by
\[%
\begin{array}
[c]{l}%
\tau_{1}^{n}=\tfrac{\partial}{\partial\theta_{n}}\\
\tau_{2}^{n}=\sin(\theta_{n}-\theta_{n-1})\tau_{1}^{n-1}+\cos(\theta
_{n}-\theta_{n-1})\tau_{2}^{n-1}.
\end{array}
\]
Observe that this definition should be understood in the sense of
Notation~\ref{not-lift}, and that the pair of vector fields $(\tau_{1}%
^{n},\tau_{2}^{n})$ coincides with that of Definition~\ref{def-trailer}.
Mechanically, to prolongate the $n$-trailer means to add one more trailer to
the system.

\subsection{Local Conversion of the N-Trailer System into a Kum\-pe\-ra-Ruiz Normal Form}

Since the $n$-trailer is a Goursat structure, it follows directly from
Kumpera-Ruiz's theorem that, in a small enough neighborhood of any point of
its configuration space, in particular at any singular configuration, the
$n$-trailer can be converted into Kumpera-Ruiz normal form. In this subsection
we describe this conversion explicitly. For regular configurations, our result
gives the transformations proposed in~\cite{sordalen-trailer} and~
\cite{tilbury-murray-sastry}; for singular configurations, our result gives a
new kind of transformations.

Denote by $\zeta$ the coordinates of $\mathbb{R}^{2}\times(S^{1})^{n+1}$, that
is
\[
\zeta=(\zeta_{1},...,\zeta_{n+3})=(\xi_{1},\xi_{2},\theta_{0},...,\theta
_{n}).
\]
Fix a point $p$ of $\mathbb{R}^{2}\times(S^{1})^{n+1}$ given in $\zeta
$-coordinates by $\zeta(p)=\zeta^{p}=(\xi_{1}^{p},\xi_{2}^{p},\theta_{0}%
^{p},...,\theta_{n}^{p})$. In order to convert, locally at~$p$, the
$n$-trailer into a Kumpera-Ruiz normal form we look for a local change of
coordinates
\[
(x_{1},\ldots,x_{n+3})=\phi^{n}(\xi_{1},\xi_{2},\theta_{0},\ldots,\theta
_{n}),
\]
a Kumpera-Ruiz normal form $(\kappa_{1}^{n+3},\kappa_{2}^{n+3})$ on
$\mathbb{R}^{n+3}$, and three smooth functions $\nu_{n}$, $\eta_{n}$, and
$\mu_{n}$ such that
\begin{equation}%
\begin{array}
[c]{lll}%
\phi_{*}^{n}(\tau_{1}^{n}) & = & (\nu_{n}\circ\psi^{n})\,\kappa_{1}^{n+3}\\
\phi_{*}^{n}(\tau_{2}^{n}) & = & (\eta_{n}\circ\psi^{n})\,\kappa_{1}%
^{n+3}+(\mu_{n}\circ\psi^{n})\,\kappa_{2}^{n+3},
\end{array}
\label{trailer-as-kumpera}%
\end{equation}
where $\psi^{n}=(\phi^{n})^{-1}$ denotes the inverse of the local
diffeomorphism $\phi^{n}$ and both $\nu_{n}(\zeta^{p})\neq0$ and $\mu
_{n}(\zeta^{p})\neq0$. Observe that we do not demand the $x$-coordinates to be
centered at $p$, and thus the point $x(p)=(\phi^{n}\circ\zeta)(p)$ will be, in
general, different from zero.

We construct inductively here a change of coordinates $\phi^{n}=(\phi
_{1},\ldots,\phi_{n+3})^{T}$ and three smooth functions $\nu_{n}$, $\eta_{n}$,
and $\mu_{n}$ that satisfy (\ref{trailer-as-kumpera}). We start with $n=0$. If
$\theta_{0}^{p}\neq\pm\pi/2\operatorname*{mod}2\pi$ then define $\phi_{1}%
=\xi_{1}$, $\phi_{2}=\xi_{2}$, and $\phi_{3}=\tan(\theta_{0})$. Moreover, take
$\mu_{0}=\cos(\theta_{0})$, $\nu_{0}=\sec^{2}(\theta_{0})$, and $\eta_{0}=0$.
If $\theta_{0}^{p}=\pm\pi/2\operatorname*{mod}2\pi$ then define $\phi_{1}%
=\xi_{2}$, $\phi_{2}=\xi_{1}$, and $\phi_{3}=\cot(\theta_{0})$.\ Moreover,
take $\mu_{0}=\sin(\theta_{0})$, $\nu_{0}=-\csc^{2}(\theta_{0})$, and
$\eta_{0}=0$. Denote $s_{i}=\sin(\theta_{i}-\theta_{i-1})$ and $c_{i}%
=\cos(\theta_{i}-\theta_{i-1})$, for $0\leq i\leq n$.

Now, consider the sequence of smooth functions defined locally around the
point$~\zeta(p)$, for $1\leq i\leq n$, by either
\begin{equation}%
\begin{array}
[c]{lll}%
\phi_{i+3} & = & \dfrac{s_{i}\nu_{i-1}+c_{i}\eta_{i-1}}{c_{i}\mu_{i-1}}\\
&  & \\
\mu_{i} & = & c_{i}\mu_{i-1}\qquad\qquad\\
\nu_{i} & = & \mathrm{L}_{\tau_{1}^{i}}\phi_{i+3}\\
\eta_{i} & = & \mathrm{L}_{\tau_{2}^{i}}\phi_{i+3},
\end{array}
\quad\label{regular-trailer-transformation}%
\end{equation}
if $\theta_{i}^{p}-\theta_{i-1}^{p}\neq\pm\pi/2\operatorname*{mod}2\pi$
(\emph{regular case}) or by
\begin{equation}%
\begin{array}
[c]{lcl}%
\phi_{i+3} & = & \dfrac{c_{i}\mu_{i-1}}{s_{i}\nu_{i-1}+c_{i}\eta_{i-1}}\\
&  & \\
\mu_{i} & = & s_{i}\nu_{i-1}+c_{i}\eta_{i-1}\\
\nu_{i} & = & \mathrm{L}_{\tau_{1}^{i}}\phi_{i+3}\\
\eta_{i} & = & \mathrm{L}_{\tau_{2}^{i}}\phi_{i+3},
\end{array}
\quad\label{singular-trailer-transformation}%
\end{equation}
if $\theta_{i}^{p}-\theta_{i-1}^{p}=\pm\pi/2\operatorname*{mod}2\pi$
(\emph{singular case}). It is easy to prove that, for $0\leq i\leq n$, the
transformations defined by $(x_{1},\ldots,x_{i+3})=\phi^{i}(\xi_{1},\xi
_{2},\theta_{0},\ldots,\theta_{i})$ are smooth changes of coordinates around
$p_{i}$ and that, moreover, we have both $\nu_{i}(\zeta^{p_{i}})\neq0$ and
$\mu_{i}(\zeta^{p_{i}})\neq0$, where~$p_{i}$ denotes the projection of $p$ on
$\mathbb{R}^{2}\times(S^{1})^{i+1}$, the product of $\mathbb{R}^{2}$ and the
first $i+1$ copies of $S^{1}$, that is $\zeta^{p_{i}}=(\xi_{1}^{p},\xi_{2}%
^{p},\theta_{0}^{p},\ldots,\theta_{i}^{p})$.

\begin{proposition}
\label{prop-trailer-conversion}For $n\geq0$, the local diffeomorphism
$\phi^{n}$ and the smooth functions $\nu_{n}$, $\eta_{n}$, and~$\mu_{n}$
satisfy \emph{(\ref{trailer-as-kumpera})}, and thus convert locally the
$n$-trailer system into a Kumpera-Ruiz normal form.
\end{proposition}

\noindent\textbf{Proof of Proposition \ref{prop-trailer-conversion}}{\ }We
will prove that the relation~(\ref{trailer-as-kumpera}) holds for $n\geq0$ by
induction on the number~$n$ of trailers. Relation~(\ref{trailer-as-kumpera})
is clearly true for $n=0$. Assume that it holds for $n-1$ trailers, that~is
\begin{align*}
\phi_{*}^{n-1}(\tau_{1}^{n-1})  & =(\nu_{n-1}\circ\psi^{n-1})\,\kappa
_{1}^{n+2}\\
\phi_{*}^{n-1}(\tau_{2}^{n-1})  & =(\eta_{n-1}\circ\psi^{n-1})\,\kappa
_{1}^{n+2}+(\mu_{n-1}\circ\psi^{n-1})\,\kappa_{2}^{n+2}.
\end{align*}
The inductive definition of the $n$-trailer gives
\begin{align*}
\tau_{1}^{n}  & =\tfrac\partial{\partial\theta_{n}}\\
\tau_{2}^{n}  & =\sin(\theta_{n}-\theta_{n-1})\tau_{1}^{n-1}+\cos(\theta
_{n}-\theta_{n-1})\tau_{2}^{n-1}.
\end{align*}

\noindent Recall (see the Proof of Theorem~\ref{thm-kumpera-ruiz}) that for a diffeomorphism $\phi^{n}=(\phi^{n-1},\phi_{n+3})^{T}$ of~$\mathbb{R}^{n+3}$,
such that~$\phi^{n-1}$ depends on the first $n+2$ coordinates only, and for a
vector field~$f$ on$~\mathbb{R}^{n+3}$ of the form $f=\alpha f^{n-1}+f_{n+3}$,
where~$\alpha$ is a smooth function on~$\mathbb{R}^{n+3}$, the vector
field~$f^{n-1}$ is the lift of a vector field on~$\mathbb{R}^{n+2}$ (see
Notation~\ref{not-lift}), and the only non-zero component of~$f_{n+3}$ is the
last one, we have
\begin{equation}
\phi_{*}^{n}(f)=(\alpha\circ\psi^{n})\phi_{*}^{n-1}(f^{n-1})+\left(
(\mathrm{L}_{f}\phi_{n+3})\circ\psi^{n}\right)  \tfrac\partial{\partial
x_{n+3}}.\label{composed-derivative}%
\end{equation}
Note that the vector field $\phi_{*}^{n-1}(f^{n-1})$ is lifted along the
$x_{n+3}$-coordinate, which is defined by~$\phi_{n+3}$.

In the regular case, that is if $\theta_{i}^{p}-\theta_{i-1}^{p}\neq\pm
\pi/2\operatorname*{mod}2\pi$, we take a regular prolongation and, using
relations~(\ref{regular-trailer-transformation})
and~(\ref{composed-derivative}), we obtain:
\begin{eqnarray*}
\phi_{*}^{n}(\tau_{2}^{n})  & = & (s_{n}\circ\psi^{n})\phi_{*}^{n-1}(\tau
_{1}^{n-1})+(c_{n}\circ\psi^{n})\phi_{*}^{n-1}(\tau_{2}^{n-1})+\left(
(\mathrm{L}_{\tau_{2}^{n}}\phi_{n+3})\circ\psi^{n}\right)  \tfrac
\partial{\partial x_{n+3}}\\
& = & \left(  (s_{n}\nu_{n-1}+c_{n}\eta_{n-1})\circ\psi^{n}\right)  \kappa
_{1}^{n+2}+\left(  (c_{n}\mu_{n-1})\circ\psi^{n}\right)  \kappa_{2}%
^{n+2}+(\eta_{n}\circ\psi^{n})\kappa_{1}^{n+3}\\
& = & \left(  c_{n}\mu_{n-1}\circ\psi^{n}\right)  \left(  \left(  \tfrac{s_{n}%
\nu_{n-1}+c_{n}\eta_{n-1}}{c_{n}\mu_{n-1}}\circ\psi^{n}\right)  \kappa
_{1}^{n+2}+\kappa_{2}^{n+2}\right)  +(\eta_{n}\circ\psi^{n})\kappa_{1}^{n+3}\\
& = & (\mu_{n}\circ\psi^{n})\left(  x_{n+3}\kappa_{1}^{n+2}+\kappa_{2}%
^{n+2}\right)  +(\eta_{n}\circ\psi^{n})\kappa_{1}^{n+3}\\
& = & (\eta_{n}\circ\psi^{n})\kappa_{1}^{n+3}+(\mu_{n}\circ\psi^{n})\kappa
_{2}^{n+3}.
\end{eqnarray*}
In the singular case, that is if $\theta_{i}^{p}-\theta_{i-1}^{p}=\pm
\pi/2\operatorname*{mod}2\pi$, we take the singular prolongation and, using
relations~(\ref{singular-trailer-transformation})
and~(\ref{composed-derivative}), we obtain:
\begin{eqnarray*}
\phi_{*}^{n}(\tau_{2}^{n})  & = &\left(  (s_{n}\nu_{n-1}+c_{n}\eta_{n-1}%
)\circ\psi^{n}\right)  \left(  \kappa_{1}^{n+2}+\left(  \tfrac{c_{n}\mu_{n-1}%
}{s_{n}\nu_{n-1}+c_{n}\eta_{n-1}}\circ\psi^{n}\right)  \kappa_{2}%
^{n+2}\right) \\
& & \mbox{}+(\eta_{n}\circ\psi^{n})\kappa_{1}^{n+3}\\
& = & (\mu_{n}\circ\psi^{n})\left(  \kappa_{1}^{n+2}+x_{n+3}\kappa_{2}%
^{n+2}\right)  +(\eta_{n}\circ\psi^{n})\kappa_{1}^{n+3}\\
& = & (\eta_{n}\circ\psi^{n})\kappa_{1}^{n+3}+(\mu_{n}\circ\psi^{n})\kappa
_{2}^{n+3}.
\end{eqnarray*}
Moreover, in both cases, we have:
\[
\phi_{*}^{n}(\tau_{1}^{n})=\left(  (\mathrm{L}_{\tau_{1}^{n}}\phi_{n+3}%
)\circ\psi^{n}\right)  \tfrac\partial{\partial x_{n+3}}=(\nu_{n}\circ\psi
^{n})\kappa_{1}^{n+3}.
\]
It follows that relation (\ref{trailer-as-kumpera}) holds for any $n\geq0$.
\hfill$\square$

\subsection{Local Conversion of an Arbitrary Goursat Structure into the
N-Trailer System}

Reversing the construction given in the Proof of
Proposition~\ref{prop-trailer-conversion} leads to the following surprising
result (already announced in~\cite{pasillas-respondek-nolcos} and proved
in~\cite{pasillas-respondek-cdc}), which states that the $n$-trailer system is
a universal local model for all Goursat structures. This theorem will play a
fundamental role in this article. Indeed, it will allow us to generalize local
results known for the $n$-trailer, like the formula for the growth vector
obtained by Jean~\cite{jean-trailer}, to all Goursat structures.

\begin{theorem}
\label{thm-trailer-universal}Any Goursat structure on a manifold $M$ of
dimension $n+3$ is locally equivalent, at any point $q$ in $M$, to the
$n$-trailer considered around a suitably chosen point$~p$ of its configuration
space $\mathbb{R}^{2}\times(S^{1})^{n+1}$.
\end{theorem}

\noindent\textbf{Proof of Theorem \ref{thm-trailer-universal}}{\quad} By
Theorem \ref{thm-kumpera-ruiz}, our Goursat structure is, in a small enough
neighborhood of any point $q\ $in $M$, equivalent to a Kumpera-Ruiz normal
form $\kappa^{n+3}$. Denote by $y=(y_{1},\ldots,y_{n+3})$ the coordinates of
$\kappa^{n+3}$ and put $(y_{1}^{q},\ldots,y_{n+3}^{q})=y(q)$.

Recall that, by definition, the pair of vector fields $\kappa^{n+3}$ is given
by a sequence of prolongations $\kappa^{i}=\sigma_{i-3}\circ\cdots\circ
\sigma_{1}(\kappa^{3})$, where $\sigma_{j}$ belongs to $\{R_{c_{j}},S\}$, for
$1\leq j\leq i-3$ and $3\leq i\leq n+3$. We call a coordinate $y_{i}$ such
that $\kappa^{i}=S(\kappa^{i-1})$ a \emph{singular coordinate}, and a
coordinate~$y_{i}$ such that $\kappa^{i}=R_{c}(\kappa^{i-1})$ a \emph{regular
coordinate}. It follows from the proof of Theorem~\ref{thm-kumpera-ruiz} (see
Section~\ref{sec-kumpera-ruiz}) that for all singular coordinates we have
$y_{i}^{q}=0$; but for regular coordinates, the constants $y_{i}^{q}$ can be
arbitrary real numbers.

To prove Theorem~\ref{thm-trailer-universal}, we will define a point $p$ of
$\mathbb{R}^{2}\times(S^{1})^{n+1}$ whose coordinates $\zeta(p)=\zeta^{p}%
=(\xi_{1}^{p},\xi_{2}^{p},\theta_{0}^{p},\ldots,\theta_{n}^{p})$ satisfy
$(x\circ\zeta)(p)=y(q)$, where $x$ and $\zeta$ denote the coordinates used in
the Proof of Proposition \ref{prop-trailer-conversion}. First, put the axle of
the last trailer at $(y_{1}^{q},y_{2}^{q})$, that is $(\xi_{1}^{p},\xi_{2}%
^{p})=(y_{1}^{q},y_{2}^{q})$, and take $\theta_{0}^{p}=\arctan(y_{3}^{q})$.
Compute $x_{3}=\tan(\theta_{0})$, $\mu_{0}=\cos(\theta_{0})$, $\nu_{0}%
=\sec^{2}(\theta_{0})$, and $\eta_{0}=0 $. Then, take for $i=1$ up to $n$, the
following values for the angles $\theta_{i}^{p}\operatorname*{mod}2\pi$. If
the coordinate $y_{i+3}$ is singular then put $\theta_{i}^{p}=\theta_{i-1}%
^{p}+\pi/2$ and compute the coordinate $x_{i+3}$ and the smooth functions
$\mu_{i}$, $\nu_{i}$, and $\eta_{i}$
using~(\ref{singular-trailer-transformation}). If $y_{i+3}$ is regular then
put
\[
\theta_{i}^{p}=\arctan\left(  \dfrac{\mu_{i-1}(\zeta^{p_{i}})y_{i+3}^{q}%
-\eta_{i-1}(\zeta^{p_{i}})}{\nu_{i-1}(\zeta^{p_{i}})}\right)  +\theta
_{i-1}^{p}
\]
and compute the coordinate $x_{i+3}$ and the smooth functions $\mu_{i}$,
$\nu_{i}$, and $\eta_{i}$ using~(\ref{regular-trailer-transformation}). The
result of this construction is that $(x\circ\zeta)(p)=y(q)$. By Proposition
\ref{prop-trailer-conversion}, the coordinates $x\circ\zeta$ convert the
$n$-trailer into a Kumpera-Ruiz normal form. By the above defined
construction, this normal form has the same singular coordinates
as~$\kappa^{n+3}$ and is defined around the same point of $\mathbb{R}^{n+3}$
(if we translate the regular coordinates in order to center them then those
two Kumpera-Ruiz normal forms have the same constants in the regular
prolongations). Hence, the diffeomorphism $\zeta^{-1}\circ x^{-1}\circ y$
gives the claimed equivalence. \hfill$\square$

\section{Singularity Type}

\label{sec-singularity-type}

\subsection{Characteristic Distributions}

A \emph{characteristic vector field} of a distribution $\mathcal{D}$ is a
vector field $f$ that belongs to $\mathcal{D}$ and satisfies $[f,\mathcal{D}%
]\subset\mathcal{D}$. The \emph{characteristic distribution} of a distribution
$\mathcal{D}$ is the module, over the ring of smooth functions, generated by
all its characteristic vector fields. A characteristic distribution need not
be of constant rank; but it follows directly from the Jacobi identity that a
characteristic distribution is always involutive. The main ingredient in the
definition of the singularity type will be the characteristic distributions
$\mathcal{C}_{i}$ defined by the following result, which is apparently due to
Cartan \cite{cartan-equivalence-absolue}, although he did not state it
explicitly in his published works. Its proof can be found in
\cite{kumpera-ruiz} and \cite{martin-rouchon-driftless} (see
also~\cite{canadas-ruiz-length-one}, \cite{kazarian-montgomery-shapiro}, and
Appendix~\ref{sec-weber}), were slightly stronger versions are proved using
the dual language of Pfaffian systems.

\begin{proposition}
[E. Cartan]\label{prop-cartan}Consider a Goursat structure$~\mathcal{D}$
defined on a manifold of dimension$~n\geq4$. Each distribution~$\mathcal{D}%
^{(i)}$, for $0\leq i\leq n-4$, contains a unique involutive
subdistribution~$\mathcal{C}_{i}\subset\mathcal{D}^{(i)}$ that is
characteristic for~$\mathcal{D}^{(i+1)}$ and has constant corank one
in~$\mathcal{D}^{(i)}$.
\end{proposition}

\subsection{Jacquard's Language}

An \emph{alphabet} is a finite set $A$ whose elements are called
\emph{letters}. A \emph{word} over $A$ is a finite sequence $w=(w_{1}%
,\ldots,w_{l})$, where each $w_{i}$ belongs to $A$ for $1\leq i\leq l$. The
integer $l$ is called the \emph{length} of the word $w$ and we denote it by
$\left|  w\right|  $. The \emph{empty} word is the only word of length $0$. We
denote it by $\epsilon$. By abuse of notation, we will often write
$w_{1}\cdots w_{l}$ instead of $(w_{1},\ldots,w_{l})$ and $a^{l}$ instead of
$(a,\ldots,a)$, for any letter $a$ repeated $l$ times. Now, call $A^{*}$ the
set of all words, including the empty word, over the alphabet$~A$. A
\emph{language} over $A$ is a subset $L$ of $A^{*}$. The \emph{concatenation}
of two words $v$ and $w$ over $A$ is the word $vw=(v_{1},\ldots,v_{l}%
,w_{1},\ldots,w_{m})$, where $l=\left|  v\right|  $ and $m=\left|  w\right|
$. The concatenation of a language $L$ and a word $w$ over $A$ is the
language
\[
Lw=\{uw:\;u\in L\}.
\]
The \emph{shift} of a word $w=(w_{1},\ldots,w_{l})$ is the word $(w)^{\prime
}=(w_{1},\ldots,w_{l-1})$. By definition, we take$~(\epsilon)^{\prime
}=\epsilon$.

\medskip\ 

We define now a sequence of languages that will play an important role in this
paper, since they will describe all possible singularity types of a Goursat
structure. For a fixed $n\geq0$, consider the alphabet $A_{n}=\{a_{0}%
,\ldots,a_{n-1}\}$ (note that $A_{0}=\emptyset$). Following \cite{jacquard}
define, inductively, the \emph{Jacquard language} $J_{n}$ by $J_{0}%
=\{\epsilon\}$, $J_{1}=\{a_{0}\}$, and
\[
J_{n}=J_{n-1}a_{0}\cup J_{n-1}a_{1}\cup J_{n-2}a_{1}a_{2}\cup\ldots\cup
J_{1}a_{1}a_{2}\cdots a_{n-1},
\]
for any integer $n\geq2$. Clearly, the elements of $J_{n}$ are words over
$A_{n} $ that all have length~$n$. For example, we have $J_{2}=\{a_{0}%
a_{0},a_{0}a_{1}\}$ and
\[
J_{3}=\{a_{0}a_{0}a_{0},a_{0}a_{0}a_{1},a_{0}a_{1}a_{0},a_{0}a_{1}a_{1}%
,a_{0}a_{1}a_{2}\}.
\]
It is easy to check, using an induction argument, that
\[
\operatorname*{card}(J_{n})=3\,\operatorname*{card}(J_{n-1}%
)-\operatorname*{card}(J_{n-2}),
\]
for $n\geq2$ (see \cite{jacquard}).

\subsection{Definition of the Singularity Type}

\label{subsec-def-singularity-type}

We start with the definition of a sequence of canonical submanifolds, which
will lead to the definition of the singularity type. Let $\mathcal{D}$ be a
Goursat structure on a manifold$~M$\ of dimension$~n$. For $0\leq i\leq n-5$,
define the subset $S_{0}^{(i)}\subset M$ by
\begin{equation}
S_{0}^{(i)}=\{p\in M:\mathcal{D}^{(i)}(p)=\mathcal{C}_{i+1}%
(p)\},\label{def-S(i)0}%
\end{equation}
where the distribution $\mathcal{C}_{i}$ denotes the characteristic
distribution of $\mathcal{D}^{(i+1)}$ (see Proposition~\ref{prop-cartan}). For
$i\geq n-4$ define $S_{0}^{(i)}=\emptyset$.

Furthermore, starting from $S_{0}^{(i)}$ define, for $1\leq j\leq i$, a
sequence of subsets $M\supset S_{0}^{(i)}\supset\cdots\supset S_{j-1}%
^{(i)}\supset S_{j}^{(i)}\supset\cdots\supset S_{i}^{(i)}$ by
\begin{equation}
S_{j}^{(i)}=\{p\in S_{j-1}^{(i)}:\mathcal{D}^{(i-j)}(p)\cap T_{p}S_{j-1}%
^{(i)}\neq\mathcal{C}_{i-j}(p)\}.\label{def-S(i)j}%
\end{equation}
For $j\geq i+1$ define $S_{j}^{(i)}=\emptyset$. The above definition is
correct since, for any non-negative integers$~i$ and~$j$, the subset
$S_{j}^{(i)}\subset M$ is a smooth submanifold of $M$ (we consider an empty
set as smooth). Indeed, we have the following result, which will be proved
later, in Subsection~\ref{subsec-kr-singularity-type}.

\begin{proposition}
\label{prop-smooth-singularity}Let $\mathcal{D}$ be a Goursat structure on a
manifold $M$.

\begin{enumerate}
\item  For any non-negative integers $i$ and~$j$, the subset $S_{j}%
^{(i)}\subset M$ defined by the relations \emph{(\ref{def-S(i)0})} and
\emph{(\ref{def-S(i)j})} is either empty or a smooth submanifold of $M$ that
has codimension $j+1$ in $M$.

\item  For any non-negative integers $i$,~$j$ and $k$ such that$~k\neq j$ we
have the following relation: $S_{j}^{(i+j)}\cap S_{k}^{(i+k)}=\emptyset$.
\end{enumerate}
\end{proposition}

The fact that a point $p$ belongs to some submanifolds $S_{j}^{(i)}$, for two
given non-negative integers $i$ and$~j$, is invariantly related to the Goursat
structure at this point$~p$. This information, however, is in general
redundant. For example, if$~p$ belongs to$~S_{j}^{(i)}$ we know, by
definition, that$~p$ belongs also to$~S_{j-k}^{(i)}$, for $0\leq k\leq j$, and
that, by Proposition~\ref{prop-smooth-singularity}, it does not belong
to$~S_{j+k}^{(i+k)}$, for $k\neq0$. In the following definition we propose a
way to encode the essential part of this information in a word of the Jacquard
language (see Corollary~\ref{cor-jacquard}, at the end of
Subsection~\ref{subsec-kr-singularity-type}).

\begin{definition}
\label{def-singularity-type}Let $\mathcal{D}$ be a Goursat structure defined
in a neighborhood of a fixed point $p$ in a manifold of dimension $n$. The
\emph{singularity type} of $\mathcal{D}$ at $p$ is the word
\[
\delta_{\mathcal{D}}(p)=w_{n-4}\cdots w_{0}
\]
defined as follows: For $0\leq i\leq n-4$, we take $w_{i}=a_{j+1}$ if there
exists some integer $j\geq0$ such that $p$ belongs to $S_{j}^{(i+j)}$; we take
$w_{i}=a_{0}$ otherwise.
\end{definition}

The correctness of this definition is assured by Proposition
\ref{prop-smooth-singularity}, which states that if there exists an integer
$j\geq0$ such that $p$ belongs to $S_{j}^{(i+j)}$ then this integer is unique.
For some readers this definition may seem rather abstract at a first glance;
but it appears quickly, once computed for concrete Goursat structures, that
the singularity type really encodes essential geometric information. In fact,
we will see that the growth vector and the abnormal curves of a Goursat
structure are given by its singularity type (see Sections~\ref{sec-growth}
and~\ref{sec-abnormal}).

The singularity type should not be confused with the codes used in
\cite{cheaito-mormul} and \cite{mormul-R9} to label Kumpera-Ruiz normal forms.
Indeed, the singularity type is, by construction, an invariant of a given
Goursat structure; while the codes of \cite{cheaito-mormul} and
\cite{mormul-R9} are not invariant: a given Goursat structure can have
different codes, depending on how it is represented by a Kumpera-Ruiz normal form.

\subsection{Low Dimensional Examples}

For any Goursat structure on a manifold of dimension three or four the
singularity type is equal, at any point, to $\epsilon$ or $a_{0}$,
respectively. That is, the singularity type of a contact or an Engel structure
does not depend on the point at which the distribution is considered. This
should be compared with the singularity type of a Goursat structure on a
five-manifold, which can be either $a_{0}a_{0}$ or $a_{0}a_{1} $ at a given
point $p$, depending on whether or not the Goursat structure can be converted
into Goursat normal form in a small enough neighborhood of $p$. Indeed, for
the Goursat structure spanned by the regular Kumpera-Ruiz normal form
(\ref{KR-normal-forms-R5a}) the canonical submanifold $S_{0}^{(0)}$ is empty,
and thus the singularity type equals $a_{0}a_{0}$ at each point of
$\mathbb{R}^{5}$; for the Goursat structure spanned by the singular
Kumpera-Ruiz normal form (\ref{KR-normal-forms-R5b}) we have $S_{0}%
^{(0)}=\{x_{5}=0\}$, and thus the singularity type equals $a_{0}a_{1}$ if
$x_{5}=0$; and $a_{0}a_{0}$ if $x_{5}\neq0$.

Let us give one more example, in dimension six. Consider the distribution
$\mathcal{D}$ spanned by the following Kumpera-Ruiz normal form on
$\mathbb{R}^{6}$:
\[
\left(
\begin{array}
[c]{c}%
\tfrac\partial{\partial x_{6}}%
\end{array}
,
\begin{array}
[c]{c}%
(x_{6}+c_{6})\tfrac\partial{\partial x_{5}}+\tfrac\partial{\partial x_{4}%
}+x_{5}x_{4}\tfrac\partial{\partial x_{3}}+x_{5}x_{3}\tfrac\partial{\partial
x_{2}}+x_{5}\tfrac\partial{\partial x_{1}}%
\end{array}
\right)  ,
\]
where the constant $c_{6}$ is either equal to $0$ or $1$. For both values of
$c_{6}$, we have $S_{0}^{(0)}=\emptyset$ and $S_{0}^{(1)}=\{x_{5}=0\}$. For
$c_{6}=1$ the submanifold $S_{1}^{(1)}$ is empty (in a small enough
neighborhood of zero); for $c_{6}=0$ we have $S_{1}^{(1)}=\{x_{5}=x_{6}=0\}$.
Therefore, the singularity type of $\mathcal{D}$ at zero equals $a_{0}%
a_{1}a_{0}$ if $c_{6}=1$; and equals $a_{0}a_{1}a_{2}$ if $c_{6}=0$.

\subsection{Singularity Type of Kumpera-Ruiz Normal Forms}

\label{subsec-kr-singularity-type}

Let $\kappa^{n}$ be a Kumpera-Ruiz normal form on $\mathbb{R}^{n}$ obtained by
the inductive procedure described in Section~\ref{sec-kumpera-ruiz}. We define
inductively the word $\delta_{\kappa^{n}}$ of $J_{n-3}$ by $\delta_{\kappa
^{3}}=\epsilon$, $\delta_{\kappa^{4}}=a_{0}$, and, for $n\geq5$, by the
relations
\begin{equation}
\left\{
\begin{array}
[c]{ll}%
\delta_{\kappa^{n}}=\delta_{\kappa^{n-1}}\,a_{1} & \text{if }\kappa
^{n}=S(\kappa^{n-1})\text{;}\\
\delta_{\kappa^{n}}=\delta_{\kappa^{n-1}}\,a_{i+1} & \text{if }\kappa
^{n}=R_{0}(\kappa^{n-1})\text{ and }\delta_{\kappa^{n-1}}=\delta_{\kappa
^{n-2}}\,a_{i}\text{, }i\geq1\text{;}\\
\delta_{\kappa^{n}}=\delta_{\kappa^{n-1}}\,a_{0} & \text{otherwise;}%
\end{array}
\right. \label{eq-kumpera-ruiz-singularity-type}%
\end{equation}
where the maps $S$ and $R_{c}$ denote respectively the singular and regular
prolongations defined in Section~\ref{sec-kumpera-ruiz}. This definition leads
to a characterization of the singularity type in the coordinates of the
Kumpera-Ruiz normal form (see
Corollary~\ref{cor-singularity-type-kumpera-ruiz} below). We start with a
Proposition that will allow us to prove the results that were announced,
without proof, in Subsection~\ref{subsec-def-singularity-type}.

\begin{proposition}
\label{prop-singularity-type-claim}Consider the distribution defined on
$\mathbb{R}^{n}$ by a Kumpera-Ruiz normal form $\kappa^{n}$. For any pair of
integers $i$ and $j$ such that $0\leq j\leq i$ and for any pair of words
$w_{1}$ and $w_{2}$ such that $w=w_{1}a_{1}a_{2}\cdots a_{j+1}w_{2}$ belongs
to $J_{n-3}$ and $\left|  w_{2}\right|  =i-j$, we have $\delta_{\kappa^{n}}=w
$ if and only if zero belongs to $S_{j}^{(i)}$. Moreover, if a submanifold
$S_{j}^{(i)}$ contains zero then, in the coordinates $(x_{1},\ldots,x_{n})$ of
the Kumpera-Ruiz normal form$~\kappa^{n}$, we have
\begin{equation}
S_{j}^{(i)}=\{x_{n-i}=0,\ldots,x_{n-i+j}=0\},\label{eq-kumpera-ruiz-S(i)j}%
\end{equation}
where $0\leq i\leq n-5$ and $0\leq j\leq i$.
\end{proposition}

The following Lemma is a direct consequence of the definition of Kumpera-Ruiz
normal forms given in Section~\ref{sec-kumpera-ruiz} (see
Definition~\ref{def-kumpera-ruiz}); its proof is straightforward. Note that
below all vector fields $\kappa_{2}^{n-i}$ should be considered as vector
fields on $\mathbb{R}^{n}$, obtained from $\kappa_{2}^{n-i}$ by $i$ successive
lifts, as defined by Notation~\ref{not-lift} (see
Section~\ref{sec-kumpera-ruiz}).

\begin{lemma}
\label{lem-singularity-type-kumpera-ruiz}Let $\mathcal{D}$ be a Goursat
structure on $\mathbb{R}^{n}$ spanned by a Kumpera-Ruiz normal form
$\kappa^{n}=(\kappa_{1}^{n},\kappa_{2}^{n})$. Suppose that $\kappa^{n}%
=\sigma_{n-3}\circ\cdots\circ\sigma_{1}(\kappa^{3})$ and denote $\kappa
^{n-i}=\sigma_{n-3-i}\circ\cdots\circ\sigma_{1}(\kappa^{3})$. The derived flag
of $\mathcal{D}$ is given by
\begin{equation}
\mathcal{D}^{(i)}=(\tfrac{\partial}{\partial x_{n}},\ldots,\tfrac{\partial
}{\partial x_{n-i}},\kappa_{2}^{n-i}),\text{\quad for }0\leq i\leq
n-3\text{.}\label{eq-d(i)}%
\end{equation}
The characteristic distribution $\mathcal{C}_{i}\subset\mathcal{D}^{(i)}$ of
$\mathcal{D}^{(i+1)}$ is given by
\begin{equation}
\mathcal{C}_{i}=(\tfrac{\partial}{\partial x_{n}},\ldots,\tfrac{\partial
}{\partial x_{n-i}}),\text{\quad for }0\leq i\leq n-4\text{.}\label{eq-ci}%
\end{equation}
Moreover, if $\delta_{\kappa^{n}}=w_{1}a_{1}a_{2}\cdots a_{j}aw_{2}$, where
$\left|  w_{2}\right|  =i-j$ and $a\in\{a_{0},a_{1},a_{j+1}\}$, then we have
\begin{equation}
\kappa_{2}^{n-i+l}=%
{\textstyle\sum_{k=1}^{l}}
x_{n-i+k}\tfrac{\partial}{\partial x_{n-i+k-1}}+\tfrac{\partial}{\partial
x_{n-i-1}}+x_{n-i}\kappa_{2}^{n-i-1},\text{\quad}\label{eq-d(i-l)}%
\end{equation}
for $0\leq l\leq j-1$.
\end{lemma}

\noindent\textbf{Proof of Proposition \ref{prop-singularity-type-claim}} For
any fixed integer $i$ such that $0\leq i\leq n-5$, we will prove the
Proposition by induction on the integer $j$, for $0\leq j\leq i$. We start
with $j=0$. Assume that $\delta_{\kappa^{n}}=w_{1}aw_{2}$, where $\left|
w_{2}\right|  =i$ and $a\in A_{n-1}$. If $a\neq a_{1}$ then the definition of
$\delta_{\kappa^{n}}$ implies
\[
\kappa_{2}^{n-i}=(x_{n-i}+c_{n-i})\tfrac{\partial}{\partial x_{n-i-1}}%
+\kappa_{2}^{n-i-1};
\]
and it follows from relation (\ref{eq-d(i)}) that
\[
\mathcal{D}^{(i)}=(\tfrac{\partial}{\partial x_{n}},\ldots,\tfrac{\partial
}{\partial x_{n-i}},(x_{n-i}+c_{n-i})\tfrac{\partial}{\partial x_{n-i-1}%
}+\kappa_{2}^{n-i-1}).
\]
Together with relation (\ref{eq-ci}) this expression implies that, for any
point $p\in\mathbb{R}^{n}$, we have $\mathcal{D}^{(i)}(p)\neq\mathcal{C}%
_{i+1}(p)$, which implies that $S_{0}^{(i)}$ is empty. In particular $0\notin
S_{0}^{(i)}$. Otherwise $a=a_{1}$, and then the definition of $\delta
_{\kappa^{n}}$ implies
\[
\kappa_{2}^{n-i}=\tfrac{\partial}{\partial x_{n-i-1}}+x_{n-i}\kappa
_{2}^{n-i-1};
\]
and it follows, again from relation (\ref{eq-d(i)}), that
\[
\mathcal{D}^{(i)}=(\tfrac{\partial}{\partial x_{n}},\ldots,\tfrac{\partial
}{\partial x_{n-i}},\tfrac{\partial}{\partial x_{n-i-1}}+x_{n-i}\kappa
_{2}^{n-i-1}).
\]
Hence, for any point $p\in\mathbb{R}^{n}$, we have $\mathcal{D}^{(i)}%
(p)=\mathcal{C}_{i+1}(p)$ if and only if $x_{n-i}(p)=0$. In particular $0\in
S_{0}^{(i)}$. It follows that Proposition \ref{prop-singularity-type-claim} is
true when $j=0$.

Now, assume that Proposition \ref{prop-singularity-type-claim} is true up to
the integer $j-1$ and that
\[
\delta_{\kappa^{n}}=w_{1}a_{1}a_{2}\cdots a_{j}aw_{2},
\]
where $\left|  w_{2}\right|  =i-j$ and $a\in\{a_{1},a_{0},a_{j+1}\}$. Since
\[
S_{j-1}^{(i)}=\{x_{n-i}=\cdots=x_{n-i+j-1}=0\},
\]
we have
\[
T_{p}S_{j-1}^{(i)}=\{\tfrac{\partial}{\partial x_{n}},\ldots,\tfrac{\partial
}{\partial x_{n-i+j}},\tfrac{\partial}{\partial x_{n-i-1}},\ldots
,\tfrac{\partial}{\partial x_{1}}\}.
\]
Observe that relation (\ref{eq-d(i)}) gives
\begin{equation}%
\begin{array}
[c]{rrl}%
\mathcal{D}^{(i-j)} & = & (\tfrac{\partial}{\partial x_{n}},\ldots
,\tfrac{\partial}{\partial x_{n-i+j}},\kappa_{2}^{n-i+j}),
\end{array}
\label{eq-d(i-j)}%
\end{equation}
Moreover, it follows from relation (\ref{eq-d(i-l)}), taken for $l=j-1$, that
\begin{equation}
\kappa_{2}^{n-i+j-1}=%
{\textstyle\sum_{k=1}^{j-1}}
x_{n-i+k}\tfrac{\partial}{\partial x_{n-i+k-1}}+\tfrac{\partial}{\partial
x_{n-i-1}}+x_{n-i}\kappa_{2}^{n-i-1}.\label{eq-k2(n-i+j-1)}%
\end{equation}
If $a=a_{1}$ then the definition of $\delta_{\kappa^{n}}$ implies
\[
\kappa_{2}^{n-i+j}=\tfrac{\partial}{\partial x_{n-i+j-1}}+x_{n-i+j}\kappa
_{2}^{n-i+j-1}.
\]
Therefore the vector field $\kappa_{2}^{n-i+j}$ that appears in the relation
(\ref{eq-d(i-j)}) is, using~(\ref{eq-k2(n-i+j-1)}), given by
\[
\kappa_{2}^{n-i+j}=\tfrac{\partial}{\partial x_{n-i+j-1}}+x_{n-i+j}(%
{\textstyle\sum_{k=1}^{j-1}}
x_{n-i+k}\tfrac{\partial}{\partial x_{n-i+k-1}}+\tfrac{\partial}{\partial
x_{n-i-1}}+x_{n-i}\kappa_{2}^{n-i-1}).
\]
Hence, for any $p\in S_{j-1}^{(i)}$, we have $\mathcal{D}^{(i-j)}(p)\cap
T_{p}S_{j-1}^{(i)}=\mathcal{C}_{i-j}(p)$, which implies that $S_{j}^{(i)}$ is
empty. In particular $0\notin S_{j}^{(i)}$. If $a=a_{0}$ then the definition
of $\delta_{\kappa^{n}}$ implies
\[
\kappa_{2}^{n-i+j}=(x_{n-i+j}+c)\tfrac{\partial}{\partial x_{n-i+j-1}}%
+\kappa_{2}^{n-i+j-1},
\]
where $c\neq0$. Therefore the vector field $\kappa_{2}^{n-i+j}$ that appears
in the relation (\ref{eq-d(i-j)}) is, using~(\ref{eq-k2(n-i+j-1)}), given by
\[
\kappa_{2}^{n-i+j}=(x_{n-i+j}+c)\tfrac{\partial}{\partial x_{n-i+j-1}}+%
{\textstyle\sum_{k=1}^{j-1}}
x_{n-i+k}\tfrac{\partial}{\partial x_{n-i+k-1}}+\tfrac{\partial}{\partial
x_{n-i-1}}+x_{n-i}\kappa_{2}^{n-i-1}.
\]
Thus for any $p\in S_{j-1}^{(i)}$ we have $\mathcal{D}^{(i-j)}(p)\cap
T_{p}S_{j-1}^{(i)}=\mathcal{C}_{i-j}(p)$, which implies that $S_{j}^{(i)}$ is
empty (at least in small enough neighborhood of zero). In particular $0\notin
S_{j}^{(i)}$. Finally, if we have $a=a_{j+1}$ then the definition of
$\delta_{\kappa^{n}}$ implies
\[
\kappa_{2}^{n-i+j}=x_{n-i+j}\tfrac{\partial}{\partial x_{n-i+j-1}}+\kappa
_{2}^{n-i+j-1}.
\]
Therefore the vector field $\kappa_{2}^{n-i+j}$ that appears in the relation
(\ref{eq-d(i-j)}) is, using~(\ref{eq-k2(n-i+j-1)}), given by
\[
\kappa_{2}^{n-i+j}=%
{\textstyle\sum_{k=1}^{j}}
x_{n-i+k}\tfrac{\partial}{\partial x_{n-i+k-1}}+\tfrac{\partial}{\partial
x_{n-i-1}}+x_{n-i}\kappa_{2}^{n-i-1}.
\]
Thus for any $p\in S_{j-1}^{(i)}$ we have $\mathcal{D}^{(i-j)}(p)\cap
T_{p}S_{j-1}^{(i)}\neq\mathcal{C}_{i-j}(p)$ if and only if
\[
x_{n-i}(p)=\cdots=x_{n-i+j}(p)=0.
\]
In particular, we have $0\in S_{j}^{(i)}$, which ends the proof.\hfill$\square$

\medskip

\noindent\textbf{Proof of Proposition \ref{prop-smooth-singularity} }Consider
a Goursat structure $\mathcal{D}$ defined on a smooth manifold $M$ of
dimension $n$. \emph{First Item:} Item~(i) of
Proposition~\ref{prop-smooth-singularity} follows directly from
Proposition~\ref{prop-singularity-type-claim}, which states that, in the
coordinates of a Kumpera-Ruiz normal form, the restriction of each set
$S_{j}^{(i)}$ to a small enough neighborhood of zero is either empty or
smooth. Indeed, by Theorem~\ref{thm-kumpera-ruiz}, the Goursat structure
$\mathcal{D}$ is locally equivalent, at any point $p$ in $M$, to a
Kumpera-Ruiz normal form centered at $p$; and hence the restriction of each
set $S_{j}^{(i)}$ to a small enough neighborhood of any point $p$ in $M$ is
either empty or smooth. This obviously implies that the whole set $S_{j}%
^{(i)}$ is either empty or smooth.\smallskip

\noindent\emph{Second Item:} We will prove Item~(ii) by contradiction.
Let$~\mathcal{D}$ be a Goursat structure such that at a point$~p$ we have
$p\in S_{j}^{(i+j)}\cap S_{k}^{(i+k)}$ for a given triple of non-negative
integers $i$, $j$, and $k$, such that $k<j$. Take a Kumpera-Ruiz normal form
$\kappa^{n}$, centered at $p$, that is locally equivalent to $\mathcal{D}$ at
$p$. In the coordinates of $\kappa^{n}$, we have $0\in S_{j}^{(i+j)}\cap
S_{k}^{(i+k)}$. Let $w=\delta_{\kappa^{n}}$ be the word uniquely attached to
$\kappa^{n}$ by~(\ref{eq-kumpera-ruiz-singularity-type}), and denote the
letters of $w$ by $w=w_{n-4}\cdots w_{0}$. By
Proposition~\ref{prop-singularity-type-claim}, we have both $w_{i+l}%
=a_{j-l+1}$, for $0\leq l\leq j$, and $w_{i+l}=a_{k-l+1}$, for $0\leq l\leq
k$. In particular, we have $w_{i}=a_{j+1}$ and $w_{i}=a_{k+1}$, which is
impossible since the inequality $k<j$ implies that $a_{k+1}\neq a_{j+1}$.\hfill$\square$

\medskip

\noindent The three Corollaries listed below follow directly from
Proposition~\ref{prop-singularity-type-claim} and from the definition of the
singularity type.

\begin{corollary}
\label{cor-singularity-type-kumpera-ruiz}Let $\mathcal{D}$ be the Goursat
structure spanned on$~\mathbb{R}^{n}$ by a Kumpera-Ruiz normal form~$\kappa
^{n} $ . We have
\[
\delta_{\mathcal{D}}(0)=\delta_{\kappa^{n}},
\]
that is the singularity type at zero of $\kappa^{n}$ is given by
$\delta_{\kappa^{n}}$.
\end{corollary}

\begin{corollary}
Let $\mathcal{D}$ be a Goursat structure defined in a neighborhood of a fixed
point $p$ of a manifold of dimension $n$. For any integers such that $0\leq
j\leq i$, the point $p$ belongs to $S_{j}^{(i)}$ if and only if the
singularity type of $\mathcal{D}$ at $p$ is of the form $\delta_{\mathcal{D}%
}(p)=w_{1}a_{1}a_{2}\cdots a_{j+1}w_{2}$, with $\left|  w_{2}\right|  =$ $i-j$.
\end{corollary}

\begin{corollary}
\label{cor-jacquard} The singularity type of any Goursat structure on a
manifold of dimension $n$ belongs to the Jacquard language $J_{n-3}$.
Conversely, any word of $J_{n-3}$ is the singularity type of some Goursat structure.
\end{corollary}

\subsection{Singularity Type of the N-Trailer System}

In this Subsection, we come back to the $n$-trailer system, for which we
compute the singularity type. Our study stays very close to that of Jean
\cite{jean-trailer}. For the $n$-trailer system we define, following
\cite{jean-trailer}, the sequence of sets $\alpha_{i}$, for $i\geq0$, of real
numbers by the relations
\begin{align*}
\alpha_{1}  & =\{-\tfrac\pi2,+\tfrac\pi2\}\\
\alpha_{i+1}  & =\{\arctan\sin(\alpha)\text{, }\arctan\sin(\alpha)+\pi:\text{
}\alpha\in\alpha_{i}\}.
\end{align*}
Note that $\operatorname*{card}(\alpha_{i})=4$, for $i\geq2$. Now, consider
the $n$-trailer system $\tau^{n}$ at a configuration $p=(\xi_{1},\xi
_{2},\theta_{0},\ldots,\theta_{n})$ of $\mathbb{R}^{2}\times(S^{1})^{n+1}$.
Define, inductively, a word $\delta_{\tau^{n}}(p)=w_{1}\cdots w_{n}$ of
$J_{n}$ by $w_{1}=a_{0} $ and, for $i=2,\ldots,n$, by the relations
\begin{equation}
\left\{
\begin{array}
[c]{lcl}%
w_{i}=a_{1} &  & \text{if\quad}\theta_{i}-\theta_{i-1}\in\alpha_{1}\text{;}\\
w_{i}=a_{k+1} &  & \text{if\quad}\theta_{i}-\theta_{i-1}\in\alpha_{k+1}\text{
and }w_{i-1}=a_{k}\text{;}\\
w_{i}=a_{0} &  & \text{otherwise.}%
\end{array}
\right. \label{eq-trailer-singularity-type}%
\end{equation}
This definition leads to a characterization of the singularity type in the
coordinates of the $n$-trailer system, which coincides with the stratification
of the singular locus given in \cite{jean-trailer}.

\begin{proposition}
\label{prop-singularity-type-trailer}Let $\mathcal{D}$ be the Goursat
structure spanned by the $n$-trailer system $\tau^{n}$ on $\mathbb{R}%
^{2}\times(S^{1})^{n+1}$. We have
\[
\delta_{\mathcal{D}}(p)=\delta_{\tau^{n}}(p).
\]
Moreover, in the coordinates $(\xi_{1},\xi_{2},\theta_{0},\ldots,\theta_{n})$
of the $n$-trailer system, we have
\begin{equation}
S_{j}^{(i)}=\{\theta_{n-i}-\theta_{n-i-1}\in\alpha_{1},\ldots,\theta
_{n-i+j}-\theta_{n-i+j-1}\in\alpha_{j+1}\},\label{eq-trailer-S(i)j}%
\end{equation}
for $0\leq i\leq n-2$ and $0\leq j\leq i$.
\end{proposition}

Like in Section \ref{sec-trailer} we will use the notation $s_{i}=\sin
(\theta_{i}-\theta_{i-1})$ and $c_{i}=\cos(\theta_{i}-\theta_{i-1})$.
Moreover, we define the product $\pi_{n\,i}^{l\,k}$ by the relation
$\pi_{n\,i}^{l\,k}=%
{\textstyle\prod_{j=k}^{l}}
c_{n-i+j}$, if $0\leq k\leq l$, and by $\pi_{n\,i}^{l\,k}=1$, if $k>l$.

The proof of the Proposition~\ref{prop-singularity-type-trailer} will use the
two following Lemmas. The first one is analogous to
Lemma~\ref{lem-singularity-type-kumpera-ruiz}, of
Subsection~\ref{subsec-kr-singularity-type}, it shows that the characteristic
distributions of the $n$-trailer are rectified in $(\xi_{1},\xi_{2},\theta
_{0},\ldots,\theta_{n})$ coordinates. Its proof is straightforward and left to
the reader.

\begin{lemma}
\label{lem-singularity-type-trailer}Let $\mathcal{D}$ be the distribution
spanned by the $n$-trailer system $(\tau_{1}^{n},\tau_{2}^{n})$. The derived
flag of $\mathcal{D}$ is given by
\begin{equation}
\mathcal{D}^{(i)}=(\tfrac{\partial}{\partial\theta_{n}},\ldots,\tfrac
{\partial}{\partial\theta_{n-i}},s_{n-i}\tfrac{\partial}{\partial
\theta_{n-i-1}}+c_{n-i}\tau_{2}^{n-i-1}),\text{\quad for }0\leq i\leq
n\text{.}\label{eq-d(i)-bis}%
\end{equation}
The characteristic distributions $\mathcal{C}_{i}$ of $\mathcal{D}^{(i+1)}$
are given by
\begin{equation}
\mathcal{C}_{i}=(\tfrac{\partial}{\partial\theta_{n}},\ldots,\tfrac{\partial
}{\partial\theta_{n-i}}),\text{\quad for }0\leq i\leq n-1\text{.}%
\label{eq-ci-bis}%
\end{equation}
Moreover, we have
\begin{equation}
\mathcal{D}^{(i-j)}=(\tfrac{\partial}{\partial\theta_{n}},\ldots
,\tfrac{\partial}{\partial\theta_{n-i+j}},\sum_{k=1}^{j+1}(s_{n-i+k-1}%
)(\pi_{n\,i}^{j\,k})\tfrac{\partial}{\partial\theta_{n-i+k-2}}+\pi
_{n\,i}^{j\,0}\tau_{2}^{n-i-1}).\text{\quad}\label{eq-d(i-l)-bis}%
\end{equation}
for $0\leq i\leq n-1$ and $0\leq j\leq i$.
\end{lemma}

The Lemma below is essentially a trigonometric identity and its proof, based
on an induction argument, is also straightforward. We also leave it to the reader.

\begin{lemma}
\label{lem-label-sum}Let $(\xi_{1},\xi_{2},\theta_{0},\ldots,\theta_{n}%
)\in\mathbb{R}^{2}\times(S^{1})^{n+1}$ be a fixed point of the configuration
space of the $n$-trailer. If $\theta_{n-i+k}-\theta_{n-i+k-1}\in\alpha_{k+1}$
for $0\leq k\leq j-1$ then
\[
\sum_{k=1}^{j}(s_{n-i+k-1})(\pi_{n\,i}^{j\,k})\tfrac{\partial}{\partial
\theta_{n-i+k-2}}=s_{n-i+j-1}c_{n-i+j}(\sum_{k=1}^{j}\tfrac{\partial}%
{\partial\theta_{n-i+k-2}}).
\]
\end{lemma}

\noindent\textbf{Proof of Proposition \ref{prop-singularity-type-trailer}} To
start with, let us prove that:
\[
S_{j}^{(i)}=\{\theta_{n-i}-\theta_{n-i-1}\in\alpha_{1},\ldots,\theta
_{n-i+j}-\theta_{n-i+j-1}\in\alpha_{j+1}\}.
\]
For any fixed $i$ we will prove the result by induction on $j$. The relations
(\ref{eq-d(i)-bis}) and (\ref{eq-ci-bis}) imply that $\mathcal{D}%
^{(i)}(p)=\mathcal{C}_{i+1}(p)$ if and only if $c_{n-i}(p)=0$. That is, if and
only if $\theta_{n-i}-\theta_{n-i-1}\in\alpha_{1}$, which implies that the
Proposition to be true for $j=0$. Now, assume the Proposition true up to
$j-1$. The relation (\ref{eq-d(i-l)-bis}) implies that $\mathcal{D}^{(i-j)}$
is given by
\[
(\tfrac\partial{\partial\theta_{n}},\ldots,\tfrac\partial{\partial
\theta_{n-i+j}},\sum_{k=1}^{j+1}(s_{n-i+k-1})(\pi_{n\,i}^{j\,k})\tfrac
\partial{\partial\theta_{n-i+k-2}}+\pi_{n\,i}^{j\,0}\tau_{2}^{n-i-1}).
\]
The induction assumption, together with Lemma \ref{lem-label-sum}, implies
that for any point $p$ that belongs to $S_{j-1}^{(i)}$ the linear subspace
$\mathcal{D}^{(i-j)}(p)$ is equal to
\[
(\tfrac\partial{\partial\theta_{n}},\ldots,\tfrac\partial{\partial
\theta_{n-i+j}},s_{n-i+j}\tfrac\partial{\partial\theta_{n-i+j-1}}%
+(s_{n-i+j-1}c_{n-i+j})(\sum_{k=1}^{j}\tfrac\partial{\partial\theta_{n-i+k-2}%
})).
\]
The induction argument says that
\[
S_{j-1}^{(i)}=\{\theta_{n-i}-\theta_{n-i-1}\in\alpha_{1},\ldots,\theta
_{n-i+j-1}-\theta_{n-i+j-2}\in\alpha_{j}\}.
\]
Since
\[
T_{p}S_{j-1}^{(i)}=(\tfrac\partial{\partial\theta_{n}},\ldots,\tfrac
\partial{\partial\theta_{n-i+j}},\sum_{k=1}^{j+1}\tfrac\partial{\partial
\theta_{n-i+k-2}},\tfrac\partial{\partial\theta_{n-i-2}},\ldots\tfrac
\partial{\partial\theta_{0}},\tfrac\partial{\partial\xi_{2}},\tfrac
\partial{\partial\xi_{1}}),
\]
we have $\mathcal{D}^{(i-j)}(p)\cap T_{p}S_{j-1}^{(i)}\neq\mathcal{C}_{i-j}(p)
$ if and only if $s_{n-i+j}(p)=s_{n-i+j-1}(p)c_{n-i+j}(p)$. That is, we have
$p\in S_{j}^{(i)}$ if and only if $p\in S_{j-1}^{(i)}$ and $\theta
_{n-i+j}-\theta_{n-i+j-1}\in\alpha_{j+1}$, which ends the induction argument.

Now, the form of $S_{j}^{(i)}$ obtained in the previous paragraph together
with the definitions of$~\delta_{\mathcal{D}}$ and$~\delta_{\tau^{n}}$, imply
directly that $\delta_{\mathcal{D}}(p)=\delta_{\tau^{n}}(p).$\hfill$\square$

\section{Growth Vector}

\label{sec-growth}

\subsection{Lie Flag and Growth Vector}

The \emph{Lie flag} of a distribution $\mathcal{D}$ is the sequence of modules
of vector fields $\mathcal{D}_{0}\subset\mathcal{D}_{1}\subset\cdots$ defined
inductively by
\begin{equation}
\mathcal{D}_{0}=\mathcal{D}\text{\quad and\quad}\mathcal{D}_{i+1}%
=\mathcal{D}_{i}+[\mathcal{D}_{0},\mathcal{D}_{i}]\text{, \quad for }%
i\geq0\text{.}\label{lie-flag}%
\end{equation}
This sequence should not be confused with the derived flag (\ref{derived-flag}%
), introduced at the beginning of the article. In general these two sequences
are different. Nevertheless, for any point $p$ in the underlying manifold~$M$,
the inclusion $\mathcal{D}_{i}(p)\subset\mathcal{D}^{(i)}(p)$ holds, for
$i\geq0$. Note that for a Goursat structure, unlike the elements of the
derived flag, the elements of the Lie flag are \emph{not} necessarily
distributions of constant rank.

A distribution $\mathcal{D}$ is \emph{completely nonholonomic} if, for each
point $p$ in $M$, there exists an integer $N(p)$ such that $\mathcal{D}%
_{N(p)}(p)=T_{p}M$. The smallest such integer is called the \emph{nonholonomy
degree} of $\mathcal{D}$ at $p$ and we denote it by $N_{p}$. For a Goursat
structure on a manifold of dimension $n$, the inequality $N_{p}\leq2^{n-3}$
holds for each point $p$ in $M$ (see e.g. \cite{laumond-singularities}). For
the $n$-trailer system, sharper bounds were obtained in~\cite{jean-trailer},
\cite{luca-risler}, and~\cite{sordalen-bound}. It follows from our
Theorem~\ref{thm-trailer-universal}, which states that any Goursat structure
is locally equivalent to the $n$-trailer system, that they hold also for any
Goursat structure.

\begin{definition}
Let $\mathcal{D}$ be a completely nonholonomic distribution. Put
$d_{i}(p)=\dim\mathcal{D}_{i}(p)$, for $0\leq i\leq N_{p}$. The \emph{growth
vector} at $p$ of the distribution $\mathcal{D}$ is the finite sequence
$(d_{0}(p),\ldots,d_{N_{p}}(p))$.
\end{definition}

Recall that if at a given point a Goursat structure can be converted into
Goursat normal form~(\ref{goursat-normal-form}) then this point is called
regular and that otherwise it is called singular (see
Section~\ref{sec-kumpera-ruiz}). The set of singular points is called the
\emph{singular locus}. An elegant characterization of this set, that
emphasizes the importance of the growth vector in the study of Goursat
structures, has been obtained by Murray \cite{murray-nilpotent}. A different
characterization can be found in~\cite{kumpera-ruiz} and~\cite{libermann}.

\begin{theorem}
[Murray]\label{thm-murray}Let $p$ be a point in a manifold $M$ of dimension
$n$. A Goursat structure on $M$ can be converted into Goursat normal form in a
small enough neighborhood of $p$ if and only if $\mathcal{D}_{i}%
(p)=\mathcal{D}^{(i)}(p)$, for $0\leq i\leq n-2$.
\end{theorem}

\subsection{Growth Vector of the N-Trailer System}

Let $d=(d_{0},\ldots,d_{N})$ be a finite sequence of integers such that
$d_{0}=2$, $d_{N}=n$, and $d_{i}\leq d_{i+1}\leq d_{i}+1$, for $0\leq i\leq
N-1$. The \emph{dual }of the sequence $d$ is the sequence $d^{\ast}%
=(d_{2}^{\ast},\ldots,d_{n}^{\ast})$ defined by
\[
d_{i}^{\ast}=\operatorname*{card}\{j\geq0:d_{j}<i\}+1\text{,\quad for }2\leq
i\leq n.
\]
In other words, the integer $d_{i}^{\ast}$ indicates the first position,
starting from the left, where the integer $i$ appears in $d$. We obviously
have $d_{2}^{\ast}=1$ and $d_{n}^{\ast}=N+1$. It is trivial to check that each
sequence $d$ is uniquely defined by its dual $d^{\ast}$. For example, we have
the following dual sequences: $(2,3,4,5,6)^{\ast}=(1,2,3,4,5)$,
$(2,3,4,5,5,5,6)^{\ast}=(1,2,3,4,7)$, and $(2,3,4,4,5,5,5,6)^{\ast
}=(1,2,3,5,8)$.

Now, following \cite{jean-trailer}, we define a set of functions that will
allow us to obtain a formula that gives the growth vector of an arbitrary
Goursat structure at an arbitrary point, as a function of its singularity type
at this point. We start with Jean's formula \cite{jean-trailer} for the
$n$-trailer. Recall that $J_{n}$ denotes the Jacquard language (see
Section~\ref{sec-singularity-type}) and that the shift of a word is defined by
$(w_{1}\cdots w_{n})^{\prime}=w_{1}\cdots w_{n-1}$ and $(\epsilon)^{\prime
}=\epsilon$ (we will denote $(w^{\prime})^{\prime}$ by $w^{\prime\prime} $).

For any $i\geq2$, we define functions $\beta_{i}:%
{\textstyle\bigcup\nolimits_{n\geq i-3}}
J_{n}\rightarrow\mathbb{N}$ . We take $\beta_{2}(w)=1,$ $\beta_{3}(w)=2,$ and
$\beta_{4}(w)=3$, for any word $w$ in $%
{\textstyle\bigcup\nolimits_{n\geq i-3}}
J_{n}$. If $i\geq5$ then we define inductively, for any word $w$ in $%
{\textstyle\bigcup\nolimits_{n\geq i-3}}
J_{n}$,
\[
\left\{
\begin{array}
[c]{lll}%
\beta_{i}(w)=\quad\beta_{i-1}(w^{\prime})+\beta_{i-2}(w^{\prime\prime}) &  &
\text{if }w=(w^{\prime})a_{1}\\
\beta_{i}(w)=2\,\,\beta_{i-1}(w^{\prime})-\beta_{i-2}(w^{\prime\prime}) &  &
\text{if }w=(w^{\prime})a_{k}\text{ and }k\geq2;\\
\beta_{i}(w)=\quad\beta_{i-1}(w^{\prime})+1 &  & \text{if }w=(w^{\prime}%
)a_{0}\text{,}%
\end{array}
\right.
\]
For example, for the word $a_{0}a_{1}a_{0}$, we have:
\begin{align*}
\beta_{5}(a_{0}a_{1}a_{0})  & =\beta_{4}(a_{0}a_{1})+1=3+1=4\\
\beta_{6}(a_{0}a_{1}a_{0})  & =\beta_{5}(a_{0}a_{1})+1=(\beta_{4}(a_{0}%
)+\beta_{3}(\epsilon))+1=(3+2)+1=6.
\end{align*}
An other example is given for the word $a_{0}a_{1}a_{2}$, for which we have:
\begin{align*}
\beta_{5}(a_{0}a_{1}a_{2})  & =\beta_{4}(a_{0}a_{1})+1=3+1=4\\
\beta_{6}(a_{0}a_{1}a_{2})  & =2\beta_{5}(a_{0}a_{1})-\beta_{4}(a_{0}%
)=2(\beta_{4}(a_{0})+\beta_{3}(\epsilon))-3=7.
\end{align*}

Having recalled the functions $\beta_{i}$ we are now able to recall the
formula, obtained by Jean~\cite{jean-trailer}, that gives the growth vector of
the $n$-trailer system.

\begin{theorem}
[Jean]\label{thm-jean}Consider the $n$-trailer system at a given point $p$ of
its configuration space $\mathbb{R}^{2}\times(S^{1})^{n+1}$ at which it has
singularity type $\delta_{\tau^{n}}(p)$. The sequence of integers
$(d_{2}^{\ast}(p),\ldots,d_{n+3}^{\ast}(p))$ dual to the growth vector of the
$n$-trailer system at$~p$ is given by $d_{i}^{\ast}(p)=\beta_{i}(\delta
_{\tau^{n}}(p))$.
\end{theorem}

\subsection{Growth Vector of Goursat Structures}

The following result is fundamental. It shows that the growth vector of any
Goursat structure is a function of its singularity type.

\begin{theorem}
\label{thm-growth-vector}Let $\mathcal{D}$ be a Goursat structure on a
manifold $M$ of dimension $n\geq3$, defined in a neighborhood of a given point
$p$ in $M$ that has singularity type $\delta_{\mathcal{D}}(p)$. The sequence
of integers $(d_{2}^{\ast}(p),\ldots,d_{n}^{\ast}(p))$ dual to its growth
vector at $p$ is given by $d_{i}^{\ast}(p)=\beta_{i}(\delta_{\mathcal{D}}%
(p))$.\medskip\noindent\noindent
\end{theorem}

\noindent\textbf{Proof of Theorem \ref{thm-growth-vector}} Let $\mathcal{D}$
be a Goursat structure on a manifold $M$ of dimension $n\geq3$, defined in a
neighborhood of a given point $p$ in $M$. By
Theorem~\ref{thm-trailer-universal}, the Goursat structure $\mathcal{D}$ is
locally equivalent at $p$ to the $n$-trailer system, considered around a well
chosen point $q$ of its configuration space. By Theorem~\ref{thm-jean}, the
sequence of integers $(d_{2}^{\ast}(q),\ldots,d_{n+3}^{\ast}(q))$ dual to the
growth vector of the $n$-trailer system at $q$ is given by $d_{i}^{\ast
}(q)=\beta_{i}(\delta_{\tau^{n}}(q))$. By
Proposition~\ref{prop-singularity-type-trailer}, the singularity type of the
$n$-trailer system at $q$ equals $\delta_{\tau^{n}}(q)$. Since the singularity
type is invariant under diffeomorphisms, we have $\delta_{\mathcal{D}%
}(p)=\delta_{\tau^{n}}(q)$. Since the growth vector is invariant under
diffeomorphisms, the sequence of integers $(d_{2}^{\ast}(p),\ldots
,d_{n+3}^{\ast}(p))$ dual to the growth vector of $\mathcal{D}$ at $p$ is
given by $d_{i}^{\ast}(p)=\beta_{i}(\delta_{\tau^{n}}(q))=\beta_{i}%
(\delta_{\mathcal{D}}(p))$.\hfill$\square$

\medskip

\noindent The latter result obviously implies the following one, which gives
the formula for the growth vector of an arbitrary Kumpera-Ruiz normal form.

\begin{corollary}
Let $\kappa^{n}$ be a Kumpera-Ruiz normal form on $\mathbb{R}^{n}$, for
$n\geq3$. The sequence $(d_{2}^{\ast},\ldots,d_{n}^{\ast})$ dual to its growth
vector at zero is given by $d_{i}^{\ast}=\beta_{i}(\delta_{\kappa^{n}})$.
\end{corollary}

\subsection{Growth Vector and Singularity Type}

We proved in the previous Subsection (Theorem \ref{thm-growth-vector}) that
the singularity type of any Goursat structure at a given point determines its
growth vector at this point. Now, we will prove the converse of this fact.

\begin{theorem}
\label{thm-singularity-type-growth-vector}Two Goursat structures have the same
growth vector at a given point if and only if they have the same singularity
type at this point.
\end{theorem}

The proof of Theorem~\ref{thm-singularity-type-growth-vector} will be based on
two Lemmas:

\begin{lemma}
\label{lem-growth-01}Let $i$ and $k$ be two integers such that $i\geq1$ and
$0\leq k\leq i-1$. For any word $w$ in $%
{\textstyle\bigcup\nolimits_{n\geq1}}
J_{n}$ we have the following relations:

\begin{enumerate}
\item $\beta_{i+4}(wa_{1}a_{2}\cdots a_{i})=2\,i+3;$

\item $\beta_{i+4}(wa_{1}a_{2}\cdots a_{i-k}a_{0}^{k})=2\,i-k+3;$

\item $\beta_{i+4}(wc_{1}\cdots c_{i})=i+3;$
\end{enumerate}

\noindent where $c_{j}$, for $1\leq j\leq i$, are any letters satisfying
$c_{j}\neq a_{1}$.
\end{lemma}

\noindent\textbf{Proof of Lemma \ref{lem-growth-01}} \emph{First
Item}\textbf{. }Item (i) is true if $i=1$ because, for any word $w$ in $%
{\textstyle\bigcup\nolimits_{n\geq1}}
J_{n}$, we have
\[%
\begin{array}
[c]{lll}%
\beta_{5}(wa_{1}) & = & \beta_{4}(w)+\beta_{3}(w^{\prime})\\
& = & 3+2=2\cdot1+3.
\end{array}
\]
It is also true if $i=2$ because, for any word $w$ in $%
{\textstyle\bigcup\nolimits_{n\geq1}}
J_{n}$, we have
\[%
\begin{array}
[c]{lll}%
\beta_{6}(wa_{1}a_{2}) & = & 2\cdot\beta_{5}(wa_{1})-\beta_{4}(w)\\
& = & 2\,\cdot5-3=2\,\cdot2+3.
\end{array}
\]
Now proceed by induction on $i\geq3$. Assume that Item (i) is true up to
$i-1$. Then, for any word $w$ in $%
{\textstyle\bigcup\nolimits_{n\geq1}}
J_{n}$, we have
\[%
\begin{array}
[c]{lll}%
\beta_{i+4}(wa_{1}a_{2}\cdots a_{i}) & = & 2\cdot\beta_{(i-1)+4}(wa_{1}%
a_{2}\cdots a_{i-1})-\beta_{(i-2)+4}(wa_{1}a_{2}\cdots a_{i-2})\\
& = & 2\cdot(2\cdot(i-1)+3)-(2\cdot(i-2)+3)\\
& = & 2\,i+3.
\end{array}
\]

\noindent\emph{Second Item.} Let us proceed by induction on $i$. It follows
from Item~(i) that, for $i\geq1$, Item (ii) is true for $i=1$ and $k=0$.
Assume that Item~(ii) is true up to $i-1$ for any $0\leq k\leq i-2$. Then we
have, for $1\leq k\leq i-1$ and for any $w$ in $%
{\textstyle\bigcup\nolimits_{n\geq1}}
J_{n}$, the following relation:
\[%
\begin{array}
[c]{lll}%
\beta_{i+4}(wa_{1}a_{2}\cdots a_{i-k}a_{0}^{k}) & = & \beta_{(i-1)+4}%
(wa_{1}a_{2}\cdots a_{(i-1)-(k-1)}a_{0}^{k-1})+1\\
& = & 2\cdot(i-1)-(k-1)+3\\
& = & 2\,i-k+3.
\end{array}
\]
Since, by Item (i), Item (ii) is true for $k=0$, it follows that Item (ii)
holds for any $i\geq1$ and any $0\leq k\leq i-1$.

\noindent\emph{Third Item. }Item (iii) is true if $i=1$ because $\beta
_{5}(wc_{1})=4$ for any word $w$ in $%
{\textstyle\bigcup\nolimits_{n\geq1}}
J_{n}$ (recall that $c_{1}\neq a_{1}$). It is also true if $i=0$. Now proceed
by induction on $i $. Assume that this Item is true up to $i-1$, then we have
either
\[%
\begin{array}
[c]{lll}%
\beta_{i+4}(wc_{1}\cdots c_{i}) & = & 2\cdot\beta_{(i-1)+4}(wc_{1}\cdots
c_{i-1})-\beta_{(i-2)+4}(wc_{1}\cdots c_{i-2})\\
& = & 2\cdot((i-1)+3)-((i-2)+3)\\
& = & \,i+3
\end{array}
\]
or
\[%
\begin{array}
[c]{lll}%
\beta_{i+4}(wc_{1}\cdots c_{i}) & = & \beta_{(i-1)+4}(wc_{1}\cdots
c_{i-1})+1\\
& = & ((i-1)+3)+1\\
& = & \,i+3,
\end{array}
\]
which ends the proof.\hfill$\square$

\begin{lemma}
\label{lem-growth-02}Let $i$ be an integer such that $i\geq5$. Consider two
words $w_{1}$ and $w_{2}$ of the Jacquard language $J_{l}$, with $l\geq i-3$,
such that:

\begin{enumerate}
\item $\beta_{i}(w_{1})>\beta_{i}(w_{2});$

\item $\beta_{i-1}(w_{1}^{\prime})\geq\beta_{i-1}(w_{2}^{\prime});$

\item $\beta_{i}(w_{1})-\beta_{i-1}(w_{1}^{\prime})\geq\beta_{i}(w_{2}%
)-\beta_{i-1}(w_{2}^{\prime}).$
\end{enumerate}

\noindent Then, for any integer $k\geq1$ and for any word $w$ such that
$w_{1}w$ and $w_{2}w$ belong to $J_{k+l}$, we have $\beta_{i+k}(w_{1}%
w)>\beta_{i+k}(w_{2}w)$.
\end{lemma}

\noindent\textbf{Proof of Lemma \ref{lem-growth-02}} Consider two words $w_{1}
$ and $w_{2}$ in $J_{l}$, with $l\geq i-3$, that satisfy conditions (i)-(iii).
Let $a$ be any letter such that $w_{1}a$ and $w_{2}a$ belong to $J_{l+1}$.
Then we have the three following cases:

\noindent If $a=a_{0}$ then
\[%
\begin{array}
[c]{lll}%
\beta_{i+1}(w_{1}a_{0}) & = & \beta_{i}(w_{1})+1\\
\beta_{i+1}(w_{2}a_{0}) & = & \beta_{i}(w_{2})+1.
\end{array}
\]

\noindent If $a=a_{1}$ then
\[%
\begin{array}
[c]{lll}%
\beta_{i+1}(w_{1}a_{1}) & = & \beta_{i}(w_{1})+\beta_{i-1}(w_{1}^{\prime})\\
\beta_{i+1}(w_{2}a_{1}) & = & \beta_{i}(w_{2})+\beta_{i-1}(w_{2}^{\prime}).
\end{array}
\]

\noindent If $a=a_{j}$ then
\[%
\begin{array}
[c]{lll}%
\beta_{i+1}(w_{1}a_{j}) & = & \beta_{i}(w_{1})+\beta_{i}(w_{1})-\beta
_{i-1}(w_{1}^{\prime})\\
\beta_{i+1}(w_{2}a_{j}) & = & \beta_{i}(w_{2})+\beta_{i}(w_{2})-\beta
_{i-1}(w_{2}^{\prime}).
\end{array}
\]

Therefore, in any case, the words $w_{1}a$ and $w_{2}a$ satisfy the three
conditions (i)-(iii); and in particular, we have $\beta_{i+1}(w_{1}%
a)>\beta_{i+1}(w_{2}a)$. Hence the Lemma is true for $k=1$. An induction
argument on the length of $w$, based on the same relations as for $k=1$, ends
the proof.\hfill$\square$

\medskip\ 

\noindent\textbf{Proof of Theorem \ref{thm-singularity-type-growth-vector}} By
Theorem~\ref{thm-growth-vector}, if two Goursat structures have the same
singularity types at $p$ and $\tilde p$, respectively, then they have the same
growth vector at $p$ and $\tilde p$, respectively. Now, we will prove the
converse. Suppose that $w$ and $\tilde w$ are the singularity types of two
distributions $\mathcal{D}$ and $\tilde{\mathcal{D}}$ at $p$ and $\tilde p$,
respectively, that is $w=\delta_{\mathcal{D}}(p)$ and $\tilde w=\delta
_{\mathcal{D}}(p)$.\ We will show that if $w\neq\tilde w$ then there exists an
integer$~i_{0}$ such that $\beta_{i_{0}}(w)\neq\beta_{i_{0}}(\tilde w)$.

It is easy to check that if $w$ and $\tilde w$ are two words of the Jacquard
language $J_{n}$ such that $w\neq\tilde w$ then there exists (after a
permutation of $w$ and $\tilde w$, if necessary) three words $z$, $v$, and
$\tilde v$ such that both $w=vz$ and $\tilde w=\tilde vz$, and which satisfy
either
\[
\left\{
\begin{array}
[c]{lll}%
v & = & ua_{1}a_{2}\cdots a_{i-k}a_{0}^{k}\\
\tilde v & = & \tilde uc_{1}c_{2}\cdots c_{i},
\end{array}
\right.
\]
where $0\leq k\leq i-1$ and $c_{j}\neq a_{1}$ for $1\leq j\leq i$, or
\[
\left\{
\begin{array}
[c]{lll}%
v & = & ua_{1}a_{2}\cdots a_{i-k}a_{0}^{k}\\
\tilde v & = & \tilde ua_{1}a_{2}\cdots a_{i-l}a_{0}^{l},
\end{array}
\right.
\]
where $k\neq l$.

For each of these two cases we can apply Lemma \ref{lem-growth-01}. In the
first case we have $\beta_{i+4}(v)=2\,i-k+3$; while $\beta_{i+4}(\tilde
v)=i+3$. Since $k\leq i-1$ we have $\beta_{i+4}(v)\neq\beta_{i+4}(\tilde v)$.
In the second case we have $\beta_{i+4}(v)=2\,i-k+3$; while $\beta
_{i+4}(\tilde v)=2\,i-l+3$. Since $k\neq l$ we have $\beta_{i+4}(v)\neq
\beta_{i+4}(\tilde v)$. Therefore, in both cases, we have $\beta_{i+4}%
(v)\neq\beta_{i+4}(\tilde v)$; but $\beta_{i+3}(v)=\beta_{i+3}(\tilde v)$,
since by the Item (iii) of Lemma \ref{lem-growth-01} they are both equal to
$i+2$. Put $i_{0}=(i+4)+\left|  z\right|  $. By Lemma \ref{lem-growth-02}, we
have $\beta_{i_{0}}(w)\neq\beta_{i_{0}}(\tilde w)$.\hfill$\square$

\subsection{Computing the Singularity Type}

Up to now, we have worked with a definition of the singularity type that uses
the submanifolds $S_{j}^{(i)}$. Although being geometric, that is independent
of a description of the Goursat structure in particular coordinates, it does
not tell us how to compute this invariant (unless we know how to compute all
$S_{j}^{(i)}$ explicitly). In order to fill this gap we give the following
Proposition, which yields to a constructive procedure to compute the
singularity type of any Goursat structure in terms of its growth vector. Its
proof is straightforward.

\begin{proposition}
Let $\mathcal{D}$ be a Goursat structure considered in a neighborhood of a
point $p$ that belongs to a manifold of dimension $n\geq5$. For $0\leq i\leq
n-5$ and $1\leq j\leq i$, the point $p$ belongs to $S_{j}^{(i)}$ if and only
if the growth vector at $p$ of the distribution $\mathcal{D}^{(i-j)} $ starts
with
\[
(i-j+2,i-j+3,\ldots,i+2,i+3,i+4,\ldots,i+4,i+5),
\]
where the integer $i+4$ is repeated $j+2$ times.
\end{proposition}

\section{Abnormal Curves}

\label{sec-abnormal}

\subsection{Integral and Abnormal Curves}

Let $M$ be a smooth manifold of dimension $n$ and let $\mathcal{A}$ be a
set-valued map $\mathcal{A}:M\rightarrow TM$ such that $\mathcal{A}(p)\subset
T_{p}M$, for each point $p$ in $M$. Note that we do not ask $\mathcal{A}(p)$
to be a linear subspace of $T_{p}M$, but just a subset of $T_{p}M $. Neither
we ask $\mathcal{A}$ to be smooth. An \emph{integral curve} of $\mathcal{A}$
is an absolutely continuous map $x:I\rightarrow M$, from an interval
$I\subset\mathbb{R}$ to $M$, such that $\dot x(t)$ belongs to $\mathcal{A}%
(x(t))$ for almost all $t$ in$~I$. A \emph{nontrivial lift} of $x(\cdot)$ is
an absolutely continuous map $P:I\rightarrow T^{*}M$ such that$~P(t)$ belongs
to $T_{x(t)}^{*}M$ and $P(t)\neq0$ for each $t$ in$~I$.

Locally, all integral curves of a rank $k$ distribution $\mathcal{D}%
=(f_{1},\ldots,f_{k})$ can be described as solutions of an (underdetermined)
ordinary differential equation. Indeed, for any given integral curve
$x(\cdot)$ of $\mathcal{D}$ we can clearly find $k$ real-valued measurable
functions $u_{i}$, for $1\leq i\leq k$, such that
\begin{equation}
\dot{x}(t)=%
{\textstyle\sum_{i=1}^{k}}
f_{i}(x(t))\,u_{i}(t)\label{control-system}%
\end{equation}
holds for almost all $t$ in $I$. These functions $u_{i}$ are called
\emph{controls}. Observe that the controls associated to an integral curve are
not uniquely defined. In control theory, an overdetermined differential
equation of the form (\ref{control-system}), where the functions $u_{i}$ for
$1\leq i\leq k$ can be taken as arbitrary measurable functions, is called a
\emph{control system}. Informally, the system (\ref{control-system}) can be
seen as a ``parametrization'' of the set of all integral curves of
$\mathcal{D}$ by~$k$ real-valued measurable functions.

Roughly speaking, a solution $x(\cdot)$ of (\ref{control-system}) is abnormal
if it is a singular point of the end-point mapping or, equivalently, if the
linearization of the control system along $x(\cdot)$ is not controllable. Many
equivalent definitions of the concept of abnormal curves are available (see
e.g. the papers \cite{agrachev-sarychev-jmsec}, \cite{bonnard-kupka},
\cite{bryant-hsu}, \cite{sussmann-liu}, \cite{zhitomirskii-nice}, the survey
article \cite{montgomery-survey}, and the references given there). The
definition that we will use is the one that appears in Pontryagin's Maximum
principle~\cite{pontryagin-book}. For further details, we refer the reader to
the above mentioned works.

Since the results of this section will be local we can work in a coordinate
chart $x:M\rightarrow\mathbb{R}^{n}$. Denote by $(x,p)$ the corresponding
coordinates on $T^{*}M$. In these coordinates, the \emph{Hamiltonian} of the
control system (\ref{control-system}) associated to a distribution
$\mathcal{D}=(f_{1},\ldots,f_{k})$ is the function defined on $\mathbb{R}%
^{n}\times\mathbb{R}^{n}\times\mathbb{R}^{k}$ by
\[
H(x,p,u)=\left\langle
\begin{array}
[c]{c}%
p
\end{array}
,
\begin{array}
[c]{c}%
{\textstyle\sum\limits_{i=1}^{k}}
f_{i}(x)u_{i}%
\end{array}
\right\rangle ,
\]
where both $x$ and $p$ belong to $\mathbb{R}^{n}$ and $u=(u_{1},\ldots,u_{k})$
belongs to $\mathbb{R}^{k}$ and $\left\langle \cdot,\cdot\right\rangle $
denotes the pairing between vector fields and differential forms.

\begin{definition}
\label{def-abnormal}An integral curve $x:I\rightarrow\mathbb{R}^{n}$,
corresponding to a measurable control $u:I\rightarrow\mathbb{R}^{k}$, of the
control system (\ref{control-system}) is called \emph{abnormal} if it admits a
nontrivial lift $(x(\cdot),p(\cdot))$ such that
\[%
\begin{array}
[c]{ccr}%
\dot{x}(t) & = & \dfrac{\partial H(x(t),p(t),u(t))}{\partial p}\\
\dot{p}(t) & = & -\dfrac{\partial H(x(t),p(t),u(t))}{\partial x}%
\end{array}
\]
and
\[%
\begin{array}
[c]{ccc}%
\dfrac{\partial H(x(t),p(t),u(t))}{\partial u} & = & 0
\end{array}
\]
for almost all $t$ in $I$.
\end{definition}

By definition, an integral curve of a distribution $\mathcal{D}=(f_{1}%
,\ldots,f_{k})$ is \emph{abnormal} if it is an abnormal curve of the
corresponding control system. It is well known that the abnormal curves of
$\mathcal{D}$ depend neither on the choice of coordinates nor on the vector
fields $f_{1},\ldots,f_{k}$ chosen to span the distribution.

Let $I\subset\mathbb{R}$ be an interval. For any $t_{0}\in I$ and for any
$\varepsilon>0$, denote by $I_{\varepsilon}(t_{0})$ the intersection
$I\cap[t_{0}-\varepsilon,t_{0}+\varepsilon]$. An integral curve
$x:I\rightarrow M$ is \emph{locally abnormal} if for each $t_{0}$ in $I$ there
exists a small enough $\varepsilon>0$ such that the restriction of $x(\cdot)$
to $I_{\varepsilon}(t_{0})$ is abnormal.

\subsection{Abnormal Curves of Goursat Structures}

Let $\mathcal{D}$ be a Goursat structure on a manifold $M$ of dimension
$n\geq3$. Recall that its singularity type can be computed using the sequence
of canonical manifolds defined, for $0\leq i\leq n-5$, by
\[
S_{0}^{(i)}=\{q\in M:\mathcal{D}^{(i)}(q)=\mathcal{C}_{i+1}(q)\}
\]
and, for $1\leq j\leq i$, by
\[
S_{j}^{(i)}=\{q\in S_{j-1}^{(i)}:\mathcal{D}^{(i-j)}(q)\cap T_{q}S_{j-1}%
^{(i)}\neq\mathcal{C}_{i-j}(q)\},
\]
where the distributions $\mathcal{C}_{i}$ are the canonical distributions of
Proposition~\ref{prop-cartan} (see Section~\ref{sec-singularity-type}). Assume
that for two given non-negative integers $i$ and $j$, such that $0\leq i+j\leq
n-5$ we have $S_{j}^{(i+j)}\neq\emptyset$. In this case, we can define on
$S_{j}^{(i+j)}$ a smooth distribution $\mathcal{A}_{j}^{(i)}$ by taking
\[
\mathcal{A}_{j}^{(i)}(q)=\mathcal{D}^{(i)}(q)\cap T_{q}S_{j}^{(i+j)},
\]
for each point $q$ in $S_{j}^{(i+j)}$. It is easy to check, using a
Kumpera-Ruiz normal form, that $\mathcal{A}_{j}^{(i)}$ is indeed a smooth
distribution and that its rank is $i+1$. Although each $\mathcal{A}_{j}^{(i)}$
is defined only on $S_{j}^{(i+j)}$, we can extend the definition of
$\mathcal{A}_{j}^{(i)}$ to $M$ by taking $\mathcal{A}_{j}^{(i)}(q)=0$ for all
points~$q$ that do not belong to $S_{j}^{(i+j)}$ and thus consider
$\mathcal{A}_{j}^{(i)}$ as a set valued map defined everywhere on $M$. This
extension allows us to define, for any $0\leq i\leq n-5$, a subset
$\mathcal{A}^{(i)}\subset TM$ by
\[
\mathcal{A}^{(i)}(q)=\mathcal{C}_{i}(q)\cup\left(
{\textstyle\bigcup\limits_{0\leq j\leq n-i-5}}
\mathcal{A}_{j}^{(i)}(q)\right)  ,
\]
for each point $q$ in $M$. Note that, usually, the subset $\mathcal{A}%
^{(i)}\subset TM$ is not a distribution.

By definition, we take $\mathcal{A}^{(n-4)}=\mathcal{C}_{n-4}$. Moreover, we
define $\mathcal{A}^{(n-3)}$ as the characteristic distribution of
$\mathcal{D}^{(n-3)}$, which is equal to $\mathcal{C}_{n-4}$ if $n\geq4$ and
equal to $\{0\}$ if $n=3$. Finally, we take $\mathcal{A}^{(n-2)}=\emptyset$.
Observe that the set-valued maps $\mathcal{A}^{(n-3)}=\{0\}$ and
$\mathcal{A}^{(n-2)}=\emptyset$ are different. Indeed, the first one has
trivial integral curves (points); while the second one has no integral curves
at all.

\begin{theorem}
\label{thm-abnormal-curves}Consider a Goursat structure $\mathcal{D}$ defined
on a manifold of dimension~$n$ and fix an integer$~i$ such that $0\leq i\leq
n-2$. An integral curve of$~\mathcal{D}^{(i)}$ is locally abnormal if and only
if it is an integral curve of$~\mathcal{A}^{(i)}$.
\end{theorem}

For $i=n-4$, $n-3$, and $n-2$, the distribution $\mathcal{D}^{(i)}$ is of rank
$n-2$, $n-1$, and $n$, respectively, and the proof of Theorem
\ref{thm-abnormal-curves} follows easily from well known results. Indeed, if
$i=n-4$ then the distribution $\mathcal{D}^{(i)}$, which is of rank $n-2$, can
be transformed into a direct generalization of Engel's normal form
(\cite{kumpera-ruiz}, \cite{martin-rouchon-driftless},
\cite{zhitomirskii-growth}, and~\cite{zhitomirskii-survey}) given by
Theorem~\ref{thm-weber-problem} (see Appendix~\ref{sec-weber}), where we have
to take $k=n-2$ and $m=2$. In this case, the abnormal curves of $\mathcal{D}%
^{(n-4)}$ are clearly the integral curves of $\mathcal{C}_{n-4}$ (see Lemma
\ref{lem-abnormal-01} below). If $i=n-3$ then the distribution $\mathcal{D}%
^{(i)}$, which is of rank $n-1$, is annihilated locally by a $1$-form$~\omega$
such that $d\omega\wedge\omega\neq0$ and $(d\omega)^{2}\wedge\omega=0$. This
property is equivalent to the fact that the characteristic distribution of
$\mathcal{D}^{(n-3)}$ is of corank $2$ in $\mathcal{D}^{(n-3)}$
(see~\cite{bryant-chern-gardner-goldschmidt-griffiths}), and it implies that
$\mathcal{D}^{(n-3)}$ is locally given by the normal form of
Theorem~\ref{thm-weber-problem}, where $k=n-1$ and $m=1$. Note, however, that
this form does not follow from Theorem~\ref{thm-weber-problem} whose
condition, when $m=1$, is only necessary but not sufficient. In this case, it
is straightforward to see that the abnormal curves of $\mathcal{D}^{(n-3)}$
are the integral curves of the characteristic distribution of $\mathcal{D}%
^{(n-3)}$, which is an involutive distribution that has rank $n-3$. Finally,
if $i=n-2$ then the situation is even simpler because $\mathcal{D}^{(n-2)}=TM
$, which implies that $\mathcal{D}^{(n-2)}$ has no abnormal curves at all.
Hence the only values of $i$ that will be considered in the proof of
Theorem~\ref{thm-abnormal-curves} are $0\leq i\leq n-5$.

In order to explain further the meaning of Theorem~\ref{thm-abnormal-curves}
we would like to emphasize the following points, relative to the geometric
structure of $\mathcal{A}^{(i)}$ and its integral curves. These facts follow
directly from our study of the singularity type (see
Section~\ref{sec-singularity-type}) and will be used in the proof of
Theorem~\ref{thm-abnormal-curves}.

(i) Although for each point $q$ in $M$ we have, by definition,
\[
\mathcal{A}^{(i)}(q)=\mathcal{C}_{i}(q)\cup\left(
{\textstyle\bigcup\limits_{0\leq j\leq n-i-5}}
\mathcal{A}_{j}^{(i)}(q)\right)  ,
\]
the relations $S_{k}^{(i+k)}\cap S_{j}^{(i+j)}=\emptyset$ for $k\neq j$ (see
Proposition \ref{prop-smooth-singularity}) imply that, for a fixed point~$q$,
only two possibilities can occur. Indeed, we have either
\[
\mathcal{A}^{(i)}(q)=\mathcal{C}_{i}(q)\text{\quad or\quad}\mathcal{A}%
^{(i)}(q)=\mathcal{C}_{i}(q)\cup\mathcal{A}_{j}^{(i)}(q),
\]
for a unique integer $j$ such that $0\leq j\leq n-i-5$. In other words, for
each point$~q$ the subset $\mathcal{A}^{(i)}(q)\subset T_{q}M$ is the union
(not the sum!) of either one or two linear subspaces of$~T_{q}M$. Note that if
$i\geq1$ then $\mathcal{C}_{i}(q)\cap\mathcal{A}_{j}^{(i)}(q)=\mathcal{C}%
_{i-1}(q)$.

(ii) For $0\leq i\leq n-5$, we define the set $K_{i}=%
{\textstyle\bigcup_{j=i}^{n-5}}
S_{0}^{(j)}$; for any other value of $i$ we take $K_{i}=\emptyset$. We will
call this set the \emph{singular locus} of $\mathcal{D}^{(i)}$. If $i=n-4$,
$n-3$, or $n-2$ then, by definition, the singular locus is empty, which
explains why these cases are simpler. If $i=0$ then this definition agrees
with the one given in Section \ref{sec-growth} for the singular locus of
$\mathcal{D}$. It follows directly from the properties of the submanifolds
$S_{0}^{(j)}$ (see Proposition \ref{prop-singularity-type-claim}) that $K_{i}$
is a stratified manifold. In fact, in Kumpera-Ruiz normal form coordinates,
this set is an algebraic variety defined by a single polynomial equation of
the form $%
{\textstyle\prod\nolimits_{r=0}^{m-1}}
x_{k_{r}}^{r}=0$, where the integer$~m$ corresponds to the number of
singularities of$~\mathcal{D}^{(i)}$, which is, in general, smaller than the
number of singularities of$~\mathcal{D}$. For any point $q$ that does not
belong to $K_{i}$ we clearly have $\mathcal{A}^{(i)}(q)=\mathcal{C}_{i}(q)$.
Note, however, that there exist points of $K_{i}$ for which we also have
$\mathcal{A}^{(i)}(q)=\mathcal{C}_{i}(q)$.

(iii) For $0\leq i\leq n-5$, define the set $L_{i}=%
{\textstyle\bigcup_{j=0}^{n-i-5}}
S_{j}^{(i+j)}$. Since for any $j$ we have $S_{j}^{(i)}\subset S_{0}^{(i)}$, it
follows that $L_{i}\subset K_{i}$. For $0\leq i\leq n-5$, the set of points
such that $\mathcal{A}^{(i)}(q)\neq\mathcal{C}_{i}(q)$ is precisely~$L_{i}$.
In other words, the set $\mathcal{A}^{(i)}(q)$ is a linear subspace of
$T_{q}M$ if and only if $q$ does not belong to $L_{i}$. Unlike $K_{i}$, the
set $L_{i}$ is always a smooth submanifold of $M$. Note, however, that~$L_{i}$
can have several connected components and that the dimensions of these
components can be different. Nevertheless, in a small enough neighborhood $U$
of any of its points, the submanifold $L_{i}$ is connected and coincides with
one and only one of the smooth manifolds $S_{j}^{(i+j)}\cap U$.

For example, in the case of a distribution spanned by a Kumpera-Ruiz normal
form on~$\mathbb{R}^{n}$, the set$~L_{i}$ is connected. If non-empty, the set
$L_{i}$ is a codimension $j+1$ linear subspace of $\mathbb{R}^{n}$, where$~j$
is the only integer such that $S_{j}^{(i+j)}$ is non-empty. In the case of the
$n$-trailer system, the situation is quite different. For example, for the
two-trailer system, the submanifold$~L_{0}$ has two connected components,
given by $\{\theta_{2}-\theta_{1}=\pi/2\}$ and $\{\theta_{2}-\theta_{1}%
=-\pi/2\}$. Each of them has codimension$~1$. For the three-trailer system,
the submanifold $L_{0}$ has six connected components given, respectively, by
$\{\theta_{3}-\theta_{2}=\pi/2\}$, $\{\theta_{3}-\theta_{2}=-\pi/2\}$,
$\{\theta_{3}-\theta_{2}=\pi/4$; $\theta_{2}-\theta_{1}=\pi/2\}$,
$\{\theta_{3}-\theta_{2}=-3\pi/4$; $\theta_{2}-\theta_{1}=\pi/2\}$,
$\{\theta_{3}-\theta_{2}=-\pi/4$; $\theta_{2}-\theta_{1}=-\pi/2\}$,
$\{\theta_{3}-\theta_{2}=3\pi/4$; $\theta_{2}-\theta_{1}=-\pi/2\}$. Two of
them have codimension$~1$; four of them have codimension$~2$.

\ \medskip\ 

We consider now a more detailed example. Let $\mathcal{D}$ be the distribution
spanned by the following Kumpera-Ruiz normal form on$~\mathbb{R}^{7}$:
\[
\left(
\begin{array}
[c]{c}%
\tfrac\partial{\partial x_{7}}%
\end{array}
,
\begin{array}
[c]{c}%
(x_{7}+c_{7})\tfrac\partial{\partial x_{6}}+\tfrac\partial{\partial x_{5}%
}+x_{6}\left(  \tfrac\partial{\partial x_{4}}+x_{5}\left(  x_{4}\tfrac
\partial{\partial x_{3}}+x_{3}\tfrac\partial{\partial x_{2}}+\tfrac
\partial{\partial x_{1}}\right)  \right)
\end{array}
\right)  ,
\]
where $c_{7}$ is either equal to $0$ or $1$. When $c_{7}=1$ the singularity
type of $\mathcal{D}$ at zero is $a_{0}a_{1}a_{1}a_{0}$ and the growth vector
at zero is $(2,3,4,5,5,6,6,6,7)$; while when $c_{7}=0$ the singularity type is
$a_{0}a_{1}a_{1}a_{2}$ and the growth vector $(2,3,4,5,5,5,6,6,6,6,7)$.

In both cases, we have
\[
S_{0}^{(0)}=\emptyset\text{, }S_{0}^{(1)}=\{x_{6}=0\}\text{, and }S_{0}%
^{(2)}=\{x_{5}=0\}.
\]
Therefore, the singular loci of the distributions $\mathcal{D}^{(0)}$,
$\mathcal{D}^{(1)}$, and $\mathcal{D}^{(2)}$ are given respectively by
\[
K_{0}=\{x_{6}x_{5}=0\}\text{, }K_{1}=\{x_{6}x_{5}=0\}\text{, and }%
K_{2}=\{x_{5}=0\}.
\]
If $c_{7}=1$ then, in a small enough neighborhood of zero, we have
$S_{1}^{(1)}=\emptyset$; but if $c_{7}=0$ then we have $S_{1}^{(1)}%
=\{x_{7}=x_{6}=0\}$. In both cases we have $S_{1}^{(2)}=\emptyset$.

If $c_{7}=1$ then, in a small enough neighborhood $U$ of zero, we have
$\mathcal{A}^{(0)}=\mathcal{C}_{0}=(\tfrac\partial{\partial x_{7}})$, which is
a smooth distribution on $U$; but if $c_{7}=0$ then the subset $\mathcal{A}%
^{(0)}$ coincides with the smooth distribution $\mathcal{C}_{0}=(\tfrac
\partial{\partial x_{7}})$ outside $L_{0}=\{x_{7}=x_{6}=0\}$ while for any
point $p$ of $L_{0}$ we have
\[
\mathcal{A}^{(0)}(p)=(\tfrac\partial{\partial x_{7}})(p)\cup(\tfrac
\partial{\partial x_{5}})(p),
\]
which is clearly not a distribution. In both cases, we have $\mathcal{A}%
^{(1)}=\mathcal{C}_{1}=(\tfrac\partial{\partial x_{7}},\tfrac\partial{\partial
x_{6}})$ outside $L_{1}=\{x_{6}=0\}$ while for any point $p$ of $L_{1}$ we
have
\[
\mathcal{A}^{(1)}(p)=(\tfrac\partial{\partial x_{7}},\tfrac\partial{\partial
x_{6}})(p)\cup(\tfrac\partial{\partial x_{7}},\tfrac\partial{\partial x_{5}%
})(p).
\]
Finally, we have $\mathcal{A}^{(2)}=(\tfrac\partial{\partial x_{7}}%
,\tfrac\partial{\partial x_{6}},\tfrac\partial{\partial x_{5}})$ outside
$L_{2}=\{x_{5}=0\}$ while for any point $p$ of $L_{2}$ we have
\[
\mathcal{A}^{(2)}(p)=(\tfrac\partial{\partial x_{7}},\tfrac\partial{\partial
x_{6}},\tfrac\partial{\partial x_{5}})(p)\cup(\tfrac\partial{\partial x_{7}%
},\tfrac\partial{\partial x_{6}},\tfrac\partial{\partial x_{4}})(p).
\]

\medskip\ 

We proceed now to the proof of Theorem \ref{thm-abnormal-curves}, which states
that an integral curve of $\mathcal{D}^{(i)}$ is locally abnormal if and only
if it is an integral curve of $\mathcal{A}^{(i)}$. The proof will be based on
the three following Lemmas.

\begin{lemma}
\label{lem-abnormal-01}Consider a Goursat structure $\mathcal{D}$ defined on a
manifold of dimension$~n$, and fix an integer $i$ such that $0\leq i\leq n-4$.
An integral curve of $\mathcal{D}^{(i)}$ that has an empty intersection with
the singular locus $K_{i}$ is locally abnormal if and only if it is an
integral curve of $\mathcal{C}_{i}$, and thus of $\mathcal{A}^{(i)}$.
\end{lemma}

\noindent\textbf{Proof of Lemma \ref{lem-abnormal-01} }Let $\gamma
:I\rightarrow M$ be an integral curve of $\mathcal{D}^{(i)}$ that does not
intersect the singular locus $K_{i}$. Since we are outside $K_{i}$ it is easy
to show, using a direct generalization of Goursat's normal form (see
\cite{kumpera-ruiz} and \cite{martin-rouchon-driftless}), given by
Theorem~\ref{thm-weber-problem}, that for any fixed $t_{0}$ in $I$ we can find
a local coordinate chart $x:U\rightarrow\mathbb{R}^{n}$ centered at
$\gamma(t_{0})$ and such that:
\[
\mathcal{D}^{(i)}=\left(
\begin{array}
[c]{c}%
\tfrac\partial{\partial x_{1}}%
\end{array}
,\ldots,
\begin{array}
[c]{c}%
\tfrac\partial{\partial x_{i+1}}%
\end{array}
,
\begin{array}
[c]{c}%
x_{i+1}\tfrac\partial{\partial x_{i+2}}+\cdots+x_{n-2}\tfrac\partial{\partial
x_{n-1}}+\tfrac\partial{\partial x_{n}}%
\end{array}
\right)  .
\]
Chose a small enough $\varepsilon>0$ such that the restriction of $\gamma$ to
$I_{\varepsilon}(t_{0})$ is completely contained in the open set $U$. Then,
the curve $x\circ\gamma:I_{\varepsilon}(t_{0})\rightarrow\mathbb{R}^{n}$,
which will be denoted shortly by $x(\cdot)$, is almost everywhere a solution
of the following control system
\begin{equation}%
\begin{array}
[c]{lcl}%
\dot x_{1} & = & u_{1}\\
& \vdots & \\
\dot x_{i} & = & u_{i}\\
\dot x_{i+1} & = & u_{i+1}\\
\dot x_{i+2} & = & x_{i+1}u_{i+2}\\
& \vdots & \\
\dot x_{n-1} & = & x_{n-2}u_{i+2}\\
\dot x_{n} & = & u_{i+2}.
\end{array}
\label{eq-lem-abnormal-02a}%
\end{equation}
Since the coordinate chart is centered at $\gamma(t_{0})$, we have
$x(t_{0})=0$. The Hamiltonian of this system is given by
\[
H(x,p,u)=%
{\textstyle\sum_{k=1}^{i+1}}
p_{k}u_{k}+%
{\textstyle\sum_{k=i+2}^{n-1}}
p_{k}x_{k-1}u_{i+2}+p_{n}u_{i+2}.
\]
Therefore, the curve $x(\cdot)$ is abnormal if and only if there exists a
non-trivial lift $(x(\cdot),p(\cdot))$ that satisfies, almost everywhere, the
following differential equation
\begin{equation}%
\begin{array}
[c]{lcl}%
\dot p_{1} & = & 0\\
& \vdots & \\
\dot p_{i} & = & 0\\
\dot p_{i+1} & = & -p_{i+2}u_{i+2}\\
& \vdots & \\
\dot p_{n-2} & = & -p_{n-1}u_{i+2}\\
\dot p_{n-1} & = & 0\\
\dot p_{n} & = & 0
\end{array}
\label{eq-lem-abnormal-02b}%
\end{equation}
and, moreover, is such that $p_{k}=0$, for $1\leq k\leq i+1$, and $p_{n}=-%
{\textstyle\sum\nolimits_{k=i+2}^{n-1}}
p_{k}x_{k-1}$. The latter condition is a consequence of $\frac{\partial
H}{\partial u}=0$.

\smallskip\ \ 

\noindent\emph{Necessity.} Assume that $x(\cdot)$ is not an integral curve of
$\mathcal{C}_{i}$. We will prove that $x(\cdot)$ is not abnormal. In the
coordinates of (\ref{eq-lem-abnormal-02a}) we have $\mathcal{C}_{i}%
=(\tfrac\partial{\partial x_{1}},\ldots,\tfrac\partial{\partial x_{i+1}}).$
Since $x(\cdot)$ is not an integral curve of $\mathcal{C}_{i}$, there exists a
measurable subset $I_{0}\subset I_{\varepsilon}(t_{0})$ such that the Lebesgue
measure of $I_{0}$ is not zero and $u_{i+2}(t)\neq0$ for each $t$ in $I_{0}$.
If $x(\cdot)$ is abnormal then $p(\cdot)$ is such that $p_{i+1}(t)=0$ for each
$t$ in $I_{\varepsilon}(t_{0})$. Therefore, we have $\dot p_{i+1}=0$ almost
everywhere on $I_{0}$. Indeed, note that if an absolutely continuous function
$f$ on $I_{0}$ is such that $f(t)=0$ for almost all $t$ in $I_{0}$ then
$f^{\prime}(t)=0$ for almost all $t$ in $I_{0}$. But $\dot p_{i+1}%
=-p_{i+2}u_{i+2}$ and $u_{i+2}\neq0$ imply $p_{i+2}=0$ almost everywhere on
$I_{0}$, which gives $\dot p_{i+2}=0$ almost everywhere on $I_{0}$. We can
repeat the previous argument to obtain $p_{k}=0$, for $1\leq k\leq n-1 $,
almost everywhere on $I_{0}$. Since $p_{n}=-%
{\textstyle\sum\nolimits_{k=i+2}^{n-1}}
p_{k}x_{k-1}$, we have also $p_{n}=0$ almost everywhere on $I_{0}$. This gives
$p_{k}=0$, almost everywhere on $I_{0}$, for $1\leq k\leq n$, which is
impossible since $p$ must be non-trivial.

\smallskip

\noindent\emph{Sufficiency.} Now, assume that $x(\cdot)$ is an integral curve
of $\mathcal{C}_{i}$. In order to prove that $x(\cdot)$ is abnormal, we will
consider the lift defined by $p_{k}=0$ for $1\leq k\leq n$, with the exception
of $p_{n-1}$, for which any non-zero real constant can be taken. Since
$x(\cdot)$ is an integral curve of $\mathcal{C}_{i}$ we must have
$u_{i+2}(t)=0$ almost everywhere on $I_{\varepsilon}(t_{0})$, which implies
that $p(\cdot)$ satisfies (\ref{eq-lem-abnormal-02b}). Moreover, since
$x(t_{0})=0$, we have $x_{k}(t)=x_{k}(t_{0})=0$, for each $t$ in
$I_{\varepsilon}(t_{0})$ and for $i+2\leq k\leq n$. Thus $p_{n}$, which was
taken to be zero, satisfies $p_{n}=-%
{\textstyle\sum\nolimits_{k=i+2}^{n-1}}
p_{k}x_{k-1}$ (recall that $p_{i+2}=0$). In other words $p(\cdot)$ satisfies
$\frac{\partial H}{\partial u}=0$. Finally, since $p_{n-1}\neq0$, our lift is
non-trivial, which implies that $x(\cdot)$ is abnormal.\hfill$\square$

\begin{lemma}
\label{lem-abnormal-02}Consider a Goursat structure $\mathcal{D}$ defined on a
manifold of dimension$~n$ and fix an integer $i$ such that $0\leq i\leq n-5 $.
Let $x(\cdot)$ be the restriction of an integral curve of $\mathcal{D}^{(i)}$
to the interval $I_{\varepsilon}(t_{0})$, where $\varepsilon>0$. If a fixed
measurable subset $I_{0}\subset\mathbb{R}$ is such that for each $t$
in$~I_{0}\cap I_{\varepsilon}(t_{0})$ we have $\dot{x}(t)\notin\mathcal{A}%
^{(i)}(x(t))$ then, for a small enough $\varepsilon>0$, we have $x(t)\notin
K_{i}$ for almost all $t$ in $I_{0}\cap I_{\varepsilon}(t_{0})$.
\end{lemma}

\noindent\textbf{Proof of Lemma \ref{lem-abnormal-02}} Let $x:I_{\varepsilon
}(t_{0})\rightarrow\mathbb{R}^{n}$ be the restriction to the interval
$I_{\varepsilon}(t_{0})$, where $\varepsilon>0$, of an integral curve of
$\mathcal{D}^{(i)}$. Firstly, if $x(t_{0})\notin K_{i}$ then there exists a
small enough$~\varepsilon$ such that $x(\cdot)$ does not intersect $K_{i}$ and
thus, in this case, the Lemma is trivially true. Secondly, if the Lebesgue
measure of $I_{0}$ is $0$ then the Lemma is also trivially true. Finally, if
the closure of $I_{0}$ does not contain $t_{0}$ then for a small
enough$~\varepsilon$ the Lebesgue measure of $I_{0}\cap I_{\varepsilon}%
(t_{0})$ will be$~0$ and thus the Lemma will be, once more, trivially true.
Hence, from now on, we will only consider curves such that $x(t_{0})$ belongs
to $K_{i}$, the Lebesgue measure of $I_{0}$ is not $0$, and the closure of
$I_{0}$ contains $t_{0}$. Moreover, once a small enough $\varepsilon>0$ has
been fixed, we will denote also by $I_{0}$ the intersection $I_{0}\cap
I_{\varepsilon}(t_{0})$. That is, we will assume that $I_{0}\subset
I_{\varepsilon}(t_{0})$.

\medskip

For any such integral curve $x(\cdot)$ of $\mathcal{D}^{(i)}$ it is easy to
prove, using a direct generalization of Kumpera-Ruiz's normal form, given by
Theorem \ref{thm-weber-kumpera-ruiz} (with a double indexation of coordinates,
like in Corollary \ref{cor-cras-indexation}), that there exist coordinates on
$\mathbb{R}^{n}$ in which $x(\cdot)$ is a solution of the following control
system:
\begin{align}
\dot x_{1}^{0}  & =u_{1}\nonumber\\
\dot x_{2}^{0}  & =u_{2}\nonumber\\
& \vdots\label{double-indexation}\\
\dot x_{i}^{0}  & =u_{i}\nonumber\\
\dot x_{i+1}^{0}  & =u_{i+1}\nonumber\\
\dot x_{q}^{p}  & =\left(
{\textstyle\prod\limits_{0\leq r\leq p-1}}
x_{k_{r}}^{r}\right)  (x_{q-1}^{p}+c_{q-1}^{p})u_{i+2}\quad%
\begin{array}
[c]{l}%
\text{for }i+2\leq q\leq k_{0}\text{ if }p=0\text{ and}\\
\text{for }2\leq q\leq k_{p}\text{ if }1\leq p\leq m
\end{array}
\nonumber\\
\dot x_{1}^{p+1}  & =\left(
{\textstyle\prod\limits_{0\leq r\leq p-1}}
x_{k_{r}}^{r}\right)  u_{i+2}\quad%
\begin{array}
[c]{l}%
\text{for }0\leq p\leq m,
\end{array}
\nonumber
\end{align}

where $x=(x_{1}^{0},x_{2}^{0},\ldots,x_{k_{0}}^{0},x_{1}^{1},\ldots,x_{k_{m}%
}^{m},x_{1}^{m+1})$ and $x(t_{0})=0$ (recall that we assume that $x(t_{0})$
belongs to the singular locus, which means that $m\geq1$). Moreover, the
integers $k_{i}$ that appear in (\ref{double-indexation}) satisfy $i+1\leq
k_{0}\leq n-4$ and $k_{1}\geq1,\ldots,k_{m-1}\geq1,k_{m}\geq3,k_{m+1}=1$ and $%
{\textstyle\sum\nolimits_{r=0}^{m+1}}
k_{r}=n$. Observe that the number $m\geq1$ is the number of singularities of
$\mathcal{D}^{(i)}$, which can be smaller than the number of singularities of
$\mathcal{D}$. In these coordinates, the singular locus is given by
\[
K_{i}=\left\{
{\textstyle\prod\limits_{r=0}^{m-1}}
x_{k_{r}}^{r}=0\right\}  .
\]
If $c_{q}^{0}=0$, for all $i+1\leq q\leq k_{0}-1$, then the only integer $j$
such that $S_{j}^{(i)}$ contains zero is $j=k_{0}-(i+1)$. Thus $L_{i}$ is
given (see Proposition~\ref{prop-singularity-type-claim}) by
\[
L_{i}=S_{k_{0}-1-i}^{(k_{0}-1)}=\left\{  x_{i+1}^{0}=x_{i+2}^{0}%
=\cdots=x_{k_{0}}^{0}=0\right\}  .
\]
Note that if for some $i+1\leq q\leq k_{0}-1$ we have $c_{q}^{0}\neq0$ then
the submanifold $L_{i}$ does not contain zero (it is locally empty at zero).

For each point $p$ of $\mathbb{R}^{n}$, we have
\[
\mathcal{C}_{i}(p)=(\tfrac\partial{\partial x_{1}^{0}},\ldots,\tfrac
\partial{\partial x_{i+1}^{0}})(p).
\]
If $c_{q}^{0}=0$, for $i+1\leq q\leq k_{0}-1$, then we have
\[
\mathcal{A}^{(i)}(p)=(\tfrac\partial{\partial x_{1}^{0}},\ldots,\tfrac
\partial{\partial x_{i+1}^{0}})(p)\cup(\tfrac\partial{\partial x_{1}^{0}%
},\ldots,\tfrac\partial{\partial x_{i}^{0}},\tfrac\partial{\partial x_{1}^{1}%
})(p),
\]
for each point $p$ in $L_{i}$ and $\mathcal{A}^{(i)}(p)=\mathcal{C}_{i}(p)$,
outside $L_{i}$. If for some $i+1\leq q\leq k_{0}-1$ we have $c_{q}^{0}\neq0$
then, in a small enough neighborhood $U$ of zero, we have $\mathcal{A}%
^{(i)}(p)=\mathcal{C}_{i}(p)$ for each point $p$ in $U$.

\medskip

Recall that the integral curve $x:I_{\varepsilon}(t_{0})\rightarrow
\mathbb{R}^{n}$ of $\mathcal{D}^{(i)}$ is such that $x(t_{0})=0$. Assume that,
at a given $t$ of $I_{\varepsilon}(t_{0})$, the velocity $\dot x(t)$ exists
and satisfies (\ref{double-indexation}). Then, if $\varepsilon$ is small
enough, the velocity $\dot x(t)$ belongs to $\mathcal{A}^{(i)}(x(t))$ if and
only if we have $u_{i+2}(t)=0$ or the three following conditions hold: (i)
$u_{i+1}(t)=0$ and (ii) $x_{q}^{0}(t)=0$, for $i+1\leq q\leq k_{0}$, and (iii)
$c_{q}^{0}=0$, for $i+1\leq q\leq k_{0}-1$.

Now, suppose that for each $t$ in $I_{0}\subset I_{\varepsilon}(t_{0})$ the
velocity $\dot{x}(t)$ exists and is such that $\dot{x}(t)\notin\mathcal{A}%
^{(i)}(x(t))$. Recall that we can assume that the Lebesgue measure of $I_{0}$
is not $0$ and the closure of $I_{0}$ contains $t_{0}$. For each $t$ in
$I_{0}$ we have $u_{i+2}(t)\neq0$. Additionally: (a) If $c_{q}^{0}=0$ for
$i+1\leq q\leq k_{0}-1$ then we can represent the subset $I_{0}$ as
$I_{0}=I_{1}\cup I_{2}$ (with, in general, a non empty intersection of $I_{1}$
and $I_{2}$), where $I_{1}$ is the set of points where (i) is not satisfied
and $I_{2}$ is the set of points where (ii) is not satisfied. (b) If there
exist an integer $i+1\leq q\leq k_{0}-1$ such that $c_{q}^{0}\neq0$ then
$I_{0}=\{t\in I_{\varepsilon}(t_{0}):u_{i+2}(t)\neq0\}$, provided that
$\varepsilon$ is small enough. We are going to show that, in both cases, we
have $x(t)\notin K_{i}$, for almost all $t$ in $I_{0}$.

\medskip

\noindent\emph{Case }(a):\emph{\ Subset }$I_{1}$. For each $t$ in $I_{1}$ we
have both $u_{i+2}(t)\neq0$ and $u_{i+1}(t)\neq0$. Therefore, we have $\dot
x_{i+1}^{0}\neq0$ almost everywhere on $I_{1}$, which implies that
$x_{i+1}^{0}\neq c_{i+1}^{0}$ almost everywhere on $I_{1}$. Indeed, note that
if an absolutely continuous function $f$ on $I_{1}$ is such that $f^{\prime
}(t)\neq0$ for almost all $t$ in $I_{1}$ then, for any constant $c$, the
measure of the set $\{t\in I_{1}:f(t)=c\}$ is zero.

Now, using an induction argument we will show, successively, that $x_{i+1}%
^{0}\neq c_{i+1}^{0}$, $x_{i+2}^{0}\neq c_{i+2}^{0}$,..., $x_{k_{0}-1}^{0}\neq
c_{k_{0}-1}^{0}$, $x_{k_{0}}^{0}\neq0$, $x_{1}^{1}\neq c_{1}^{1}$,...,
$x_{k_{1}-1}^{1}\neq c_{k_{1}-1}^{1}$, $x_{k_{1}}^{1}\neq0$,..., $x_{k_{m-1}%
}^{m-1}\neq0$, almost everywhere on $I_{1}$. Suppose that this assumption is
true up to $x_{q-1}^{p}$. We have two cases: either $q\leq k_{p}$ or
$q=k_{p}+1$. If $q\leq k_{p}$ then $\dot{x}_{q}^{p}=%
{\textstyle\prod\nolimits_{0\leq r\leq p-1}}
(x_{k_{r}}^{r})(x_{q-1}^{p}+c_{q-1}^{p})u_{i+2}$. Since $x_{k_{r}}^{r}\neq0$,
for $0\leq r\leq p-1$, and $x_{q-1}^{p}\neq c_{q-1}^{p} $ and $u_{i+2}\neq0$,
almost everywhere on $I_{1}$, we have $\dot{x}_{q}^{p}(t)\neq0$ for almost all
$t$ in $I_{1}$. This implies, almost everywhere on $I_{1}$, that $x_{q}%
^{p}\neq c_{q}^{p}$ if $q\leq k_{p}-1$ or that $x_{q}^{p}\neq0$ if $q=k_{p}$.
Otherwise $q=k_{p}+1$ and in this case $\dot{x}_{1}^{p+1}=%
{\textstyle\prod\nolimits_{0\leq r\leq p-1}}
(x_{k_{r}}^{r})u_{i+2}$. Since $x_{k_{r}}^{r}\neq0$, for $0\leq r\leq p-1$,
and $u_{i+2}\neq0$, almost everywhere on $I_{1}$, we have $\dot{x}_{1}%
^{p+1}(t)\neq0$ for almost all $t$ in $I_{1}$, which implies $x_{1}^{p+1}\neq
c_{1}^{p+1}$ almost everywhere on $I_{1}$. This ends the induction argument.
In particular, we have proved that $x_{k_{r}}^{r}(t)\neq0$ for almost all $t$
in $I_{1}$, for each $0\leq r\leq m-1$. Now, recall that the singular locus is
given by the relation $%
{\textstyle\prod\nolimits_{r=0}^{m-1}}
x_{k_{r}}^{r}=0$. It thus follows that we have $x(t)\notin K_{i}$ for almost
all $t$ in $I_{1}$.

\medskip

\noindent\emph{Case }(a):\emph{\ Subset }$I_{2}$. We can represent the subset
$I_{2}$ as $I_{2}=I_{2}^{i+1}\cup\cdots\cup I_{2}^{k_{0}}$ , where $I_{2}%
^{q}=\{t\in I_{0}:x_{q}^{0}(t)\neq0\}$. Observe that, in general, the
intersection of these subsets will be non empty. Now on each subset $I_{2}%
^{q}$, of positive Lebesgue measure, we can follow the same proof as for the
subset $I_{1}$, starting the induction argument with $x_{q}^{0}$. For each one
of these subsets the conclusion is the same: we have $x(t)\notin K_{i}$ for
almost all $t$ in $I_{2}^{q} $.

\medskip

\noindent\emph{Case }(b): We have $u_{i+2}(t)\neq0$ for each $t$ in $I_{0}$
and, moreover, there exists $q$ such that $c_{q}^{0}\neq0$, where $i+1\leq
q\leq k_{0}-1$. Since $c_{q}^{0}\neq0$, we can take a smaller $\varepsilon>0$,
if necessary, in order to have $\dot x_{q+1}^{0}(t)=(x_{q}^{0}(t)+c_{q}%
^{0})u_{i+2}\neq0$ for each $t$ in $I_{0}$. This implies that $x_{q+1}^{0}\neq
c_{q+1}^{0}$ almost everywhere on $I_{0}$. The rest of the proof follows like
for the subset $I_{1}$ and the conclusion is the same: we have $x(t)\notin
K_{i}$ for almost all $t$ in $I_{0}$. The only difference being that the
induction argument starts with $x_{q+1}^{0}$ instead of $x_{i+1}^{0}$%
.\hfill$\square$

\begin{lemma}
\label{lem-abnormal-03}Consider a Goursat structure $\mathcal{D}$ defined on a
manifold of dimension$~n$ and fix an integer $i$ such that $0\leq i\leq n-5 $.
An integral curve of $\mathcal{A}^{(i)}$ that intersects the singular locus
$K_{i}$ is locally, in a small enough neighborhood of any point of
intersection, an abnormal curve of $\mathcal{D}^{(i)}$.
\end{lemma}

\noindent\textbf{Proof of Lemma \ref{lem-abnormal-03}} Let $x:I_{\varepsilon
}(t_{0})\rightarrow\mathbb{R}^{n}$ be the restriction to the interval
$I_{\varepsilon}(t_{0})$, where $\varepsilon>0$, of an integral curve of
$\mathcal{D}^{(i)}$ that intersects the singular locus at $x(t_{0})$. It is
easy to prove (see Lemma~\ref{lem-abnormal-normal-form}) that, for a small
enough $\varepsilon>0$, there exist coordinates such that the integral curve
$x(\cdot)$ is a solution of the following control system:
\begin{equation}%
\begin{array}
[c]{lcl}%
\dot x_{1} & = & u_{1}\\
& \vdots & \\
\dot x_{i} & = & u_{i}\\
\dot x_{i+1} & = & u_{i+1}\\
\dot x_{i+2} & = & (x_{i+1}+c_{i+1})u_{i+2}\\
& \vdots & \\
\dot x_{k_{0}} & = & (x_{k_{0}-1}+c_{k_{0}-1})u_{i+2}\\
\dot x_{k_{0}+1} & = & x_{k_{0}}u_{i+2}\\
\dot x_{k_{0}+2} & = & u_{i+2}\\
\dot x_{k_{0}+3} & = & x_{k_{0}}x_{k_{0}+2}u_{i+2}\\
\dot x_{k_{0}+4} & = & x_{k_{0}}x_{k_{0}+3}u_{i+2}\\
\dot x_{j} & = & x_{k_{0}}\phi_{j}(\overline{x}_{k_{0}+3})u_{i+2}\text{\quad
for }k_{0}+5\leq j\leq n,
\end{array}
\label{smart-form}%
\end{equation}
where $x=(x_{1},x_{2},\ldots,x_{n})$ and $\overline{x}_{k_{0}+3}=(x_{k_{0}%
+3},x_{k_{0}+4},\ldots,x_{n})$. Since $x(t_{0})$ belongs to $K_{i}$, we can
assume that $x(t_{0})=0$. Moreover, like in the proof of the previous Lemma,
we have $i+1\leq k_{0}\leq n-4$.

The Hamiltonian of this system is given by
\begin{align*}
H(x,p,u)  & =%
{\textstyle\sum_{j=1}^{i+1}}
p_{j}u_{j}+%
{\textstyle\sum_{j=i+2}^{k_{0}}}
p_{j}(x_{j-1}+c_{j-1})u_{i+2}+p_{k_{0}+1}x_{k_{0}}u_{i+2}+p_{k_{0}+2}u_{i+2}\\
& \ +\left(  p_{k_{0}+3}x_{k_{0}+2}+p_{k_{0}+4}x_{k_{0}+3}+%
{\textstyle\sum_{j=k_{0}+5}^{n}}
p_{j}\phi_{j}(\overline{x}_{k_{0}+3})\right)  x_{k_{0}}u_{i+2},
\end{align*}
which implies that any abnormal lift $(x(\cdot),p(\cdot))$ of $x(\cdot)$ must
satisfy
\begin{align}
\dot p_{j}  & =0\quad\quad\quad\quad\quad\quad\text{for }1\leq j\leq
i\label{co-adjoint}\\
\dot p_{j}  & =-p_{j+1}u_{i+2}\quad\quad\text{for }i+1\leq j\leq
k_{0}-1\nonumber\\
\dot p_{k_{0}}  & =\left(  -p_{k_{0}+1}-p_{k_{0}+3}x_{k_{0}+2}-p_{k_{0}%
+4}x_{k_{0}+3}-%
{\textstyle\sum_{j=k_{0}+5}^{n}}
p_{j}\phi_{j}(\overline{x}_{k_{0}+3})\right)  u_{i+2}\nonumber\\
\dot p_{k_{0}+1}  & =0\nonumber\\
\dot p_{j}  & =-\psi_{j}(\overline{x}_{k_{0}+3},p)x_{k_{0}}u_{i+2}\quad
\quad\text{for }k_{0}+2\leq j\leq n,\nonumber
\end{align}
where the $\psi_{j}$'s are some functions of $x_{k_{0}+3},x_{k_{0}+4}%
,\ldots,x_{n}$ and $p_{k_{0}+3},p_{k_{0}+4},\ldots,p_{n}$, for $k_{0}+2\leq
j\leq n$. The exact form of these functions is irrelevant for our purpose. Any
abnormal lift $(x(\cdot),p(\cdot))$ of $x(\cdot)$ must also satisfy the
relation $\tfrac{\partial H}{\partial u}=0$, which implies $p_{j}=0$, for
$1\leq j\leq i+1$, and
\begin{align*}
p_{k_{0}+2}  & =-%
{\textstyle\sum_{j=i+2}^{k_{0}}}
p_{j}(x_{j-1}+c_{j-1})+\\
& -\left(  p_{k_{0}+1}+p_{k_{0}+3}x_{k_{0}+2}+p_{k_{0}+4}x_{k_{0}+3}+%
{\textstyle\sum_{j=k_{0}+5}^{n}}
p_{j}\phi_{j}(\overline{x}_{k_{0}+3})\right)  x_{k_{0}}.
\end{align*}

Recall that $x(\cdot)$ is an integral curve of $\mathcal{D}^{(i)}$. Like in
the proof of Lemma~\ref{lem-abnormal-02}, if for a given $t$ the velocity
$\dot{x}(t)$ exists and satisfies~(\ref{smart-form}) then we have $\dot
{x}(t)\in\mathcal{A}^{(i)}(x(t))$ if and only if $u_{i+2}(t)=0$ or the three
following conditions hold: (i) $u_{i+1}(t)=0$ and (ii) $x_{j}(t)=0$, for
$i+1\leq j\leq k_{0}$, and (iii) $c_{j}=0$, for $i+1\leq j\leq k_{0}-1$. Now,
assume that$~x(\cdot)$ is an integral curve of $\mathcal{A}^{(i)}$, that is
$\dot{x}(t)\in\mathcal{A}^{(i)}(x(t))$ for almost all$~t$ in $I_{\varepsilon
}(t_{0})$. In order to prove that$~x(\cdot)$ is abnormal we must construct a
non-trivial abnormal lift $(x(\cdot),p(\cdot))$ of $x(\cdot)$. Take $p_{j}=0$,
for $1\leq j\leq n$, with the exception of$~p_{k_{0}+4}$, for which we take
any non-zero real constant. It is straightforward to check that our lift
satisfies (\ref{co-adjoint}). Indeed, the coordinate $x_{k_{0}+3}$ is constant
because we have $\dot{x}_{k_{0}+3}=x_{k_{0}}x_{k_{0}+2}u_{i+2}$; and
$u_{i+2}(t)=0$ or $x_{k_{0}}(t)=0$ for almost all$~t$. Moreover, since
$x(0)=0$, we have $x_{k_{0}+3}(t)=0$, for each$~t$ in $I_{\varepsilon}(t_{0}%
)$. It is also trivial to check that our lift satisfies $\tfrac{\partial
H}{\partial u}=0$. Since $p_{k_{0}+4}\neq0$ our lift in non-trivial. It
follows that the integral curve $x(\cdot)$ is abnormal.\hfill$\square$

\medskip\ 

\noindent\textbf{Proof of Theorem \ref{thm-abnormal-curves}} Let
$x:I_{\varepsilon}(t_{0})\rightarrow\mathbb{R}^{n}$ be the restriction to the
interval $I_{\varepsilon}(t_{0})$ of an integral curve of $\mathcal{D}^{(i)}$.
For$~\varepsilon>0$ small enough, we can apply both
Lemma~\ref{lem-abnormal-01} and Lemma~\ref{lem-abnormal-03}, which imply that
if the curve $x(\cdot)$ is such that $\dot x(t)$ belongs to $\mathcal{A}%
^{(i)}(x(t))$ for almost all $t$ in $I_{\varepsilon}(t_{0})$ then $x(\cdot)$
is abnormal. In other words, the integral curves of $\mathcal{A}^{(i)}$ are
locally abnormal.

Now assume that, for a fixed interval $I_{\varepsilon}(t_{0})$, the curve
$x(\cdot)$ is abnormal. Define the subset $I_{0}\subset I_{\varepsilon}%
(t_{0})$ by
\[
I_{0}=\{t\in I_{\varepsilon}(t_{0}):\dot{x}(t)\text{ exists and }\dot
{x}(t)\notin\mathcal{A}^{(i)}(x(t)\}
\]
We will show that if$~\varepsilon$ is small enough then the Lebesgue measure
of $I_{0}$ is zero. We can decompose $I_{0}$ into $I_{0}=I_{1}\cup I_{2}$,
where
\[
I_{1}=\{t\in I_{0}:x(t)\in K_{i}\}\text{\quad and\quad}I_{2}=\{t\in
I_{0}:x(t)\notin K_{i}\}.
\]
On the one hand, the measure of $I_{1}$ is equal to zero. Indeed, since for
each$~t$ in$~I_{0}$ we have $\dot{x}(t)\notin\mathcal{A}^{(i)}(x(t))$,
for$~\varepsilon$ small enough Lemma~\ref{lem-abnormal-02} implies that we
have $x(t)\notin K_{i}$ for almost all$~t$ in$~I_{0}$. But, on the other hand,
the measure of$~I_{2}$ is also equal to zero. To see this, let us
write$~I_{2}$ as $I_{2}=I_{0}\cap I_{3}$, where
\[
I_{3}=\{t\in I_{\varepsilon}(t_{0}):x(t)\notin K_{i}\}
\]
Since $K_{i}$ is closed (see the discussion following the statement of Theorem
\ref{thm-abnormal-curves}) and $x(\cdot)$ is continuous, it is clear that we
can decompose $I_{3}$ into a union of disjoint open intervals $I_{3}=%
{\textstyle\bigcup}
J_{\alpha}$ such that, on each of them, the curve $x(\cdot)$ has an empty
intersection with the singular locus $K_{i}$. Moreover, since the set $I_{3}$
is an open subset of $\mathbb{R}$ the union can be taken to be countable. Now
Lemma \ref{lem-abnormal-01} implies that, for each$~\alpha$, we have $\dot
{x}(t)\in\mathcal{A}^{(i)}(x(t))$ for almost all $t$ in $J_{\alpha}$ because
$x(\cdot)$ is abnormal and we are outside the singular locus. Hence, since the
measure of $I_{2}$ is the sum of the measures of the sets $I_{0}\cap
J_{\alpha}$ (the union is countable) and the measure of each of these sets is
zero, the measure of $I_{2}$ equals zero.\hfill$\square$

\subsection{Abnormal Curves and Singularity Type}

\begin{theorem}
\label{thm-singularity-type-abnormal}Let $\mathcal{D}$ and $\tilde
{\mathcal{D}}$ be two Goursat structures defined respectively on two
manifolds$~M$ and$~\tilde{M}$, both of dimension$~n$. Fix two points$~p$
and$~\widetilde{p}$ of$~M$ and$~\tilde{M}$, respectively. There exists a
diffeomorphism$~\varphi$, with $\tilde{p}=\varphi(p)$, between two small
enough neighborhoods of$~p$ and$~\tilde{p}$ that transforms, for $0\leq i\leq
n-4$, the abnormal curves of$~\mathcal{D}^{(i)}$ into the abnormal curves
of$~\tilde{\mathcal{D}}^{(i)} $ if and only if the singularity type
of$~\mathcal{D}$ at$~p$ equals the singularity type of$~\tilde{\mathcal{D}}$
at$~\tilde{p}$.
\end{theorem}

\noindent\textbf{Proof of Theorem \ref{thm-singularity-type-abnormal}}
\emph{Necessity:} Consider two distributions $\mathcal{D}$ and $\tilde
{\mathcal{D}}$, defined on two manifolds $M$ and $\tilde M$, respectively,
that have different singularity types $w$ and $\tilde w$ at $p$ and $\tilde p
$, respectively, that is $w=\delta_{\mathcal{D}}(p)$ and $w=\delta
_{\tilde{\mathcal{D}}}(\tilde p)$. We have already pointed out (see the Proof
of Theorem~\ref{thm-singularity-type-growth-vector}) that if $w$ and $\tilde
w$ are two words of the Jacquard language $J_{n}$ such that $w\neq\tilde w$
then there exists (after a permutation of $w$ and $\tilde w$, if necessary)
three words $z$, $v$, and $\tilde v$ such that both $w=vz$ and $\tilde
w=\tilde vz$, and which satisfy either
\[
\left\{
\begin{array}
[c]{lll}%
v & = & ua_{1}a_{2}\cdots a_{i-k}a_{0}^{k}\\
\tilde v & = & \tilde uc_{1}c_{2}\cdots c_{i},
\end{array}
\right.
\]
where $0\leq k\leq i-1$ and $c_{j}\neq a_{1}$ for $1\leq j\leq i$, or
\[
\left\{
\begin{array}
[c]{lll}%
v & = & ua_{1}a_{2}\cdots a_{i-k}a_{0}^{k}\\
\tilde v & = & \tilde ua_{1}a_{2}\cdots a_{i-l}a_{0}^{l},
\end{array}
\right.
\]
where $k<l$.

In both cases, consider the abnormal curves of $\mathcal{D}^{(i_{0}+k)}$ and
$\tilde{\mathcal{D}}^{(i_{0}+k)}$, where $i_{0}=\left|  z\right|  $. It
follows directly from the definition of the singularity type (see
Definition~\ref{def-singularity-type}) that for $\mathcal{D}$ we have $p\in
S_{i-k-1}^{(i_{0}+i-1)}$ while for $\tilde{\mathcal{D}}$ the point $\tilde{p}
$ does not belong to any submanifold $S_{j}^{(i_{0}+k+j)}$. Therefore, the
subset $\mathcal{A}^{(i_{0}+k)}(p)$ is not a linear subspace of $T_{p}M$ while
the subset $\tilde{\mathcal{A}}^{(i_{0}+k)}(\tilde{p})$ is a linear subspace
of $T_{\tilde{p}}\tilde{M}$. For each vector $\tau_{p}$ of $\mathcal{A}%
^{(i_{0}+k)}(p)$ there exist an abnormal curve of $\mathcal{D}^{(i_{0}+k)}$
that is tangent to $\tau_{p}$; for each vector $\tilde{\tau}_{p}$ of
$\tilde{\mathcal{A}}^{(i_{0}+k)}(\tilde{p})$ there exist an abnormal curve of
$\tilde{\mathcal{D}}^{(i_{0}+k)}$ that is tangent to$~\tilde{\tau}_{p}$. It
follows that no diffeomorphism can transform the abnormal curves of
$\mathcal{D}^{(i_{0}+k)}$ into the abnormal curves of $\tilde{\mathcal{D}%
}^{(i_{0}+k)}$, locally at$~p$ and$~\tilde{p}$.

\smallskip\ 

\noindent\emph{Sufficiency:} Now, assume that the singularity type
$\delta_{\mathcal{D}}(p)$ of $\mathcal{D}$ at $p$ and $\delta_{\tilde
{\mathcal{D}}}(\tilde p)$ of$~\tilde{\mathcal{D}}$ at$~\tilde p$ coincide and
are equal to $w$. The distribution $\mathcal{D}$ (respectively $\tilde
{\mathcal{D}}$) can be converted into a Kumpera-Ruiz normal form $\kappa^{n}$
(respectively $\tilde\kappa^{n}$) centered at$~p$ (respectively$~\tilde p$)
via a diffeomorphisms$~\phi$ (respectively $\tilde\phi$). Let $x=(x_{1}%
,\ldots,x_{n})$ (respectively $\tilde x=(\tilde x_{1},\ldots,\tilde x_{n})$)
denote the coordinates in which $\kappa^{n}$ (respectively $\tilde\kappa^{n}$)
is expressed. By Corollary~\ref{cor-singularity-type-kumpera-ruiz} and the
invariance of the singularity type we have $\delta_{\kappa^{n}}=\delta
_{\tilde\kappa^{n}}=w$. Moreover, by
Proposition~\ref{prop-singularity-type-claim}, a submanifold $S_{j}^{(i+j)}$
contains zero if and only if the submanifold $\tilde S_{j}^{(i+j)}$ contains
zero, which is the case if and only if $w=w_{1}a_{1}\cdots a_{j+1}w_{2}$, for
some words $w_{1}$ and $w_{2}$ such that $\left|  w_{2}\right|  =i$. If those
manifolds contain zero then, once again by
Proposition~\ref{prop-singularity-type-claim}, they are respectively given by
\[
S_{j}^{(i+j)}=\{x_{n-i-j}=0,\ldots,x_{n-i}=0\}\text{\quad and\quad}\tilde
S_{j}^{(i+j)}=\{\tilde x_{n-i-j}=0,\ldots,\tilde x_{n-i}=0\}\text{.}
\]

Now, for each integer $i$, we must distinguish two cases. \emph{First case:}
If for each integer $j$ the submanifolds $S_{j}^{(i+j)}$ and $\tilde{S}%
_{j}^{(i+j)}$ are empty, in a small enough neighborhood of zero, then, by
Theorem~\ref{thm-abnormal-curves}, the abnormal curves of $\mathcal{D}^{(i)}$
(respectively $\tilde{\mathcal{D}}^{(i)}$) are, in a small enough neighborhood
of zero, the integral curves of $\mathcal{C}_{i}$ (respectively$~\tilde
{\mathcal{C}}_{i}$). Moreover, we have
\[
\mathcal{C}_{i}=\left(
\begin{array}
[c]{c}%
\tfrac{\partial}{\partial x_{n}}%
\end{array}
,\ldots,
\begin{array}
[c]{c}%
\tfrac{\partial}{\partial x_{n-i}}%
\end{array}
\right)  \text{\quad and\quad}\tilde{\mathcal{C}}_{i}=\left(
\begin{array}
[c]{c}%
\tfrac{\partial}{\partial\tilde{x}_{n}}%
\end{array}
,\ldots,
\begin{array}
[c]{c}%
\tfrac{\partial}{\partial\tilde{x}_{n-i}}%
\end{array}
\right)  \text{.}
\]
\emph{Second case:} If for some integer $j$ the submanifolds $S_{j}^{(i+j)}$
and $\tilde{S}_{j}^{(i+j)}$ contain zero then, by
Proposition~\ref{prop-smooth-singularity}, this integer $j$ is unique. By
Theorem~\ref{thm-abnormal-curves}, the abnormal curves of $\mathcal{D}^{(i)}$
(respectively $\tilde{\mathcal{D}}^{(i)}$) are, in a small enough neighborhood
of zero, the integral curves of $\mathcal{A}_{j}^{(i)}$ (respectively
$\tilde{\mathcal{A}}_{j}^{(i)}$). Moreover, we have
\[
\mathcal{A}_{j}^{(i)}(q)=\left(
\begin{array}
[c]{c}%
\tfrac{\partial}{\partial x_{n}}%
\end{array}
,\ldots,
\begin{array}
[c]{c}%
\tfrac{\partial}{\partial x_{n-i}}%
\end{array}
\right)  (q)\text{\quad and\quad}\tilde{\mathcal{A}}_{j}^{(i)}(\tilde
{q})=\left(
\begin{array}
[c]{c}%
\tfrac{\partial}{\partial\tilde{x}_{n}}%
\end{array}
,\ldots,
\begin{array}
[c]{c}%
\tfrac{\partial}{\partial\tilde{x}_{n-i}}%
\end{array}
\right)  (\tilde{q})\text{,}
\]
for each point $q$ (respectively $\tilde{q}$) that does not belong to
$S_{j}^{(i+j)}$ (respectively $\tilde{S}_{j}^{(i+j)}$), and
\[
\mathcal{A}_{j}^{(i)}(q)=\left(
\begin{array}
[c]{c}%
\tfrac{\partial}{\partial x_{n}}%
\end{array}
,\ldots,
\begin{array}
[c]{c}%
\tfrac{\partial}{\partial x_{n-i+1}}%
\end{array}
,
\begin{array}
[c]{c}%
\tfrac{\partial}{\partial x_{n-i-j-1}}%
\end{array}
\right)  (q)\cup\left(
\begin{array}
[c]{c}%
\tfrac{\partial}{\partial x_{n}}%
\end{array}
,\ldots,
\begin{array}
[c]{c}%
\tfrac{\partial}{\partial x_{n-i}}%
\end{array}
\right)  (q)
\]
\[
\tilde{\mathcal{A}}_{j}^{(i)}(\tilde{q})=\left(
\begin{array}
[c]{c}%
\tfrac{\partial}{\partial\tilde{x}_{n}}%
\end{array}
,\ldots,
\begin{array}
[c]{c}%
\tfrac{\partial}{\partial\tilde{x}_{n-i+1}}%
\end{array}
,
\begin{array}
[c]{c}%
\tfrac{\partial}{\partial\tilde{x}_{n-i-j-1}}%
\end{array}
\right)  (\tilde{q})\cup\left(
\begin{array}
[c]{c}%
\tfrac{\partial}{\partial\tilde{x}_{n}}%
\end{array}
,\ldots,
\begin{array}
[c]{c}%
\tfrac{\partial}{\partial\tilde{x}_{n-i}}%
\end{array}
\right)  (\tilde{q}),
\]
for each point $q$ (respectively $\tilde{q}$) that belongs to $S_{j}^{(i+j)}$
(respectively $\tilde{S}_{j}^{(i+j)}$).

Let $\Phi$ be the local diffeomorphism of $\mathbb{R}^{n}$ defined by $\tilde
x_{i}=x_{i}$, for $1\leq i\leq n$. In both cases, the diffeomorphism
\[
\varphi=\phi^{-1}\circ\Phi\circ\phi
\]
transforms the integral curves of $\mathcal{A}^{(i)}$ into the integral curves
of $\tilde{\mathcal{A}}^{(i)}$, and thus, by Theorem~\ref{thm-abnormal-curves}%
, the abnormal curves of $\mathcal{D}^{(i)}$ into the abnormal curves
of$~\tilde{\mathcal{D}}^{(i)}$.\hfill$\square$

\subsection{Rigid Curves of Goursat Structures}

The concept of rigidity for integral curves of distributions was introduced by
Bryant and Hsu~\cite{bryant-hsu}. Rigid curves are always abnormal but there
exist abnormal curves that are not rigid (see e.g.
\cite{agrachev-sarychev-jmsec}, \cite{bryant-hsu}, and
\cite{zhitomirskii-nice}). Nevertheless, we will prove that in the case of
Goursat structures these two concepts coincide (for $C^{1}$ immersed curves).

\begin{definition}
Let $\mathcal{D}$ be a completely nonholonomic distribution defined on a
manifold $M$. Fix a closed interval $[a,b]$ and two points $p$ and $q$ in $M$.
Denote by $\mathcal{O}_{p,q}$ the space of all $C^{1}$ integral curves
$x:[a,b]\rightarrow M$ of $\mathcal{D}$ such that $x(a)=p$ and $x(b)=q$,
endowed with the $C^{1}$-topology. An integral curve $x(\cdot)$ that belongs
to $\mathcal{O}_{p,q}$ is \emph{rigid} if there exists a small enough
neighborhood $\mathcal{V}$ of $x(\cdot)$ in $\mathcal{O}_{p,q}$ such that any
curve $\widetilde{x}:[a,b]\rightarrow M$ contained in $\mathcal{V}$ is a
reparametrization of$~x(\cdot)$.
\end{definition}

Roughly speaking, a curve $x:[a,b]\rightarrow M$ is rigid if it is an isolated
point of $\mathcal{O}_{x(a),x(b)}$. Our study of abnormal curves leads easily
to the following result, which characterizes immersed rigid curves. This
result gives also, for Goursat structures, a more intuitive view of the
concept of abnormal curve.

Let $I\subset\mathbb{R}$ be a \emph{closed} interval. For any $t_{0}\in I$ and
for any $\varepsilon>0$, denote by $I_{\varepsilon}(t_{0})$ the intersection
$I\cap[t_{0}-\varepsilon,t_{0}+\varepsilon]$. An integral curve
$x:I\rightarrow M$ of $\mathcal{D}$ is \emph{locally rigid} if for each
$t_{0}$ in $I$ there exists a small enough $\varepsilon>0$ such that the
restriction of $x(\cdot)$ to $I_{\varepsilon}(t_{0})$ is rigid.

\begin{theorem}
\label{thm-rigidity}Let $x(\cdot)$ be a$~C^{1}$ immersed integral curve of a
Goursat structure $\mathcal{D}$, defined on a manifold of dimension$~n$. The
three following conditions are equivalent:

\begin{enumerate}
\item  The curve $x(\cdot)$ is locally abnormal;

\item  The curve $x(\cdot)$ is locally rigid;

\item  The curve $x(\cdot)$ is either an integral curve of $\mathcal{C}_{0}$
or an integral curve of $\mathcal{A}_{k_{0}-1}^{(0)}$, for some $1\leq
k_{0}\leq n-4$.
\end{enumerate}
\end{theorem}

We supposed in this Theorem that the integral curve is \emph{immersed}, which
means that its velocity (defined everywhere, since the curve is $C^{1}$) never
vanish. This assumption is fundamental. Indeed, an immersed rigid curve can
loose its rigidity if we change its parametrization in such a way that it is
not immersed anymore (see e.g. \cite{zelenko-zhitomirskii}). Observe also that
the Theorem is stated for integral curves of $\mathcal{D}^{(0)}$ but not for
those of $\mathcal{D}^{(i)}$, if $i\geq1$. In fact, the abnormal curves of
$\mathcal{D}^{(i)}$ such that their velocity does not belong to $\mathcal{C}%
_{i-1}$ have only a weaker form of rigidity: all curves that are close enough
to them in the $C^{1}$ topology stay in a submanifold of the original
manifold. We will consider this situation in a forthcoming work.

Our proof of Theorem~\ref{thm-rigidity} is mainly based on the ideas
introduced by Bryant and Hsu~\cite{bryant-hsu} and Zhitomirski\u
\i\ \cite{zhitomirskii-nice}. In particular, it is a direct consequence of
Zhitomirski\u\i's work that the immersed integral curves of $\mathcal{C}_{0}$
are rigid.

To prove the rigidity of the integral curves of $\mathcal{A}_{k_{0}-1}^{(0)}$
we follow the main ideas of~\cite{zhitomirskii-nice}. Note, however, that the
statement for $\mathcal{A}_{k_{0}-1}^{(0)}$ is not implied by any of the
results of \cite{agrachev-sarychev-jmsec}, \cite{bryant-hsu},
\cite{sussmann-liu}, or \cite{zhitomirskii-nice} because Goursat structures
are highly non-generic and do not fit into the large categories of (generic)
rank two distributions studied in those papers. We would like to point out
that, in the particular case of dimension five, the rigidity of the immersed
integral curves of $\mathcal{A}_{0}^{(0)}$ was already observed
in~\cite{mormul-fake}.\ Moreover, the equivalence of Items (ii) and (iiii) of
Theorem~\ref{thm-rigidity} has already been announced in
\cite{pasillas-respondek-nolcos}.

Our proof of Theorem~\ref{thm-rigidity} will use the following Lemma, which
will be proved later in Appendix \ref{sec-additional}. The normal form that we
introduce in it is analogous to the one used in \cite{zhitomirskii-nice} to
prove that the integral curves of $\mathcal{C}_{0}$ are rigid.

\begin{lemma}
\label{lem-rigid}Let $\mathcal{D}$ be a Goursat structure on a manifold $M$ of
dimension $n\geq5$. If the singularity type of $\mathcal{D}$ at $p$ is equal
to $\delta_{\mathcal{D}}(p)=wa_{1}\cdots a_{k_{0}}$ for some $1\leq k_{0}\leq
n-4$, where $w$ is an arbitrary word of $J_{n-k_{0}-4}$, then $\mathcal{D}$ is
locally equivalent at $p$ to the distribution spanned by a pair of vector
fields that has the following form:
\begin{align*}
\xi_{1}  & =\tfrac{\partial}{\partial y_{1}}\\
\xi_{2}  & =y_{1}\tfrac{\partial}{\partial y_{2}}+\cdots+y_{k_{0}}%
\tfrac{\partial}{\partial y_{k_{0}+1}}+\tfrac{\partial}{\partial y_{k_{0}+2}%
}+\tfrac{1}{2}y_{k_{0}+1}^{2}\tfrac{\partial}{\partial y_{k_{0}+3}}+%
{\textstyle\sum\limits_{i=k_{0}+4}^{n}}
\varphi_{i}(y)\tfrac{\partial}{\partial y_{i}},
\end{align*}
where the coordinates $y_{1},\ldots,y_{n}$ are centered at $p$. In these
coordinates, the canonical submanifold $S_{k_{0}-1}^{(k_{0}-1)}$ is given by
\[
S_{k_{0}-1}^{(k_{0}-1)}=\{y_{1}=0,\ldots,y_{k_{0}}=0\}.
\]
Moreover, we have $\mathcal{C}_{0}=(\xi_{1})$, for any point $p$ of
$\mathbb{R}^{n}$, and $\mathcal{A}_{k_{0}-1}^{(0)}(p)=(\xi_{2})(p)$, for any
point$~p$ of $S_{k_{0}-1}^{(k_{0}-1)}$.\ 
\end{lemma}

\noindent\textbf{Proof of Theorem \ref{thm-rigidity}} It is well known that
rigidity implies abnormality (see \cite{agrachev-sarychev-jmsec},
\cite{bryant-hsu}, and \cite{zhitomirskii-nice}) and thus, that (ii) implies
(i). By Theorem~\ref{thm-abnormal-curves}, any abnormal curve of
$\mathcal{D}^{(0)}$ is an integral curve of $\mathcal{A}^{(0)}$. Recall that
$\mathcal{A}^{(0)}(p)=\mathcal{C}_{0}(p)\cup\mathcal{A}_{k_{0}-1}^{(0)}(p)$,
for a unique $1\leq k_{0}\leq n-4$, and that $\mathcal{C}_{0}(p)\cap
\mathcal{A}_{k_{0}-1}^{(0)}(p)=0$. Therefore any $C^{1}$ immersed abnormal
curve of $\mathcal{D}^{(0)}$ is either an integral curve of $\mathcal{C}_{0}$
or an integral curve of $\mathcal{A}_{k_{0}-1}^{(0)}$. Hence (i) implies (iii).

What remains to prove is that if a $C^{1}$ immersed integral curve
$y:I_{\varepsilon}(t_{0})\rightarrow\mathbb{R}^{n}$ of $\mathcal{D}^{(0)}$ is
an integral curve of either $\mathcal{C}_{0}$ or $\mathcal{A}_{k_{0}-1}^{(0)}$
then it is rigid. This result is known~\cite{zhitomirskii-nice} for the
integral curves of $\mathcal{C}_{0}$. We can thus assume that $y(\cdot)$ is an
immersed integral curve of $\mathcal{A}_{k_{0}-1}^{(0)}$ (which then, by
definition, stays in $S_{k_{0}-1}^{(k_{0}-1)}$). It follows from Lemma
\ref{lem-rigid} that we can find coordinates such that $y(\cdot)$ satisfies
$y(t_{0}-\varepsilon)=0$ and is a solution of the following control system:
\begin{align}
\dot y_{1}  & =u_{1}\nonumber\\
\dot y_{2}  & =y_{1}u_{2}\nonumber\\
& \ \vdots\nonumber\\
\dot y_{k_{0}+1}  & =y_{k_{0}}u_{2}\label{rigid-control-system}\\
\dot y_{k_{0}+2}  & =u_{2}\nonumber\\
\dot y_{k_{0}+3}  & =\tfrac12y_{k_{0}+1}^{2}u_{2}\nonumber\\
\dot y_{i}  & =\varphi_{i}(y)u_{2}\text{\quad for }k_{0}+4\leq i\leq
n\text{,}\nonumber
\end{align}
with $u_{1}(t)=0$ (because $y_{1}(t)=0$ on $S_{k_{0}-1}^{(k_{0}-1)}$) and
$u_{2}(t)\neq0$ for each $t$ in $I_{\varepsilon}(t_{0})$ (because the curve is
immersed). Since the coordinates of Lemma~\ref{lem-rigid} are chosen to be
centered at $y(t_{0}-\varepsilon)=0$, from $y_{k_{0}}(t)=0$ we conclude that
$y_{k_{0}+1}(t)=0$, and thus that $y_{k_{0}+3}(t)=0$, for each $t$ in
$I_{\varepsilon}(t_{0})$.

Now, consider a $C^{1}$ immersed integral curve $\tilde{y}:I_{\varepsilon
}(t_{0})\rightarrow\mathbb{R}^{n}$ of $\mathcal{D}^{(0)}$ that has the same
end-points as the curve $y(\cdot)$. In particular, we have $\tilde{y}%
_{k_{0}+3}(t_{0}-\varepsilon)=0$ and $\tilde{y}_{k_{0}+3}(t_{0}+\varepsilon
)=0$. By taking a small enough neighborhood of $y(\cdot)$ in $\mathcal{O}%
_{y(t_{0}-\varepsilon),y(t_{0}+\varepsilon)}$ (which is not the same as taking
a smaller $\varepsilon>0$), we can assume that $\tilde{u}_{2}(t)\neq0$ for
each $t$ in $I_{\varepsilon}(t_{0})$, where $\tilde{u}_{1}$ and $\tilde{u}%
_{2}$ denote the controls for which $\tilde{y}(\cdot)$ is a solution of
(\ref{rigid-control-system}). Without loss of generality, we can assume that
$\tilde{u}_{2}(t)>0$ (the proof for $\tilde{u}_{2}(t)<0$ is identical). Since
we have $\tfrac{1}{2}\tilde{y}_{k_{0}+1}^{2}\tilde{u}_{2}(t)\geq0$ for each
$t$ in $I_{\varepsilon}(t_{0})$ and both $\tilde{y}_{k_{0}+3}(t_{0}%
-\varepsilon)=0$ and $\tilde{y}_{k_{0}+3}(t_{0}+\varepsilon)=0$, we must have
$\tilde{y}_{k_{0}+1}^{2}(t)=0$, for each $t$ in $I_{\varepsilon}(t_{0})$.
Together with $\tilde{u}_{2}>0$, the latter relation implies that, for $1\leq
i\leq k_{0}+1$, we have $\tilde{y}_{i}(t)=0 $ for each $t$ in $I_{\varepsilon
}(t_{0})$, which clearly implies $\tilde{u}_{1}(t)=0$ for each $t$ in
$I_{\varepsilon}(t_{0})$. Hence, the curve $\tilde{y}(\cdot)$ is a
reparametrization of the original curve $y(\cdot)$. Indeed, these two curves
are $C^{1}$ immersed integral curves of $\xi_{2}$ and have the same end points
(see \cite{zhitomirskii-nice} for more details about this last point).\hfill
$\square$

\subsection{Rigid Curves of the N-Trailer System}

Let us illustrate Theorem \ref{thm-rigidity} by applying it to the $n$-trailer
system. Let $\mathcal{D}$ be the Goursat structure spanned by the $n $-trailer
system $\tau^{n}$ on $\mathbb{R}^{2}\times(S^{1})^{n+1}$. Recall that Jean's
sequence of sets of real numbers $\alpha_{i}$, for $i\geq0$, is defined by the
relations
\[%
\begin{array}
[c]{lcl}%
\alpha_{1} & = & \{-\frac\pi2,+\frac\pi2\}\\
\alpha_{i+1} & = & \{\arctan\sin(\alpha),\text{ }\arctan\sin(\alpha
)+\pi:\text{ }\alpha\in\alpha_{i}\}.
\end{array}
\]
By Proposition~\ref{prop-singularity-type-trailer}, we have
\[
S_{j}^{(i)}=\{\theta_{n-i}-\theta_{n-i-1}\in\alpha_{1},\ldots,\theta
_{n-i+j}-\theta_{n-i+j-1}\in\alpha_{j+1}\},
\]
for $0\leq i\leq n-2$ and $0\leq j\leq i$ (recall that $n$ is the number of
trailers, not the dimension of the configuration space!). It obviously follows
that
\[
S_{j}^{(j)}=\{\theta_{n-j}-\theta_{n-j-1}\in\alpha_{1},\ldots,\theta
_{n}-\theta_{n-1}\in\alpha_{j+1}\},
\]
for $0\leq j\leq n-2$. Each submanifold $S_{j}^{(j)}$ has clearly codimension
$j+1$. Moreover, these manifolds are pairwise disjoint. Thus a given point $p$
is either in none of the submanifold $S_{j}^{(j)}$ at all or in one and only
one of them.

Recall also that $\mathcal{A}_{j}^{(0)}(p)=\mathcal{D}^{(i)}(p)\cap T_{p}%
S_{j}^{(j)}$ and that $\mathcal{A}^{(0)}(p)=\mathcal{C}_{0}(p)\cup
\mathcal{A}_{j}^{(0)}(p)$, for a unique $0\leq j\leq n-2$. The canonical line
field $\mathcal{C}_{0}$ is given on $\mathbb{R}^{2}\times(S^{1})^{n+1}$ by
$(\tfrac\partial{\partial\theta_{n}})$. A simple computation shows that, on
each submanifold $S_{j}^{(j)}$, the line field $\mathcal{A}_{j}^{(0)}$ is
given by $\mathcal{A}_{j}^{(0)}=(\tfrac\partial{\partial\theta_{n}}%
+\cdots+\tfrac\partial{\partial\theta_{n-j-1}})$.

By Theorem \ref{thm-rigidity}, a $C^{1}$ motion of the $n$-trailer for which
the velocity never vanishes is rigid if and only if (i) it is an integral
curve of $\mathcal{C}_{0}$ or (ii) it is an integral curve of $\mathcal{A}%
_{j}^{(0)}$. In the second case, the motion lies in $S_{j}^{(j)}$. In fact,
there is an easy way to visualize these rigid trajectories:

\begin{corollary}
An immersed motion of the $n$-trailer system is locally rigid if and only if
it fixes the positions in the $(\xi_{1},\xi_{2})$-plane of the centers of the
axles of at least two trailers.
\end{corollary}

For example, there passes through any configuration of $\mathbb{R}^{2}%
\times(S^{1})^{n+1}$ an integral curve of $\mathcal{C}_{0}$. The corresponding
motion fixes the positions in the $(\xi_{1},\xi_{2})$-plane of all trailers
(we just turn the front wheels). If a configuration is such that $\theta
_{n}-\theta_{n-1}\in\alpha_{1}$ (it belongs to $S_{0}^{(0)}$) then, besides
the motions associated to $\mathcal{C}_{0}$, there is an additional motion
given by $\mathcal{A}_{0}^{(0)}$ for which the positions in the $(\xi_{1}%
,\xi_{2})$-plane of all trailers, excepted the first one, are fixed. For these
motions, the center of the first trailer moves on a circle around the center
of the second trailer, which turns with its center fixed (see e.g.
Figure~\ref{fig-rigid-trailer-R5}). Observe that such a motion is possible if
and only if $\theta_{n}-\theta_{n-1}=\pm\pi/2$.

\section{Contact Transformations}

\label{sec-contact}

\subsection{A Singular Version of B\"acklund's Theorem}

Let $\mathcal{D}$ and $\tilde{\mathcal{D}}$ be two Goursat structures defined
on two manifolds $M$ and $\tilde{M}$, respectively, of dimension $n\geq3$. A
\emph{(generalized) contact transformation} (of order $n-2$) is a smooth
diffeomorphism $\phi$ between $M$ and $\tilde{M}$ such that $(\phi_{\ast
}\mathcal{D})(\tilde{p})=\tilde{\mathcal{D}}(\tilde{p})$, for each point
$\tilde{p}$ in $\tilde{M}$. Such transformations are called
\emph{automorphisms} in the work of Kumpera and Ruiz \cite{kumpera-ruiz} (see
also \cite{gaspar} and \cite{mormul-dijon}). In a neighborhood of a regular
point our definition coincides with the classical definition of a contact
transformation on the space $J^{n-2}(\mathbb{R},\mathbb{R})$ of $(n-2)$-jets
of functions that have one dependent and one independent variable (see
\cite{bryant-chern-gardner-goldschmidt-griffiths} and \cite{olver-equivalence}%
). From now on, unless we want to distinguish generalized contact
transformations from the classical ones, we will omit the word ``generalized''.

Fix two points $p$ and $\tilde{p}$ of $M$ and $\tilde{M}$, respectively. Let
$\phi$ be a local contact transformation between$~\mathcal{D}$ and$~\tilde
{\mathcal{D}}$ such that $\phi(p)=\tilde{p}$. Fix two small enough
neighborhoods$~U$ and$~\tilde{U}$ of$~p$ and$~\tilde{p}$, respectively, such
that $\tilde{U}=\phi(U)$ and such that$~\mathcal{D}$ on$~U$ and$~\tilde
{\mathcal{D}}$ on$~\tilde{U}$ are equivalent to two Kumpera-Ruiz normal
forms$~\kappa^{n}$ and$~\tilde{\kappa}^{n}$ centered at$~p$ and$~\tilde{p}$,
respectively, and defined on two open subsets $x(U)$ and $\tilde{x}(\tilde
{U})$ of$~\mathbb{R}^{n}$, where$~x$ and$~\tilde{x}$ denote coordinates that
transform the Goursat structures$~\mathcal{D}$ and$~\tilde{\mathcal{D}}$ into
their Kumpera-Ruiz normal forms$~\kappa^{n}$ and$~\tilde{\kappa}^{n}$,
respectively. We can assume, without loss of generality, that the first
prolongation (in the sequence of prolongations that define$~\kappa^{n}$
and$~\tilde{\kappa}^{n}$) is regular. Namely $\kappa^{4}=R_{0}(\kappa^{3})$
and $\tilde{\kappa}^{4}=R_{0}(\tilde{\kappa}^{3})$. Once such a pair of
Kumpera-Ruiz charts $(x,U)$ and $(\tilde{x},\tilde{U})$ has been fixed, we can
associate to the contact transformation$~\phi$ a unique contact
transformation$~\Phi$, between$~\kappa^{n}$ on$~x(U)$ and$~\tilde{\kappa}^{n}$
on$~\tilde{x}(\tilde{U})$, by taking
\[
\Phi=\tilde{x}\circ\phi\circ x^{-1}.
\]
In other words $\tilde{x}=(\Phi\circ x)\circ\psi$, where$~\psi$ denotes the
inverse of the diffeomorphism$~\phi$. Observe that, since the Kumpera-Ruiz
charts$~x$ and$~\tilde{x}$ are centered at$~p$ and$~\tilde{p}$, respectively,
we have $\Phi(0)=0$. We will denote by$~\Phi_{i}$ the$~i^{\mathrm{th}}$
component of$~\Phi$.\ 

In the next two Propositions we will assume that all the above defined data
(the Goursat structures$~\mathcal{D}$ and$~\tilde{\mathcal{D}}$, the
diffeomorphism$~\phi$, the coordinates$~x$ and$~\tilde{x}$, and the
Kumpera-Ruiz normal forms$~\kappa^{n}$ and$~\tilde{\kappa}^{n}$) have been
fixed and, therefore, that the diffeomorphism$~\Phi$ is uniquely defined. The
following result is a direct consequence of the obvious relations
\[
\phi_{\ast}(\mathcal{C}_{i})=\tilde{\mathcal{C}}_{i},
\]
for $0\leq i\leq n-4$, where $\mathcal{C}_{i}\subset\mathcal{D}^{(i)}$ denotes
the characteristic distribution of $\mathcal{D}^{(i+1)}$ and $\tilde
{\mathcal{C}}_{i}$ that of $\tilde{\mathcal{D}}^{(i)}$ (see
Proposition~\ref{prop-cartan}).

\begin{proposition}
\label{prop-automorphisms-01}For each $1\leq i\leq3$ we have $\Phi_{i}%
(x)=\Phi_{i}(x_{1},x_{2},x_{3})$. For each $4\leq i\leq n$ we have $\Phi
_{i}(x)=\Phi_{i}(x_{1},\ldots,x_{i})$.
\end{proposition}

Decompose$~\mathbb{R}^{n}$ into a direct product $\mathbb{R}^{n}%
=\mathbb{R}^{i}\times\mathbb{R}^{n-i}$. It follows directly from
Proposition~\ref{prop-automorphisms-01} that for each $3\leq i\leq n$ we can
build a diffeomorphism$~\Phi^{(i)}$, between the projection of $x(U)$
on$~\mathbb{R}^{i}$ and the projection of $\tilde x(\tilde U)$ on$~\mathbb{R}%
^{i}$, by taking the components$~\Phi_{j}$, for $1\leq j\leq i$, as the
components of $\Phi^{(i)}$. Denote by $\Psi^{(i)}$ the inverse of $\Phi^{(i)}%
$. We obviously have $\Phi^{(n)}=\Phi$. The following result is a direct
consequence of Proposition~\ref{prop-automorphisms-01} and the obvious
relations
\[
\phi_{*}(\mathcal{D}^{(i)})=\tilde{\mathcal{D}}^{(i)},
\]
which hold for $0\leq i\leq n-2$. Recall that, by definition, the two
Kumpera-Ruiz normal forms$~\kappa^{n}$ and$~\tilde\kappa^{n}$ are given by two
sequences of prolongations. We will denote by $\kappa^{3},\ldots,\kappa^{n}$
and $\tilde\kappa^{3},\ldots,\tilde\kappa^{n}$, respectively, the Kumpera-Ruiz
normal forms obtained as intermediate steps of these successive prolongations.

\begin{proposition}
\label{prop-automorphisms-02}There exist four smooth functions, denoted by
$\nu_{3}$, $\eta_{3}$, $\mu_{3}$, and $\lambda_{3}$, that depend on the
coordinates $x_{1}$, $x_{2}$, and $x_{3}$ only, such that
\[%
\begin{array}
[c]{ccl}%
\Phi_{\ast}^{(3)}(\kappa_{1}^{3}) & = & (\nu_{3}\circ\Psi^{(3)})\tilde{\kappa
}_{1}^{3}+(\lambda_{3}\circ\Psi^{(3)})\tilde{\kappa}_{2}^{3}\\
\Phi_{\ast}^{(3)}(\kappa_{2}^{3}) & = & (\eta_{3}\circ\Psi^{(3)})\tilde
{\kappa}_{1}^{3}+(\mu_{3}\circ\Psi^{(3)})\tilde{\kappa}_{2}^{3}.
\end{array}
\]
Moreover, for each $i\geq4$, there exist three smooth functions, denoted by
$\nu_{i}$, $\eta_{i}$, and $\mu_{i}$, that depend on the coordinates
$x_{1},\ldots,x_{i}$ only, such that
\[%
\begin{array}
[c]{ccl}%
\Phi_{\ast}^{(i)}(\kappa_{1}^{i}) & = & (\nu_{i}\circ\Psi^{(i)})\tilde{\kappa
}_{1}^{i}\\
\Phi_{\ast}^{(i)}(\kappa_{2}^{i}) & = & (\eta_{i}\circ\Psi^{(i)})\tilde
{\kappa}_{1}^{i}+(\mu_{i}\circ\Psi^{(i)})\tilde{\kappa}_{2}^{i}.
\end{array}
\]
The functions $\nu_{i}$, $\eta_{i}$, $\mu_{i}$, and $\lambda_{i}$ are uniquely
defined, for each $i\geq3$, once the diffeomorphism $\Phi$ has been fixed.
They obviously satisfy $(\nu_{3}\mu_{3}-\lambda_{3}\eta_{3})(0)\neq0$ and
$(\nu_{i}\mu_{i})(0)\neq0$, for $i\geq4$.
\end{proposition}

The following result can be considered as a singular version of B\"{a}cklund's
theorem~\cite{backlund} (see \cite{olver-equivalence} for a modern approach).
It shows that any contact transformation is the ``prolongation'' of a first
order contact transformation. Though the case $n=4$ is classical~
\cite{bryant-hsu}, it seems that our result for $n\geq5$ is new. Notice that a
weaker version of Theorem~\ref{thm-automorphisms} has already been announced
in~\cite{cheaito-mormul-pasillas-respondek}. Independently, an infinitesimal
version of Theorem~\ref{thm-automorphisms} has been announced
in~\cite{mormul-poisson} and proved in \cite{mormul-dijon}.

\begin{theorem}
\label{thm-automorphisms}Let $\phi$ be a local (generalized) contact
transformation between two Goursat structures$~\mathcal{D}$ and$~\tilde
{\mathcal{D}}$, defined locally at$~p$ and$~\tilde{p}$, respectively. Let$~x$
and$~\tilde{x}$ be local coordinates that transform$~\mathcal{D}$
and$~\tilde{\mathcal{D}}$ into their Kumpera-Ruiz normal forms$~\kappa^{n}$
and$~\tilde{\kappa}^{n}$, respectively, and let $\delta_{\mathcal{D}}%
(p)=w_{0}\cdots w_{n-4}$ be the singularity type of$~\mathcal{D}$ at$~p$,
which equals $\delta_{\tilde{\mathcal{D}}}(\tilde{p})$ since$~\mathcal{D}$
at$~p$ and$~\tilde{\mathcal{D}}$ at$~\tilde{p}$ are locally equivalent. The
constants$~c_{i}$ and$~\tilde{c}_{i}$ that appear in$~\kappa^{n}$
and$~\tilde{\kappa}^{n}$, respectively, and the contact transformation$~\Phi$
associated to$~\phi$ and to the coordinates$~x$ and$~\tilde{x}$ fulfill the
following relations:

\begin{enumerate}
\item  The diffeomorphism $\Phi^{(3)}$ is a first order contact transformation
and the functions $\nu_{3}$, $\eta_{3}$, $\mu_{3}$, and $\lambda_{3}$ are
uniquely determined by $\Phi^{(3)}$.

\item  The diffeomorphism $\Phi^{(4)}$ is uniquely defined by
\[%
\begin{array}
[c]{ccl}%
\Phi_{4}(x) & = & \dfrac{\nu_{3}+x_{4}\eta_{3}}{\mu_{3}+x_{4}\lambda_{3}}\\
&  & \\
\mu_{4} & = & \mu_{3}+x_{4}\lambda_{3}\\
\nu_{4} & = & \mathrm{L}_{\kappa_{1}^{4}}\Phi_{4}=(\mu_{3}\eta_{3}-\lambda
_{3}\nu_{3})/(\nu_{3}+x_{4}\eta_{3})^{2}\\
\eta_{4} & = & \mathrm{L}_{\kappa_{2}^{4}}\Phi_{4}.
\end{array}
\]

\item  If $i\geq5$ and $w_{i-4}\neq a_{1}$ then $\Phi^{(i)}$ is uniquely
defined by
\[%
\begin{array}
[c]{ccl}%
\tilde{c}_{i} & = & c_{i}\dfrac{\nu_{i-1}(0)}{\mu_{i-1}(0)}+\dfrac{\eta
_{i-1}(0)}{\mu_{i-1}(0)}\\
\Phi_{i}(x) & = & \dfrac{1}{\mu_{i-1}}\left(  (x_{i}+c_{i})\nu_{i-1}%
+\eta_{i-1}\right)  -\tilde{c}_{i}\\
&  & \\
\mu_{i} & = & \mu_{i-1}\\
\nu_{i} & = & \mathrm{L}_{\kappa_{1}^{i}}\Phi_{i}=\nu_{i-1}/\mu_{i-1}\\
\eta_{i} & = & \mathrm{L}_{\kappa_{2}^{i}}\Phi_{i}.
\end{array}
\]

\item  If $i\geq5$ and $w_{i-4}=a_{1}$ then $\Phi^{(i)}$ is uniquely defined
by
\[%
\begin{array}
[c]{ccl}%
\Phi_{i}(x) & = & \dfrac{x_{i}\mu_{i-1}}{\nu_{i-1}+x_{i}\eta_{i-1}}\\
&  & \\
\mu_{i} & = & \nu_{i-1}+x_{i}\eta_{i-1}\\
\nu_{i} & = & \mathrm{L}_{\kappa_{1}^{i}}\Phi_{i}=(\mu_{i-1}\nu_{i-1}%
)/(\nu_{i-1}+x_{i}\eta_{i-1})^{2}\\
\eta_{i} & = & \mathrm{L}_{\kappa_{2}^{i}}\Phi_{i}.
\end{array}
\]
\end{enumerate}

Therefore, the (generalized) contact transformation$~\Phi$ is uniquely
determined by the first order contact transformation$~\Phi^{(3)}$.
\end{theorem}

This Theorem says that any (generalized) contact transformation between two
Goursat structures is uniquely defined by a first order contact transformation
$\Phi^{(3)}$. In fact, the component $\Phi_{4}$ of $\Phi^{(4)}$ is a linear
fractional transformation (M\"{o}bius transformation) whose coefficients are
uniquely determined by the components of $\Phi^{(3)}$ (compare
\cite{bryant-hsu}). For $i\geq5$, successively, the component $\Phi_{i}$ of
$\Phi^{(i)}$ is either, in the case of a singular prolongation, a zero
preserving linear fractional transformation with $x_{i}=0$ being fixed by the
fact that $\Phi^{(i)}$ preserves the hypersurface $\{x_{i}=0\}$ or, in the
case of a regular prolongation, by an affine transformation. In both cases,
the coefficients of the linear fractional transformation or of the affine
transformation are uniquely determined by $\Phi^{(i-1)}$.

\medskip\ 

\noindent\textbf{Proof of Theorem \ref{thm-automorphisms}} If $n=3$ then there
is nothing to prove. If $n=4$ then the result is well known (see e.g.
\cite{bryant-hsu}). Therefore, we can proceed by induction on the integer
$n\geq5$. Assume that the Theorem is true for $n-1$. By
Proposition~\ref{prop-automorphisms-02}, we have
\[%
\begin{array}
[c]{ccl}%
\Phi_{\ast}^{(n-1)}(\kappa_{1}^{n-1}) & = & (\nu_{n-1}\circ\Psi^{(n-1)}%
)\tilde{\kappa}_{1}^{n-1}\\
\Phi_{\ast}^{(n-1)}(\kappa_{2}^{n-1}) & = & (\eta_{n-1}\circ\Psi
^{(n-1)})\tilde{\kappa}_{1}^{n-1}+(\mu_{n-1}\circ\Psi^{(n-1)})\tilde{\kappa
}_{2}^{n-1}.
\end{array}
\]
In other words, the restriction of $\Phi$ to $x(U)\cap\mathbb{R}^{n-1}$,
equipped with coordinates $x_{1},\ldots,x_{n-1}$, is a contact transformation
between $\kappa^{n-1}$ and $\tilde{\kappa}^{n-1}$. Since the Theorem is
assumed to be true for $n-1$, each component $\Phi_{i}$, for $1\leq i\leq
n-1$, satisfies the relations given by the Theorem, as do, for $3\leq i\leq
n-1$, the smooth functions $\nu_{i}$, $\mu_{i}$, $\lambda_{i}$, and $\eta_{i}%
$, given by Proposition~\ref{prop-automorphisms-02}. What remains to check is
that $\Phi_{n}$, $\nu_{n}$, $\mu_{n}$, and $\eta_{n}$ satisfy our conditions.

Recall that for any diffeomorphism $\Phi^{(n)}=(\Phi^{(n-1)},\Phi_{n})^{\top}$
of $\mathbb{R}^{n}$, such that $\Phi^{(n-1)}$ depends on the first $n-1$
coordinates $x_{1},\ldots,x_{n-1}$ only, and for any vector field $f=\alpha
f^{n-1}+f_{n}$ on~$\mathbb{R}^{n}$, where $\alpha$ is a smooth function on
$\mathbb{R}^{n}$, the vector field $f^{n-1}$ is the lift of a vector field
on~$\mathbb{R}^{n-1}$, and the only non-zero component of $f_{n}$ is the last
one, we have:
\begin{equation}
\Phi_{\ast}^{(n)}(f)=(\alpha\circ\Psi^{(n)})\Phi_{\ast}^{(n-1)}(f^{n-1}%
)+\left(  (\mathrm{L}_{f}\Phi_{n})\circ\Psi^{(n)}\right)  \tfrac{\,\partial
}{\partial\tilde{x}_{n}}.\label{triangular-tangent-map}%
\end{equation}
Observe that the vector field $\Phi_{\ast}^{(n-1)}(f^{n-1})$ is lifted (see
Notation~\ref{not-lift}) along the coordinate~$\tilde{x}_{n}$, which is given
by$~\Phi_{n}$.

\medskip\ 

\noindent\emph{Regular case\/}: If $w_{n-4}\neq a_{1}$ then we have
$\kappa_{2}^{n}=(x_{n}+c_{n})\kappa_{1}^{n-1}+\kappa_{2}^{n-1}$. This
relation, together with (\ref{triangular-tangent-map}) and the induction
hypothesis leads to:
\begin{align*}
\Phi_{*}^{(n)}(\kappa_{2}^{n})  & =\left(  (x_{n}+c_{n})\circ\Psi
^{(n)}\right)  \Phi_{*}^{(n-1)}(\kappa_{1}^{n-1})+\Phi_{*}^{(n-1)}(\kappa
_{2}^{n-1})\\
& +\left(  (\mathrm{L}_{\kappa_{2}^{n}}\Phi_{n})\circ\Psi^{(n)}\right)
\,\tilde\kappa_{1}^{n}\\
& =\left(  \left(  (x_{n}+c_{n})\nu_{n-1}+\eta_{n-1}\right)  \circ\Psi
^{(n)}\right)  \tilde\kappa_{1}^{n-1}+\left(  \mu_{n-1}\circ\Psi^{(n)}\right)
\tilde\kappa_{2}^{n-1}\\
& +\left(  (\mathrm{L}_{\kappa_{2}^{n}}\Phi_{n})\circ\Psi^{(n)}\right)
\,\tilde\kappa_{1}^{n}\\
& =\left(  \mu_{n-1}\circ\Psi^{(n)}\right)  \left(  \left(  \tfrac
{(x_{n}+c_{n})\nu_{n-1}+\eta_{n-1}}{\mu_{n-1}}\circ\Psi^{(n)}\right)
\tilde\kappa_{1}^{n-1}+\tilde\kappa_{2}^{n-1}\right) \\
& +\left(  (\mathrm{L}_{\kappa_{2}^{n}}\Phi_{n})\circ\Psi^{(n)}\right)
\,\tilde\kappa_{1}^{n}.
\end{align*}
By Proposition \ref{prop-automorphisms-02}, we know that there exist two
smooth functions $\mu_{n}$ and $\eta_{n}$ (with $\mu_{n}\neq0$) such that
\[
\Phi_{*}^{(n)}(\kappa_{2}^{n})=\left(  \eta_{n}\circ\Psi^{(n)}\right)
\tilde\kappa_{1}^{n}+\left(  \mu_{n}\circ\Psi^{(n)}\right)  \tilde\kappa
_{2}^{n}.
\]
Comparing the last two relations and taking into account that $\tilde
\kappa_{1}^{n}=\tfrac\partial{\partial\tilde x_{n}}$ while $\tilde\kappa
_{1}^{n-1}$, $\tilde\kappa_{2}^{n-1}$, and $\tilde\kappa_{2}^{n}$ have zeros
as components multiplying $\tfrac\partial{\partial\tilde x_{n}}$ we see that
$\eta_{n}=\mathrm{L}_{\kappa_{2}^{n}}\Phi_{n}$. From the inductive definition
of Kumpera-Ruiz normal forms (regular prolongation) given in
Section~\ref{sec-kumpera-ruiz}, we have
\[
\tilde\kappa_{2}^{n}=(\tilde x_{n}+\tilde c_{n})\tilde\kappa_{1}^{n-1}%
+\tilde\kappa_{2}^{n-1}.
\]
We can now conclude that $\mu_{n}=\mu_{n-1}$ and that
\[
\Phi_{n}(x)=\dfrac1{\mu_{n-1}}\left(  (x_{n}+c_{n})\nu_{n-1}+\eta
_{n-1}\right)  -\tilde c_{n},
\]
where
\[
\tilde c_{n}=c_{n}\dfrac{\nu_{n-1}(0)}{\mu_{n-1}(0)}+\dfrac{\eta_{n-1}(0)}%
{\mu_{n-1}(0)}.
\]
Now consider $\kappa_{1}^{n}$. Relation (\ref{triangular-tangent-map}) gives
$\Phi_{*}^{(n)}(\kappa_{1}^{n})=\left(  (\mathrm{L}_{\kappa_{1}^{n}}\Phi
_{n})\circ\Psi^{(n)}\right)  \,\tilde\kappa_{1}^{n}$, which implies $\nu
_{n}=\mathrm{L}_{\kappa_{1}^{n}}\Phi_{n}$. This obviously gives $\nu_{n}%
=\nu_{i-1}/\mu_{i-1}$.

\medskip\ 

\noindent\emph{Singular case\/}: If $w_{n-4}=a_{1}$ then we have $\kappa
_{2}^{n}=\kappa_{1}^{n-1}+x_{n}\kappa_{2}^{n-1}$. Together with relation
(\ref{triangular-tangent-map}) and with the induction hypothesis, this
relation leads to:
\begin{align*}
\Phi_{*}^{(n)}(\kappa_{2}^{n})  & =\Phi_{*}^{(n-1)}(\kappa_{1}^{n-1}%
)+(x_{n}\circ\Psi^{(n)})\Phi_{*}^{(n-1)}(\kappa_{2}^{n-1})\\
& +\left(  (\mathrm{L}_{\kappa_{2}^{n}}\Phi_{n})\circ\Psi^{(n)}\right)
\,\tilde\kappa_{1}^{n}\\
& =\left(  \left(  \nu_{n-1}+x_{n}\eta_{n-1}\right)  \circ\Psi^{(n)}\right)
\tilde\kappa_{1}^{n-1}+\left(  x_{n}\mu_{n-1}\circ\Psi^{(n)}\right)
\tilde\kappa_{2}^{n-1}\\
& +\left(  (\mathrm{L}_{\kappa_{2}^{n}}\Phi_{n})\circ\Psi^{(n)}\right)
\,\tilde\kappa_{1}^{n}\\
& =\left(  \left(  \nu_{n-1}+x_{n}\eta_{n-1}\right)  \circ\Psi^{(n)}\right)
\left(  \tilde\kappa_{1}^{n-1}+\left(  \tfrac{x_{n}\mu_{n-1}}{\nu_{n-1}%
+x_{n}\eta_{n-1}}\circ\Psi^{(n)}\right)  \tilde\kappa_{2}^{n-1}\right) \\
& +\left(  (\mathrm{L}_{\kappa_{2}^{n}}\Phi_{n})\circ\Psi^{(n)}\right)
\,\tilde\kappa_{1}^{n}.
\end{align*}
By Proposition \ref{prop-automorphisms-02}, we know that there exist two
functions $\mu_{n}$ and $\eta_{n}$ such that
\[
\Phi_{*}^{(n)}(\kappa_{2}^{n})=\left(  \eta_{n}\circ\Psi^{(n)}\right)
\tilde\kappa_{1}^{n}+\left(  \mu_{n}\circ\Psi^{(n)}\right)  \tilde\kappa
_{2}^{n}.
\]
The same argument as in the regular case implies $\eta_{n}=\mathrm{L}%
_{\kappa_{2}^{n}}\Phi_{n}$, $\mu_{n}=\nu_{n-1}+x_{n}\eta_{n-1}$, and
\[
\Phi_{n}(x)=\dfrac{x_{n}\mu_{n-1}}{\nu_{n-1}+x_{n}\eta_{n-1}}.
\]
Moreover, like in the regular case, the relation
\[
\Phi_{*}^{(n)}(\kappa_{1}^{n})=\left(  (\mathrm{L}_{\kappa_{1}^{n}}\Phi
_{n})\circ\Psi^{(n)}\right)  \,\tilde\kappa_{1}^{n}
\]
implies $\nu_{n}=\mathrm{L}_{\kappa_{1}^{n}}\Phi_{n}=(\mu_{i-1}\nu_{i-1}%
)/(\nu_{i-1}+x_{i}\eta_{i-1})^{2}$.\hfill$\square$

\subsection{Are Goursat Structures Locally Determined by Their Abnormals?}

In this Subsection we will be interested, in the case of Goursat structures,
in the following question asked by Jakubczyk: ``Are nonholonomic distributions
determined by their abnormal curves?''. Several results have been obtained
giving a positive answer to this question: for stable degenerations of Engel
structures by Zhitomirsk\u{\i}i~\cite{zhitomirskii-engel}, for singular
contact structures by Jakubczyk and Zhitomirsk\u{\i}%
i~\cite{jakubczyk-zhitomirskii}, for generic distributions of corank at least
equal to three, at typical points, by Montgomery~\cite{montgomery-survey}.
Recently, Jakubczyk~\cite{jakubczyk-complex-abnormals} has proved that the
answer is positive if we consider abnormal curves of the complexified problem,
for all distributions with the exception of a small subclass.\ We will show in
this Subsection that this subclass contains Goursat structures.

To start with, let us be more precise on what we mean by the statement that
distributions are determined by their abnormal curves. We will follow the
definitions given in~\cite{montgomery-survey}. Distributions that belong to a
class $\mathcal{Q}$ of distributions are \emph{strongly determined} by their
abnormal curves if, for any pair of distributions $\mathcal{D}$ and
$\tilde{\mathcal{D}}$ that belong to $\mathcal{Q}$, any local diffeomorphism
that transforms each abnormal curve of $\mathcal{D}$ into an abnormal curve of
$\tilde{\mathcal{D}}$, and the other way around, transforms also $\mathcal{D}$
into $\tilde{\mathcal{D}}$. It is clear that Goursat structures are not
strongly determined by their abnormal curves because they have very few
abnormal curves. For example, contact structures do not have any non-trivial
abnormal curve.

A weaker property can be defined as follows. Distributions that belong to a
class$~\mathcal{Q}$ of distributions are \emph{weakly determined} by their
abnormal curves if, for any pair of distributions $\mathcal{D}$ and
$\tilde{\mathcal{D}}$ that belong to$~\mathcal{Q}$, the existence of a local
diffeomorphism that transforms each abnormal curve of $\mathcal{D}$ into an
abnormal curve of $\tilde{\mathcal{D}}$, and the other way around, implies the
local equivalence of $\mathcal{D}$ and $\tilde{\mathcal{D}}$.

\begin{proposition}
\label{prop-weakly-determined}Goursat structures on $n$-manifolds are not
weakly determined by their abnormal curves if $n\geq6$.
\end{proposition}

\noindent\textbf{Proof of Proposition \ref{prop-weakly-determined} }Consider
the two following Kumpera-Ruiz normal forms defined on$~\mathbb{R}^{6}$ by
\[
\left(
\begin{array}
[c]{l}%
\tfrac\partial{\partial x_{6}}%
\end{array}
,
\begin{array}
[c]{l}%
x_{6}\tfrac\partial{\partial x_{5}}+x_{5}\tfrac\partial{\partial x_{4}}%
+x_{4}\tfrac\partial{\partial x_{3}}+x_{3}\tfrac\partial{\partial x_{2}%
}+\tfrac\partial{\partial x_{1}}%
\end{array}
\right)
\]
and
\[
\left(
\begin{array}
[c]{l}%
\tfrac\partial{\partial x_{6}}%
\end{array}
,
\begin{array}
[c]{l}%
(x_{6}+1)\tfrac\partial{\partial x_{5}}+\tfrac\partial{\partial x_{4}}%
+x_{5}\left(  x_{4}\tfrac\partial{\partial x_{3}}+x_{3}\tfrac\partial{\partial
x_{2}}+\tfrac\partial{\partial x_{1}}\right)
\end{array}
\right)  .
\]

On the one hand, by Theorem~\ref{thm-abnormal-curves}, the distributions
spanned by these two Kumpera-Ruiz normal forms have the same abnormal curves,
locally at zero. Indeed, for each of them, the submanifolds $S_{j}^{(j)}$, for
$j=0$ and $1$, are empty in a small enough neighborhood of zero (see
Proposition~\ref{prop-singularity-type-claim}); and thus their abnormal curves
are given, in both cases, by $\mathcal{A}^{(0)}=\mathcal{C}_{0}=(\tfrac
{\partial}{\partial x_{6}})$, in a small enough neighborhood of zero. But on
the other hand, it has been shown by Kumpera and Ruiz \cite{kumpera-ruiz} that
these two distributions are not locally equivalent at zero. Indeed, the first
one has singularity type $a_{0}a_{0}a_{0}$ at zero while the second one has
singularity type $a_{0}a_{1}a_{0}$ at zero. Analogous examples can be
constructed for any $n\geq6$.\hfill$\square$

\medskip

\ Our study of relations between abnormal curves and their singularity type
shows that the geometry of a Goursat structure is reflected by abnormal curves
of all elements of the derived flag. It is thus natural to introduce the
following definition. Distributions that belong to a class $\mathcal{Q}$ of
distributions are \emph{weakly determined by abnormal curves of their derived
flags} if, for any pair of distributions $\mathcal{D}$ and $\tilde
{\mathcal{D}}$ that belong to $\mathcal{Q}$, the existence of a local
diffeomorphism that transforms each abnormal curve of $\mathcal{D}^{(i)}$ into
an abnormal curve of $\tilde{\mathcal{D}}^{(i)}$, and the other way around,
for each $i\geq0$, implies the local equivalence of $\mathcal{D}$
and$~\tilde{\mathcal{D}}$. It is a direct consequence of
Theorem~\ref{thm-abnormal-curves} and of the classification obtained
in~\cite{cheaito-mormul}, \cite{gaspar}, and~\cite{kumpera-ruiz}, that Goursat
structures on $\mathbb{R}^{n}$, for $3\leq n\leq8$ are determined by abnormal
curves of their derived flags. It is surprising that in higher dimensions it
is not the case. Indeed, we have the following result which is a direct
consequence of Theorem~\ref{thm-abnormal-curves} and the Theorem announced
in~\cite{cheaito-mormul-pasillas-respondek}.

\begin{proposition}
\label{prop-very-weakly-determined}Goursat structures on $n$-manifolds are not
determined by abnormal curves of their derived flags if $n\geq9$.
\end{proposition}

It has already been announced in \cite{cheaito-mormul-pasillas-respondek} that
the growth vector is not a complete invariant for Goursat structures on
$\mathbb{R}^{n}$, for $n\geq9$ (which, together with
Theorem~\ref{thm-abnormal-curves}, implies the above result). We will give in
this Subsection our proof of this latter fact. An alternative proof can be
found in \cite{mormul-R9}. It is important to stress that the method used in
\cite{mormul-R9} and the method that we will present in this Section are
different. It seems that both methods apply, in general, to different cases of non-equivalence.

\medskip\ 

Our aim now is to prove Proposition~\ref{prop-very-weakly-determined}. This
will be done by giving an example (Proposition~\ref{prop-R9}) of two Goursat
structures $\mathcal{D}$ and $\tilde{\mathcal{D}}$ that are locally
non-equivalent but that have the same singularity type (which, by
Theorem~\ref{thm-abnormal-curves}, implies the existence of a diffeomorphism
between the abnormal curves of $\mathcal{D}^{(i)}$ and those of $\tilde
{\mathcal{D}}^{(i)}$, for $i\geq0$). Then, this example will be improved
(Proposition~\ref{prop-R11}) by constructing, instead of a pair of
distributions, a continuous family (parametrized by a real number) of locally
non-equivalent Goursat structures that have the same singularity type, and
thus diffeomorphic collections of abnormal curves for all elements of their
derived flags.

Consider two Kumpera-Ruiz normal forms $(\kappa_{1}^{n},\kappa_{2}^{n})$ and
$(\tilde\kappa_{1}^{n},\tilde\kappa_{2}^{n})$, defined on$~\mathbb{R}^{n}$,
centered at zero, and given, respectively, in coordinates $x=(x_{1}%
,\ldots,x_{n})$ and $\tilde x=(\tilde x_{1},\ldots,\tilde x_{n})$. Assume that
they have been obtained from $(\kappa_{1}^{i},\kappa_{2}^{i})$ and
$(\tilde\kappa_{1}^{i},\tilde\kappa_{2}^{i})$, respectively, by a sequence of
regular prolongations, for $i\geq3$. Suppose, moreover, that the Goursat
structures spanned by $(\kappa_{1}^{n},\kappa_{2}^{n})$ and $(\tilde\kappa
_{1}^{n},\tilde\kappa_{2}^{n})$ are locally equivalent and let $\tilde
x=\Phi(x)$ be a (generalized) contact transformation, of order $n-2$, that
establishes this equivalence. We have $\tilde x_{j}=\Phi_{j}(x)$, for $1\leq
j\leq n$. We are going to prove that the components $\Phi_{j}$, for $i+1\leq
j\leq n$, can be obtained by a sequence of derivations (with respect to a well
chosen vector field) from the component $\Phi_{i}$. To start with, apply
Theorem \ref{thm-automorphisms} to the component $\Phi_{i+1}$. We have
\[
\Phi_{i+1}=\dfrac1{\mu_{i}}\left(  (x_{i+1}+c_{i+1})\nu_{i}+\eta_{i}\right)
-\tilde c_{i+1}.
\]
In follows also from Theorem \ref{thm-automorphisms} (regular case) that this
expression can be written in the following form:
\begin{align*}
\Phi_{i+1}  & =\dfrac1{\mu_{i}}\left(  (x_{i+1}+c_{i+1})\mathrm{L}_{\kappa
_{1}^{i}}\Phi_{i}+\mathrm{L}_{\kappa_{2}^{i}}\Phi_{i}\right)  -\tilde
c_{i+1}\\
& =\dfrac1{\mu_{i}}\mathrm{L}_{(x_{i+1}+c_{i+1})\kappa_{1}^{i}+\kappa_{2}^{i}%
}\Phi_{i}-\tilde c_{i+1}\\
& =\mathrm{L}_{\tfrac1{\mu_{i}}\kappa_{2}^{i+1}}\Phi_{i}-\tilde c_{i+1}.
\end{align*}
But since $\Phi_{i}$ is a function of $x_{1},\ldots,x_{i}$ only, the latter
expression can be rewritten as
\[
\Phi_{i+1}=\mathrm{L}_{\tfrac1{\mu_{i}}\kappa_{2}^{n}}\Phi_{i}-\tilde
c_{i+1}.
\]
Theorem \ref{thm-automorphisms} implies, moreover, that $\mu_{j}=\mu_{i}$, for
$i+1\leq j\leq n$. Thus, the previous argument can be repeated to obtain, for
$1\leq k\leq n-i$, the following relations:
\[
\Phi_{i+k}=\mathrm{L}_{\tfrac1{\mu_{i}}\kappa_{2}^{n}}^{k}\Phi_{i}-\tilde
c_{i+k},
\]
which imply that
\[
\tilde c_{i+k}=\left(  \mathrm{L}_{\tfrac1{\mu_{i}}\kappa_{2}^{n}}^{k}\Phi
_{i}\right)  (0)
\]
because the coordinates are centered. Therefore, in the case of a sequence of
regular prolongations, the constants $\tilde c_{i+k}$ can be obtained by
computing the successive derivatives $\mathrm{L}_{(1/\mu_{i})\kappa_{2}^{n}%
}^{k}\Phi_{i}$ of the component $\Phi_{i}$ (that defines the coordinate
$\tilde x_{i}$) and by taking their values at zero.

The following definition is natural and will simplify the proofs of the next
results given in this Subsection. Let $\gamma$ be a smooth function defined on
$\mathbb{R}^{n}$ and let $g$ be a smooth vector field, also defined on
$\mathbb{R}^{n}$. The \emph{degree} of the function $\gamma$, with respect to
the vector field $g$, is the smallest integer $k$ (maybe infinite) such that
$\mathrm{L}_{g}^{k}(\gamma)(0)\neq0$. Note that if the degree of $\gamma_{1}$
is $i_{1}$ and the degree of $\gamma_{2}$ is $i_{2}$ then the degree of
$\gamma_{1}\gamma_{2}$ is obviously $i_{1}+i_{2}$.

\begin{proposition}
\label{prop-R9}Consider the two following Kumpera-Ruiz normal forms defined
on$~\mathbb{R}^{9}$ by
\begin{eqnarray*}
\kappa_{1}^{9}  & = & \tfrac{\partial}{\partial x_{9}}\\
\kappa_{2}^{9}(c_{9})  & = & (x_{9}+c_{9})\tfrac{\partial}{\partial x_{8}}%
+(x_{8}+1)\tfrac{\partial}{\partial x_{7}}+x_{7}\tfrac{\partial}{\partial
x_{6}}+\tfrac{\partial}{\partial x_{5}}\\
& & \mbox{} + x_{6}\left(  x_{5}\tfrac{\partial}{\partial x_{4}}+x_{4}\tfrac{\partial
}{\partial x_{3}}+x_{3}\tfrac{\partial}{\partial x_{2}}+\tfrac{\partial
}{\partial x_{1}}\right)  ,
\end{eqnarray*}
where $c_{9}=0$ or $1$. They are locally non-equivalent at zero, although both
of them have the same singularity type $a_{0}a_{0}a_{1}a_{2}a_{0}a_{0}$ at zero.
\end{proposition}

\noindent\textbf{Proof of Proposition \ref{prop-R9}} Denote by $\kappa^{9}$
the Kumpera-Ruiz normal form given by $(\kappa_{1}^{9},\kappa_{2}^{9}(0))$, in
$(x_{1},\ldots,x_{9})$-coordinates, and denote by $\tilde\kappa^{9}$ the
Kumpera-Ruiz normal form given by $(\kappa_{1}^{9},\kappa_{2}^{9}(\tilde
c_{9}))$, in $(\tilde x_{1},\ldots,\tilde x_{9})$-coordinates. We are going to
show that if a (generalized) contact transformation $\tilde x=\Phi(x)$
converts the Goursat structure generated by $\kappa^{9}$ into the one
generated by $\tilde\kappa^{9}$ then we must have $\tilde c_{9}=0$.

Denote by $\kappa^{4},\ldots,\kappa^{9}$ and by $\tilde{\kappa}^{4}%
,\ldots,\tilde{\kappa}^{9}$ the elements of the two sequences of Kumpera-Ruiz
normal forms used to construct, via prolongations, the normal forms $\kappa^{9}$
and~$\tilde{\kappa}^{9}$, respectively. Since $\kappa^{5}=R_{0}(\kappa^{4})$
we have, by the regular case of Theorem~\ref{thm-automorphisms}, the following
relations:
\[%
\begin{array}
[c]{ccc}%
\mu_{5}=\mu_{4} & \text{ and } & \nu_{5}=\dfrac{\nu_{4}}{\mu_{4}}%
\end{array}
.
\]
Hence $\mu_{5}$ and $\nu_{5}$ are functions of $x_{1},\ldots,x_{4}$ only.
Denote $\mu=\mu_{4}$, $\nu=\nu_{4}$, and $\eta=\eta_{5}$. Since $\kappa
^{6}=S(\kappa^{5})$ we have, by the singular case of Theorem
\ref{thm-automorphisms}, the following relations:
\begin{align*}
\Phi_{6}(x)  & =\frac{x_{6}\mu}{\dfrac{\nu}{\mu}+x_{6}\eta}\\
\mu_{6}  & =\dfrac{\nu}{\mu}+x_{6}\eta.
\end{align*}
Denote $\alpha=1/\mu_{6}$ and $g=\alpha\kappa_{2}^{9}$. Since both $\kappa
^{9}$ and $\tilde{\kappa}^{9}$ are obtained by a sequence of regular
prolongations from $\kappa^{6}$ and $\tilde{\kappa}^{6}$, respectively, it
follows from the discussion given at the beginning of this Subsection that the
new constant $\tilde{c}_{9}$ can be calculated by computing the successive
derivatives of $\Phi_{6}$, in the direction of the vector field $g=(1/\mu
_{6})\kappa_{2}^{9}$. Namely
\[
\tilde{c}_{9}=\left(  \mathrm{L}_{g}^{3}\Phi_{6}\right)  (0).
\]

Instead of computing the successive derivatives of $\Phi_{6}$ directly, take
the Taylor series expansion of$~\Phi_{6}$. The terms of this expansion that
contain coordinate functions of degree $d\geq4$, with respect to$~g$, can
obviously be discarded. To this aim, we will start by computing the degree,
with respect to$~g$, of the functions $x_{1},\ldots,x_{6}$, that is of the
variables on which$~\Phi_{6}$ depends.

For $x_{6}$, we have:
\begin{align*}
\mathrm{L}_{g}x_{6}  & =\alpha x_{7}\\
\mathrm{L}_{g}^{2}x_{6}  & =\alpha^{2}(x_{8}+1)+\left(  \mathrm{L}_{g}%
\alpha\right)  x_{7}\\
\mathrm{L}_{g}^{3}x_{6}  & =\alpha^{3}x_{9}+3\alpha\left(  \mathrm{L}%
_{g}\alpha\right)  (x_{8}+1)+\left(  \mathrm{L}_{g}^{2}\alpha\right)
x_{7}\text{.}%
\end{align*}
Since $\left(  \mathrm{L}_{g}^{2}x_{6}\right)  (0)=\alpha^{2}(0)=\left(
\mu(0)/\nu(0)\right)  ^{2}\neq0$, the degree of $x_{6}$ is $2$. We have
$\mathrm{L}_{g}x_{5}=\alpha$. Therefore the degree of $x_{5}$ is $1$. We have
$\mathrm{L}_{g}x_{4}=\alpha x_{6}x_{5}$. Thus the degree of $x_{4}$ is $4$.
Analogously, the degree of $x_{3}$ is $7$, the degree of $x_{2}$ is $10$, and
the degree of $x_{1}$ is $3$.

Now observe that $\Phi_{6}(x)=x_{6}\varphi(x_{1},\ldots,x_{6})$, for a
suitable function$~\varphi$. This implies that each term of the Taylor series
expansion of $\Phi_{6}$ is of the form $x_{6}x_{1}^{k_{1}}\cdots x_{6}^{k_{6}%
}$, for some integers $k_{1},\ldots,k_{6}$. Since $\tilde{c}_{9}=\left(
\mathrm{L}_{g}^{3}\Phi_{6}\right)  (0)$, we consider only terms of degree
$d\leq3$ with respect to$~g$. Therefore we have:
\[
\Phi_{6}(x)=Ax_{6}+Bx_{6}x_{5},
\]
up to terms of degree $d\geq4$ with respect to$~g$. Recall that neither $\mu$
nor $\nu$ depend on the variables$~x_{5}$ and$~x_{6}$. Hence
\begin{align*}
\dfrac{\partial\Phi_{6}}{\partial x_{6}}  & =\frac{\nu}{\left(  \dfrac{\nu
}{\mu}+x_{6}\eta\right)  ^{2}}\\
\frac{\partial^{2}\Phi_{6}}{\partial x_{5}\partial x_{6}}  & =\frac
{-2x_{6}\eta_{x_{5}}}{\left(  \dfrac{\nu}{\mu}+x_{6}\eta\right)  ^{3}}.
\end{align*}
Thus $A=\mu(0)$ and $B=0$. This implies that $\Phi_{6}(x)=\mu(0)x_{6}$, up to
terms of degree $d\geq4$. Since we have already computed the successive
derivatives of $x_{6}$, it is easy to obtain that:
\begin{align*}
\left(  \mathrm{L}_{g}\Phi_{6}\right)  (0)  & =0\\
\left(  \mathrm{L}_{g}^{2}\Phi_{6}\right)  (0)  & =\mu(0)\alpha^{2}(0)\\
\left(  \mathrm{L}_{g}^{3}\Phi_{6}\right)  (0)  & =3\mu(0)\alpha(0)\left(
\mathrm{L}_{g}\alpha\right)  (0)\text{.}%
\end{align*}
But $\left(  \mathrm{L}_{g}\alpha\right)  (0)=0$. Hence, since $\tilde{c}%
_{9}=\left(  \mathrm{L}_{g}^{3}\Phi_{6}\right)  (0)$, we have $\tilde{c}_{9}=0
$.\hfill$\square$

\begin{proposition}
\label{prop-R11}Consider the following family of Kumpera-Ruiz normal forms
defined on $\mathbb{R}^{11}$ by
\begin{eqnarray*}
\kappa_{1}^{11}  & = & \tfrac{\partial}{\partial x_{11}}\\
\kappa_{2}^{11}(c_{11})  & = & (x_{11}+c_{11})\tfrac{\partial}{\partial x_{10}%
}+(x_{10}+1)\tfrac{\partial}{\partial x_{9}}+(x_{9}+1)\tfrac{\partial
}{\partial x_{8}}+x_{8}\tfrac{\partial}{\partial x_{7}}\\
& & \mbox{} + x_{7}\tfrac{\partial}{\partial x_{6}}+\tfrac{\partial}{\partial x_{5}%
}+x_{6}\left(  x_{5}\tfrac{\partial}{\partial x_{4}}+x_{4}\tfrac{\partial
}{\partial x_{3}}+x_{3}\tfrac{\partial}{\partial x_{2}}+\tfrac{\partial
}{\partial x_{1}}\right)  .
\end{eqnarray*}
where $c_{11}$ is an arbitrary real constant. Two Kumpera-Ruiz normal forms
that belong to this family are locally equivalent at zero if and only if they
have the same constant~$c_{11}$, although all of them have the same
singularity type $a_{0}a_{0}a_{1}a_{2}a_{3}a_{0}a_{0}a_{0}$ at zero.
\end{proposition}

\noindent\textbf{Proof of Proposition \ref{prop-R11}} Denote by $\kappa^{11}$
the Kumpera-Ruiz normal form given by $(\kappa_{1}^{11},\kappa_{2}^{11}%
(c_{11}))$, in $(x_{1},\ldots,x_{11})$-coordinates, and denote by
$\tilde\kappa^{11}$ the Kumpera-Ruiz normal form given by $(\kappa_{1}%
^{11},\kappa_{2}^{11}(\tilde c_{11}))$, in $(\tilde x_{1},\ldots,\tilde
x_{11})$-coordinates. We are going to show that if a (generalized) contact
transformation $\tilde x=\Phi(x)$ converts the Goursat structure generated by
$\kappa^{11}$ into the one generated by $\tilde\kappa^{11}$ then we must have
$\tilde c_{11}=c_{11}$.

Denote by $\kappa^{4},\ldots,\kappa^{11}$ and by $\tilde{\kappa}^{4}%
,\ldots,\tilde{\kappa}^{11}$ the elements of the two sequences of Kumpera-Ruiz
normal forms used to construct, via prolongations, the normal forms~$\kappa^{11}$
and~$\tilde{\kappa}^{11}$, respectively. Since $\kappa^{5}=R_{0}(\kappa^{4})$
we have, by the regular case of Theorem~\ref{thm-automorphisms}, the following
relations:
\[%
\begin{array}
[c]{ccc}%
\mu_{5}=\mu_{4} & \text{ and } & \nu_{5}=\dfrac{\nu_{4}}{\mu_{4}}%
\end{array}
.
\]
Hence $\mu_{5}$ and $\nu_{5}$ are functions of $x_{1},\ldots,x_{4}$ only.
Denote $\mu=\mu_{4}$, $\nu=\nu_{4}$, and $\eta=\eta_{5}$. Since $\kappa
^{6}=S(\kappa^{5})$ we have, by the singular case of Theorem
\ref{thm-automorphisms}, the following relations:
\begin{align*}
\Phi_{6}(x)  & =\frac{x_{6}\mu}{\dfrac{\nu}{\mu}+x_{6}\eta}\\
\mu_{6}  & =\dfrac{\nu}{\mu}+x_{6}\eta.
\end{align*}
Denote $\alpha=1/\mu_{6}$ and $g=\alpha\kappa_{2}^{11}$. Since both
$\kappa^{11}$ and $\tilde{\kappa}^{11}$ are obtained by a sequence of regular
prolongations from $\kappa^{6}$ and $\tilde{\kappa}^{6}$, respectively, it
follows from the discussion given at the beginning of this Subsection that the
new constant $\tilde{c}_{11}$ can be obtained by computing the successive
derivatives of $\Phi_{6}$, in the direction of the vector field $g$. Namely
\[
\tilde{c}_{11}=\left(  \mathrm{L}_{g}^{5}\Phi_{6}\right)  (0).
\]

Let us consider the Taylor series expansion of $\Phi_{6}$. Again, the terms of
this expansion that contain coordinate functions of degree $d\geq6$, with
respect to$~g$, will be discarded. The successive derivatives of $x_{6}$ are
given by:
\begin{align*}
\mathrm{L}_{g}x_{6}  & =\alpha x_{7}\\
\mathrm{L}_{g}^{2}x_{6}  & =\alpha^{2}x_{8}+(\mathrm{L}_{g}\alpha)x_{7}\\
\mathrm{L}_{g}^{3}x_{6}  & =\alpha^{3}(x_{9}+1)+3\alpha(\mathrm{L}_{g}%
\alpha)x_{8}+(\mathrm{L}_{g}^{2}\alpha)x_{7}\\
\mathrm{L}_{g}^{4}x_{6}  & =\alpha^{4}(x_{10}+1)+6\alpha^{2}(\mathrm{L}%
_{g}\alpha)(x_{9}+1)\\
& +\left(  3(\mathrm{L}_{g}\alpha)^{2}+4\alpha(\mathrm{L}_{g}^{2}%
\alpha)\right)  x_{8}+(\mathrm{L}_{g}^{3}\alpha)x_{7}\\
\mathrm{L}_{g}^{5}x_{6}  & =\alpha^{5}(x_{11}+c_{11})+10\alpha^{3}%
(\mathrm{L}_{g}\alpha)(x_{10}+1)\\
& +\left(  15\alpha(\mathrm{L}_{g}\alpha)^{2}+10\alpha^{2}(\mathrm{L}_{g}%
^{2}\alpha)\right)  (x_{9}+1)\\
& +\left(  10(\mathrm{L}_{g}\alpha)(\mathrm{L}_{g}^{2}\alpha)+5\alpha
(\mathrm{L}_{g}^{3}\alpha)\right)  x_{8}+(\mathrm{L}_{g}^{4}\alpha)x_{7}.
\end{align*}
Since $\left(  \mathrm{L}_{g}^{3}x_{6}\right)  (0)=\alpha^{3}(0)=\left(
\mu(0)/\nu(0)\right)  ^{3}\neq0$, the degree of $x_{6}$ is $3$. The degree of~$x_{5}$
is~$1$, the degree of $x_{4}$ is $5$, the degree of $x_{3}$ is $9$,
the degree of $x_{2}$ is $13$, and the degree of $x_{1}$ is $4$ (all degrees
are with respect to$~g$).

Now observe that $\Phi_{6}(x)=x_{6}\varphi(x_{1},\ldots,x_{6})$, for a
suitable function$~\varphi$. The Taylor series expansion of $\Phi_{6}$, up to
terms of degree $d\leq5$ with respect to$~g$, is given by:
\[
\Phi_{6}(x)=Ax_{6}+Bx_{6}x_{5}+Cx_{6}x_{5}^{2}.
\]
Recall that neither $\mu$ nor $\nu$ depend on the variable~$x_{5}$. Therefore
\begin{align*}
\dfrac{\partial\Phi_{6}}{\partial x_{6}}(0)  & =\mu(0)\\
\frac{\partial^{2}\Phi_{6}}{\partial x_{5}\partial x_{6}}(0)  & =\frac
{\partial^{3}\Phi_{6}}{\partial^{2}x_{5}\partial x_{6}}(0)=0.
\end{align*}
Hence $A=\mu(0)$ and both $B$ and $C$ are equal to $0$. This implies that
$\Phi_{6}(x)=\mu(0)x_{6}$, up to terms of degree $d\geq6$. Since we have
already computed the successive derivatives of $x_{6}$, it is easy to obtain
that:
\begin{align*}
\left(  \mathrm{L}_{g}\Phi_{6}\right)  (0)  & =0\\
\left(  \mathrm{L}_{g}^{2}\Phi_{6}\right)  (0)  & =0\\
\left(  \mathrm{L}_{g}^{3}\Phi_{6}\right)  (0)  & =\mu(0)\alpha^{3}(0)\\
\left(  \mathrm{L}_{g}^{4}\Phi_{6}\right)  (0)  & =\mu(0)\alpha^{4}(0)\\
\left(  \mathrm{L}_{g}^{5}\Phi_{6}\right)  (0)  & =\mu(0)\alpha^{5}(0)c_{11}.
\end{align*}
Since $\Phi$ transforms $\kappa^{11}$ into $\tilde\kappa^{11}$, we must have
both $\tilde c_{9}=1$ and $\tilde c_{10}=1$. But $\tilde c_{9}=\left(
\mathrm{L}_{g}^{3}\Phi_{6}\right)  (0)$ and $\tilde c_{10}=\left(
\mathrm{L}_{g}^{4}\Phi_{6}\right)  (0)$. Therefore, $\mu^{4}(0)/\nu^{3}(0)=1$
and $\mu^{5}(0)/\nu^{4}(0)=1$. This obviously implies $\mu(0)=\nu(0)=1$.
Hence, since $\tilde c_{11}=c_{11}\mu^{6}(0)/\nu^{5}(0)$, we have $\tilde
c_{11}=c_{11}$.\hfill$\square$

\newpage

\appendix

\section{Weber's Problem}

\label{sec-weber}

Our proof of Kumpera-Ruiz's Theorem was based on the following fact: If a rank
two distribution$~\mathcal{D}$ on a manifold $M$ of dimension $n\geq4$
satisfies $\dim\mathcal{D}^{(1)}(p)=3$ and $\dim\mathcal{D}^{(2)}(p)=4$, for
each point $p$ in $M$, then there exists a canonical line field $\mathcal{L}%
\subset\mathcal{D}$ that satisfies $[\mathcal{L},\mathcal{D}^{(1)}%
]\subset\mathcal{D}^{(1)}$. This observation has a natural generalization: If
a rank $k\geq2$ distribution $\mathcal{D}$ on a manifold $M$ of dimension
$n\geq k+2$ satisfies $\dim\mathcal{D}^{(1)}(p)=k+1$ and $\dim\mathcal{D}%
^{(2)}(p)=k+2$, for each point $p$ in $M$, then there exists (i) a canonical
involutive distribution $\mathcal{L}_{1}\subset\mathcal{D}^{(0)}$ that has
rank $k-1$ and is uniquely characterized by $[\mathcal{L}_{1},\mathcal{D}%
^{(1)}]\subset\mathcal{D}^{(1)}$; and (ii) a canonical involutive distribution
$\mathcal{L}_{0}\subset\mathcal{D}^{(0)}$ that has rank $k-2$ and is uniquely
characterized by $[\mathcal{L}_{0},\mathcal{D}^{(0)}]\subset\mathcal{D}^{(0)}$
(see \cite{kumpera-ruiz} and \cite{martin-rouchon-driftless} for an approach
based on Pfaffian systems; see also \cite{kazarian-montgomery-shapiro} and
Proposition~\ref{prop-cartan}).

Though the above observation appears more or less clearly in the work of E.
Cartan (see e.g.~\cite{cartan-equivalence-absolue}; see also \cite{goursat}),
its origin can be found in the pioneering work of F.~Engel~\cite{engel}, for
$n=k+2$, and E.~von~Weber~\cite{weber-article}, for $n\geq k+2$ (see
also~\cite{cartan-weber}). This observation is closely related to the
following result, which is clearly stated in~Weber's article~\cite[Theorem
V]{weber-article} (using the dual language of Pfaffian systems).

\begin{theorem}
[E. von Weber]\label{thm-weber}Let $\mathcal{D}$ be a rank $k\geq2$
distribution on a manifold $M$ of dimension $n=m+k-2\geq4$. Assume that
$\dim\mathcal{D}^{(1)}(p)=k+1$ and $\dim\mathcal{D}^{(2)}(p)=k+2$, for each
point $p$ in $M$. Then, in a small enough neighborhood of any point $p$ in
$M$, the distribution~$\mathcal{D}$ is equivalent to a distribution spanned by
a family of vector fields that has the following form:
\begin{equation}
\left(
\begin{array}
[c]{c}%
\tfrac{\partial}{\partial x_{m+k-2}}%
\end{array}
,\ldots,
\begin{array}
[c]{c}%
\tfrac{\partial}{\partial x_{m+1}}%
\end{array}
,
\begin{array}
[c]{c}%
\tfrac{\partial}{\partial x_{m}}%
\end{array}
,
\begin{array}
[c]{c}%
x_{m}\tfrac{\partial}{\partial x_{m-1}}+%
{\textstyle\sum\limits_{i=2}^{m-2}}
\varphi_{i}(\overline{x}_{m-1})\tfrac{\partial}{\partial x_{i}}+\tfrac
{\partial}{\partial x_{1}}%
\end{array}
\right)  ,\label{weber-preliminary-normal-form}%
\end{equation}
where the functions $\varphi_{i}$, for $2\leq i\leq m-2$, depend on the
variables $x_{1},\ldots,x_{m-1}$ only.
\end{theorem}

The following result is a direct consequence of Theorem~\ref{thm-weber}.

\begin{proposition}
\label{prop-preliminary-rigid}Any Goursat structure on a manifold $M$ of
dimension $n\geq4$ is equivalent, in a small enough neighborhood of any point
$p$ in $M$, to a distribution spanned by a pair of vector fields that has the
following form:
\begin{equation}
\left(
\begin{array}
[c]{c}%
\tfrac{\partial}{\partial x_{n}}%
\end{array}
,
\begin{array}
[c]{c}%
x_{n}\tfrac{\partial}{\partial x_{n-1}}+x_{n-1}\tfrac{\partial}{\partial
x_{n-2}}+%
{\textstyle\sum\limits_{i=2}^{n-3}}
\varphi_{i}(\overline{x}_{n-1})\tfrac{\partial}{\partial x_{i}}+\tfrac
{\partial}{\partial x_{1}}%
\end{array}
\right)  ,\label{weber-form-a}%
\end{equation}
where the coordinates $x_{1},\ldots,x_{n}$ are centered at $p$ and the
functions $\varphi_{i}$, for $2\leq i\leq n-3$, depend on the variables
$x_{1},\ldots,x_{n-1}$ only.
\end{proposition}

In the particular case of four-manifolds the last result gives:

\begin{corollary}
[Engel's Theorem]\label{thm-engel}Any Goursat structure on a four-manifold $M$
is equivalent, in a small enough neighborhood of any point $p$ in $M$, to the
distribution spanned by the following pair of vector fields (Engel's normal
form):
\[
\left(
\begin{array}
[c]{c}%
\tfrac{\partial}{\partial x_{4}}%
\end{array}
,
\begin{array}
[c]{c}%
x_{4}\tfrac{\partial}{\partial x_{3}}+x_{3}\tfrac{\partial}{\partial x_{2}%
}+\tfrac{\partial}{\partial x_{1}}%
\end{array}
\right)  ,
\]
where the $x$-coordinates are centered at $p$.
\end{corollary}

The following Theorem can be considered as a rigorous version of Weber's
result~\cite[Theorem VI]{weber-article}.\ Although it is a direct consequence
of the work of Kumpera and Ruiz~\cite{kumpera-ruiz}, Martin and
Rouchon~\cite{martin-rouchon-driftless}, and Murray~\cite{murray-nilpotent},
it seems that it has never been stated in the following explicit form.

\begin{theorem}
[Weber's Problem]\label{thm-weber-problem} A rank $k\geq2$ distribution
$\mathcal{D}$ on a manifold $M$ of dimension $n=m+k-2\geq4$ is equivalent, in
a small enough neighborhood of a given point $p$ in $M$, to the distribution
spanned by the following family of vector fields (Weber's normal form)
\begin{equation}
\left(
\begin{array}
[c]{c}%
\tfrac{\partial}{\partial x_{m+k-2}}%
\end{array}
,\ldots,
\begin{array}
[c]{c}%
\tfrac{\partial}{\partial x_{m+1}}%
\end{array}
,
\begin{array}
[c]{c}%
\tfrac{\partial}{\partial x_{m}}%
\end{array}
,
\begin{array}
[c]{c}%
x_{m}\tfrac{\partial}{\partial x_{m-1}}+\cdots+x_{3}\tfrac{\partial}{\partial
x_{2}}+\tfrac{\partial}{\partial x_{1}}%
\end{array}
\right) \label{weber-normal-form}%
\end{equation}
if and only if $\dim\mathcal{D}_{i}(p)=\dim\mathcal{D}^{(i)}(p)=k+i$, for
$0\leq i\leq m-2$, in a small enough neighborhood of $p$.
\end{theorem}

If we have $\dim\mathcal{D}^{(i)}(p)=k+i$, for $0\leq i\leq n-2$, but we do
not impose any condition on $\dim\mathcal{D}_{i}(p)$ then we still have the
following result, which is a direct consequence of Theorem~\ref{thm-weber} and
Theorem~\ref{thm-kumpera-ruiz}, applied to the last two vectors fields
of~(\ref{weber-preliminary-normal-form}).

\begin{theorem}
[Kumpera-Ruiz]\label{thm-weber-kumpera-ruiz} Let $\mathcal{D}$ be a rank
$k\geq2$ distribution on a manifold $M$ of dimension $n=m+k-2\geq4$, such that
for any point $p$ in $M$ we have $\dim\mathcal{D}^{(i)}(p)=k+i$, for $0\leq
i\leq m-2$. Then, the distribution $\mathcal{D}$ is equivalent, in a small
enough neighborhood of any point $p$ in $M$, to the distribution spanned by
the following family of vector fields:
\[
\left(
\begin{array}
[c]{c}%
\tfrac{\partial}{\partial x_{m+k-2}}%
\end{array}
,\ldots,
\begin{array}
[c]{c}%
\tfrac{\partial}{\partial x_{m+1}}%
\end{array}
,
\begin{array}
[c]{c}%
\kappa_{1}^{m}%
\end{array}
,
\begin{array}
[c]{c}%
\kappa_{2}^{m}%
\end{array}
\right)  ,
\]
where the pair of vector fields $(\kappa_{1}^{m},\kappa_{2}^{m})$ denotes a
Kumpera-Ruiz normal form on$~\mathbb{R}^{m}$.
\end{theorem}

\section{Additional Normal Forms}

\label{sec-additional}

Let $\xi^{m}=(\xi_{1}^{m},\xi_{2}^{m})$ be a pair of vector fields defined on
$\mathbb{R}^{m}$ that has the following form:
\begin{equation}
\left(
\begin{array}
[c]{c}%
\tfrac{\partial}{\partial x_{m}}%
\end{array}
,
\begin{array}
[c]{c}%
x_{m}\tfrac{\partial}{\partial x_{m-1}}+x_{m-1}\tfrac{\partial}{\partial
x_{m-2}}+%
{\textstyle\sum\limits_{i=2}^{m-3}}
\varphi_{i}(\overline{x}_{m-1})\tfrac{\partial}{\partial x_{i}}+\tfrac
{\partial}{\partial x_{1}}%
\end{array}
\right)  .\label{weber-form-bis}%
\end{equation}
A pair of vector fields $\xi^{m+l}=(\xi_{1}^{m+l},\xi_{2}^{m+l})$ defined on
$\mathbb{R}^{m+l}$, for $l\geq0$, is called a \emph{prolongation of order} $l$
of $\xi^{m}$ if we have $\xi^{m+l}=\sigma_{l}\circ\cdots\circ\sigma_{1}%
(\xi^{m})$, where each $\sigma_{i}$, for $1\leq i\leq l$, equals either $S$ or
$R_{c_{i}}$, for some real constants~$c_{i}$ (recall that the singular and
regular prolongations $S$ and $R_{c_{i}}$ have been defined in
Section~\ref{sec-kumpera-ruiz}).

The following Lemma is a natural generalization of
Proposition~\ref{prop-preliminary-rigid}.

\begin{lemma}
\label{lem-cartan-goursat-form}Let $\mathcal{D}$ be a rank $k\geq2$
distribution on a manifold $M$ of dimension $n=m+l+k-2\geq4$, where $l$ and
$m$ are two non-negative integers. Assume that for each point $p$ in $M$ we
have $\dim\mathcal{D}^{(i)}(p)=k+i$, for $0\leq i\leq l+2$. Then, in a small
enough neighborhood of any point $p$ in$~M$, the distribution~$\mathcal{D}$ is
equivalent to a distribution spanned by a family of vector fields that has the
following form:
\begin{equation}
\left(
\begin{array}
[c]{c}%
\tfrac{\partial}{\partial x_{m+l+k-2}}%
\end{array}
,\ldots,
\begin{array}
[c]{c}%
\tfrac{\partial}{\partial x_{m+l+1}}%
\end{array}
,
\begin{array}
[c]{c}%
\xi_{1}^{m+l}%
\end{array}
,
\begin{array}
[c]{c}%
\xi_{2}^{m+l}%
\end{array}
\right)  ,\label{cartan-goursat-form}%
\end{equation}
where the pair of vector fields $(\xi_{1}^{m+l},\xi_{2}^{m+l})$ is a
prolongation of order $l$ of a pair of vector fields $(\xi_{1}^{m},\xi_{2}%
^{m})$ of the form \emph{(\ref{weber-form-bis})}.
\end{lemma}

The proof of Lemma~\ref{lem-cartan-goursat-form} is left to the reader. For
generic points, the Lemma is stated and proved in the work of
Cartan~\cite{cartan-equivalence-absolue} and Goursat~ \cite{goursat}. For singular
points, the Lemma is a direct consequence of the results obtained by Kumpera and
Ruiz~\cite{kumpera-ruiz} and its proof is almost the same as that of
Theorem~\ref{thm-kumpera-ruiz} but there are essentially two differences. The
first difference is that instead of using Proposition~\ref{prop-extended-engel},
as it is done in the Proof of Theorem~\ref{thm-kumpera-ruiz}, one uses
Theorem~\ref{thm-weber}; the second difference is that instead of starting the
induction argument, for $l=0$, with the Pfaff-Darboux normal form, as it is done 
in the Proof of Theorem~\ref{thm-kumpera-ruiz}, one starts it with Weber's
preliminary normal form (\ref{weber-form-bis}). 

Let $\mathcal{D}$ be a rank $k\geq2$ distribution on a manifold $M$ of
dimension $n=m+k-2\geq4$, such that for any point $p$ in $M$ we have
$\dim\mathcal{D}^{(i)}(p)=k+i$, for $0\leq i\leq m-2$. It is easy to check
that each distribution $\mathcal{D}^{(i)}$, for $0\leq i\leq m-4$, contains a
unique involutive subdistribution $\mathcal{C}_{i}\subset\mathcal{D}^{(i)}$
that has constant corank one in $\mathcal{D}^{(i)}$ and is characteristic for
$\mathcal{D}^{(i+1)}$. We can generalize the canonical submanifolds
$S_{0}^{(i)}$ of Section~\ref{sec-singularity-type} by the following
definition:
\[
S_{0}^{(i)}=\{p\in M:\mathcal{D}^{(i)}(p)=\mathcal{C}_{i+1}(p)\},
\]
where $0\leq i\leq m-5$. We say that a point $p$ of $M$ is \emph{singular} if
there exists some $0\leq i\leq m-5$ such that $p\in S_{0}^{(i)}$. For a
singular point$~p$, we denote by $k_{0}$ the smallest integer $1\leq k_{0}\leq
m-4$ such that $p\in S_{0}^{(k_{0}-1)}$.

\begin{lemma}
\label{lem-abnormal-normal-form}Let $\mathcal{D}$ be a rank $k\geq2$
distribution on a manifold $M$ of dimension $n=m+k_{0}+k-2\geq4$, such that
for any point $p$ in $M$ we have $\dim\mathcal{D}^{(i)}(p)=k+i$, for $0\leq
i\leq k_{0}+2$. Assume, moreover, that $k_{0}$ is the smallest integer such
that $\mathcal{D}^{(k_{0}-1)}(p)=\mathcal{C}_{k_{0}}(p)$. Then, in a small
enough neighborhood of $p$, the distribution~$\mathcal{D}$ is equivalent to a
distribution spanned by a family of vector fields that has the following
form:
\begin{equation}
\left(
\begin{array}
[c]{c}%
\tfrac{\partial}{\partial x_{m+k_{0}+k-2}}%
\end{array}
,\ldots,
\begin{array}
[c]{c}%
\tfrac{\partial}{\partial x_{m+k_{0}+1}}%
\end{array}
,
\begin{array}
[c]{c}%
\xi_{1}^{m+k_{0}}%
\end{array}
,
\begin{array}
[c]{c}%
\xi_{2}^{m+k_{0}}%
\end{array}
\right)  ,\label{cartan-goursat-abnormal-form}%
\end{equation}
where the pair of vector fields $\xi^{m+k_{0}}=(\xi_{1}^{m+k_{0}},\xi
_{2}^{m+k_{0}})$ is a prolongation of order $k_{0}$ of a pair of vector fields
$\xi^{m}=(\xi_{1}^{m},\xi_{2}^{m})$ of the form \emph{(\ref{weber-form-bis})}.
Moreover, we have $\xi^{m+k_{0}}=\sigma_{k_{0}}\circ\cdots\circ\sigma_{1}%
(\xi^{m})$, where $\sigma_{1}=S$ and each $\sigma_{j}$, for $2\leq j\leq
k_{0}$, equals $R_{c_{j}}$, for some real constants~$c_{i}$.
\end{lemma}

The proof of Lemma~\ref{lem-abnormal-normal-form} follows the same line as that 
of Proposition~\ref{prop-singularity-type-claim}. Though instead of
considering a Kumpera-Ruiz normal form we consider now a family of vector fields
of the form (\ref{cartan-goursat-abnormal-form}), the idea is the same.
Firstly, we compute the distributions $\mathcal{D}^{(i)}$ and $\mathcal{C}%
_{i}$, and the submanifolds $S_{0}^{(i)}$. Secondly, we observe that if
$\sigma_{1}=R_{c}$, for some real constant~$c$, then $p\notin S_{0}^{(k_{0}-1)}%
$; since $p\in S_{0}^{(k_{0}-1)}$ we must have $\sigma_{1}=S$. Thirdly, we
observe that if $\sigma_{j}=S$ for some $2\leq j\leq k_{0}$ then $p\in
S_{0}^{(k_{0}-j)}$; since $k_{0}$ is by definition the smallest integer such
that $p\in S_{0}^{(k_{0}-1)}$ we must have $\sigma_{j}=R_{c_{j}}$, for $2\leq
j\leq k_{0}$.

For Goursat structures, using the singularity type leads to the following
stronger result, which states that if the singularity type is of the form
$wa_{1}a_{2}\cdots a_{k_{0}}$ then the constants that appear in all regular
prolongations in the above Lemma equal zero.

\begin{lemma}
\label{lem-pre-rigid}Let $\mathcal{D}$ be a Goursat structure on a manifold
$M$ of dimension $n\geq5$ and let $p$ be a point in $M$. If the singularity
type of $\mathcal{D}$ at $p$ is of the form $wa_{1}a_{2}\cdots a_{k_{0}}$, for
some $1\leq k_{0}\leq n-4$, where $w$ is an arbitrary word of $J_{n-k_{0}-3}$,
then$~\mathcal{D}$ is locally equivalent to a distribution spanned by a pair
of vector fields that has the following form:
\begin{eqnarray*}
\xi_{1}  & = & \tfrac{\partial}{\partial x_{1}}\\
\xi_{2}  & = & x_{1}\tfrac{\partial}{\partial x_{2}}+\cdots+x_{k_{0}}%
\tfrac{\partial}{\partial x_{k_{0}+1}}+\tfrac{\partial}{\partial x_{k_{0}+2}%
} \\
 & & \mbox{} + x_{k_{0}}\left(  x_{k_{0}+2}\tfrac{\partial}{\partial x_{k_{0}+3}}%
+x_{k_{0}+3}\tfrac{\partial}{\partial x_{k_{0}+4}}+%
{\textstyle\sum\limits_{i=k_{0}+5}^{n}}
\varphi_{i}(x)\tfrac{\partial}{\partial x_{i}}\right)  ,
\end{eqnarray*}
where the coordinates $x_{1},\ldots,x_{n}$ are centered at $p$.
\end{lemma}

The proof of the last Lemma follows also the same line as the Proof of
Proposition~\ref{prop-singularity-type-claim}. Again, we leave details to the
reader. The main interest of the last Lemma is that it gives directly the
proof of Lemma~\ref{lem-rigid}.

\medskip\ 

\noindent\textbf{Proof of Lemma \ref{lem-rigid} }It is straightforward to
check that, in the coordinates of Lemma~\ref{lem-pre-rigid}, the canonical
submanifold $S_{k_{0}-1}^{(k_{0}-1)}$ is given by
$$S_{k_{0}-1}^{(k_{0}-1)}=\{x_{1}=0,\ldots,x_{k_{0}}=0\}$$
and that, moreover, we have
$\mathcal{C}_{0}=(\xi_{1})$ on $M$ and $\mathcal{A}_{k_{0}-1}^{(0)}%
(p)=(\xi_{2})(p)$ for each point~$p$ on $S_{k_{0}-1}^{(k_{0}-1)}$. In order to
obtain the required normal form, we only have to change two coordinates. For
$1\leq i\leq k_{0}+2$ and $k_{0}+5\leq i\leq n$, take~$y_{i}=x_{i}$. Moreover,
take~$y_{k_{0}+4}=x_{k_{0}+3}$ and $y_{k_{0}+3}=x_{k_{0}+4}-x_{k_{0}%
+3}x_{k_{0}+1}+\tfrac{1}{2}x_{k_{0}+2}x_{k_{0}+1}^{2}$.\hfill$\square$

\newpage

\section{Figures of Low-Dimensional Trailer Systems}

\label{sec-figures}

\subsection{The Unicycle and the Car}

\medskip

\begin{center}
\begin{figure}[h]
\begin{minipage}{0.5\hsize}
\begin{center}
\includegraphics[width=\hsize]{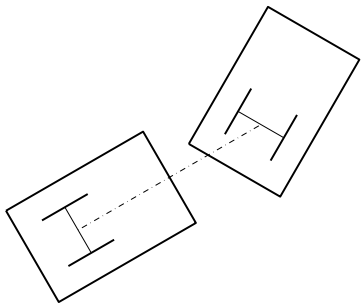}
\end{center}
\end{minipage}
\hspace{0.1\hsize} \begin{minipage}{0.3\hsize}
\begin{center}
$$
\begin{array}{ccl}
\dot x_1 & = & u_1 \\
\dot x_2 & = & x_1\,u_2 \\
\dot x_3 & = & u_2
\end{array}
$$
\end{center}
\end{minipage}
\caption{The unicycle and its normal form. Growth vector: $(2,3)$. Singularity
type:~$\epsilon$.}%
\end{figure}

\bigskip

\begin{figure}[h]
\begin{minipage}[h]{0.5\hsize}
\begin{center}
\includegraphics[width=\hsize]{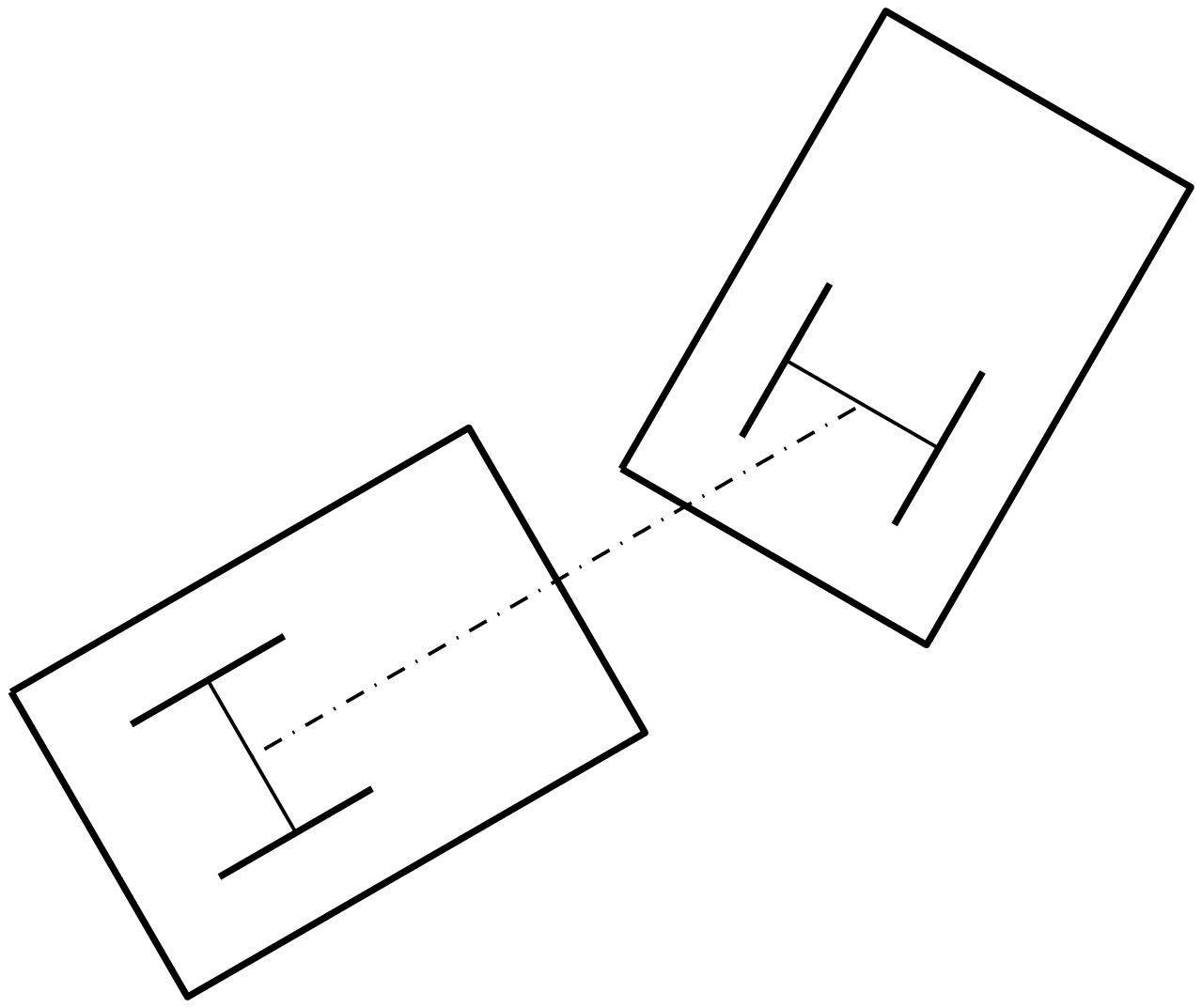}
\end{center}
\end{minipage}
\hspace{0.1\hsize} \begin{minipage}[h]{0.3\hsize}
\begin{center}
$$
\begin{array}{ccl}
\dot x_1 & = & u_1 \\
\dot x_2 & = & x_1\,u_2 \\
\dot x_3 & = & x_2\,u_2 \\
\dot x_4 & = & u_2
\end{array}
$$
\end{center}
\end{minipage}
\caption{The car and its normal form. Growth vector: $(2,3,4)$. Singularity
type:~$a_{0}$.}%
\end{figure}
\end{center}

\newpage

\subsection{The Two-Trailer System}

\begin{center}
\begin{figure}[h]
\begin{minipage}[h]{0.5\hsize}
\begin{center}
\includegraphics[width=\hsize]{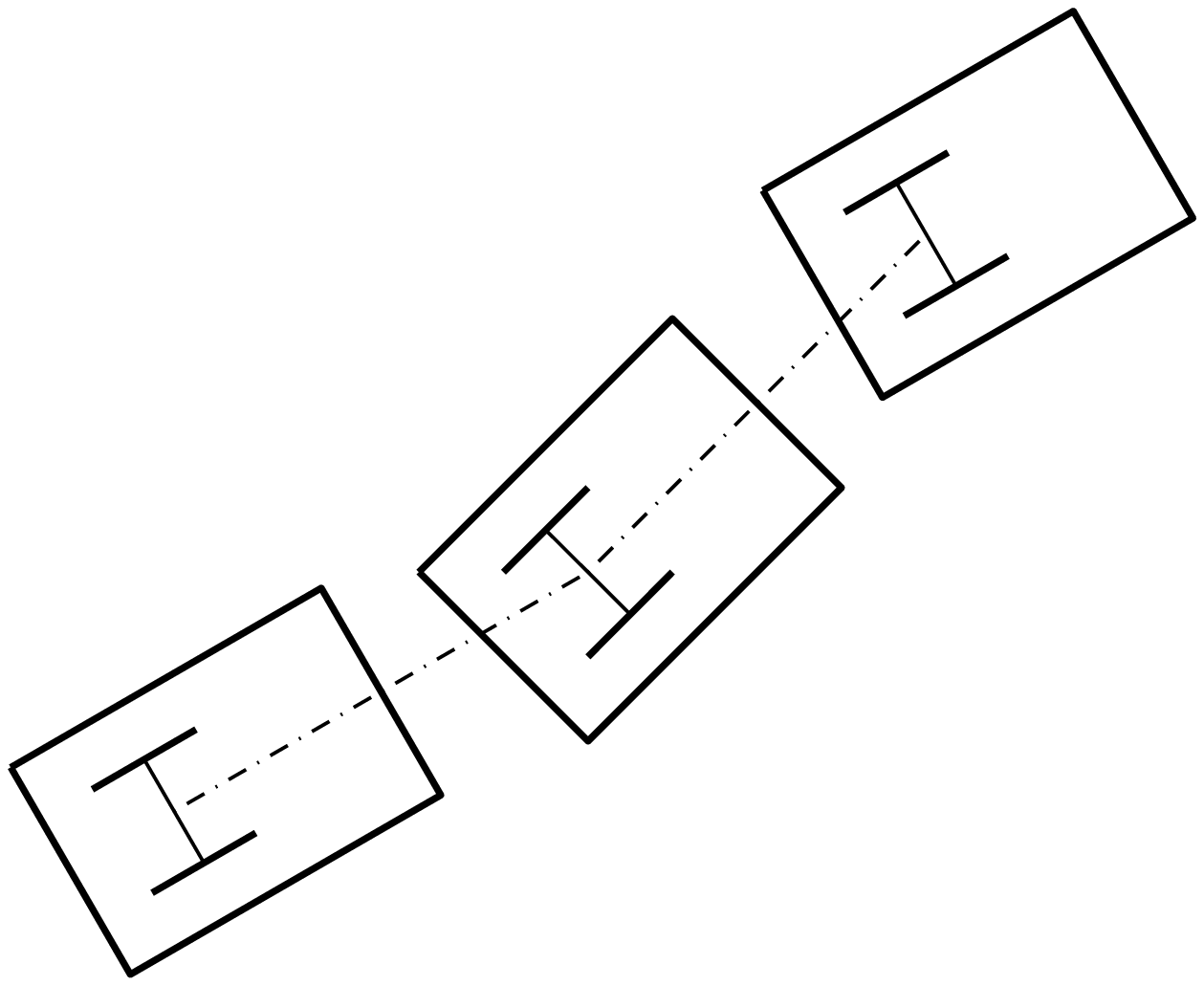}
\end{center}
\end{minipage}
\hspace{0.1\hsize} \begin{minipage}[h]{0.3\hsize}
\begin{center}
$$
\begin{array}{ccl}
\dot x_1 & = & u_1 \\
\dot x_2 & = & x_1\,u_2 \\
\dot x_3 & = & x_2\,u_2 \\
\dot x_4 & = & x_3\,u_2 \\
\dot x_5 & = & u_2
\end{array}
$$
\end{center}
\end{minipage}
\caption{A two-trailer and its normal form. Growth vector: $(2,3,4,5)$.
Singularity type: $a_{0}a_{0}$.}%
\end{figure}

\bigskip

\begin{figure}[h]
\begin{minipage}[h]{0.5\hsize}
\begin{center}
\includegraphics[width=\hsize]{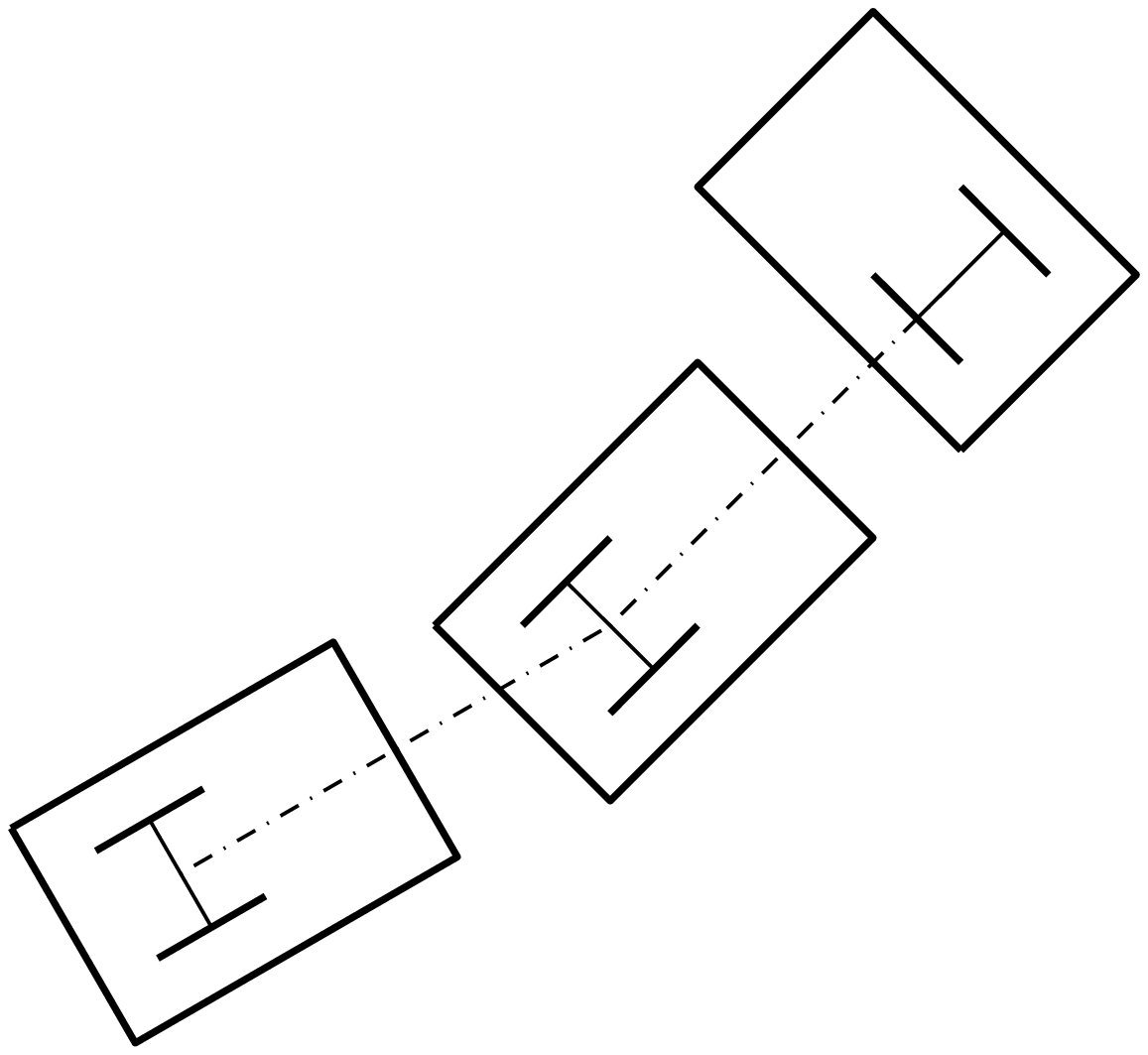}
\end{center}
\end{minipage}
\hspace{0.1\hsize} \begin{minipage}[h]{0.3\hsize}
\begin{center}
$$
\begin{array}{ccl}
\dot x_1 & = & u_1 \\
\dot x_2 & = & u_2 \\
\dot x_3 & = & x_1\,x_2\,u_2 \\
\dot x_4 & = & x_1\,x_3\,u_2 \\
\dot x_5 & = & x_1\,u_2
\end{array}
$$
\end{center}
\end{minipage}
\caption{A two-trailer and its normal form. Growth vector: $(2,3,4,4,5)$.
Singularity type: $a_{0}a_{1}$.}%
\label{fig-rigid-trailer-R5}%
\end{figure}
\end{center}

\newpage

\subsection{The Three-Trailer System}

\begin{center}
\begin{figure}[h]
\begin{minipage}[h]{0.5\hsize}
\begin{center}
\includegraphics[width=\hsize]{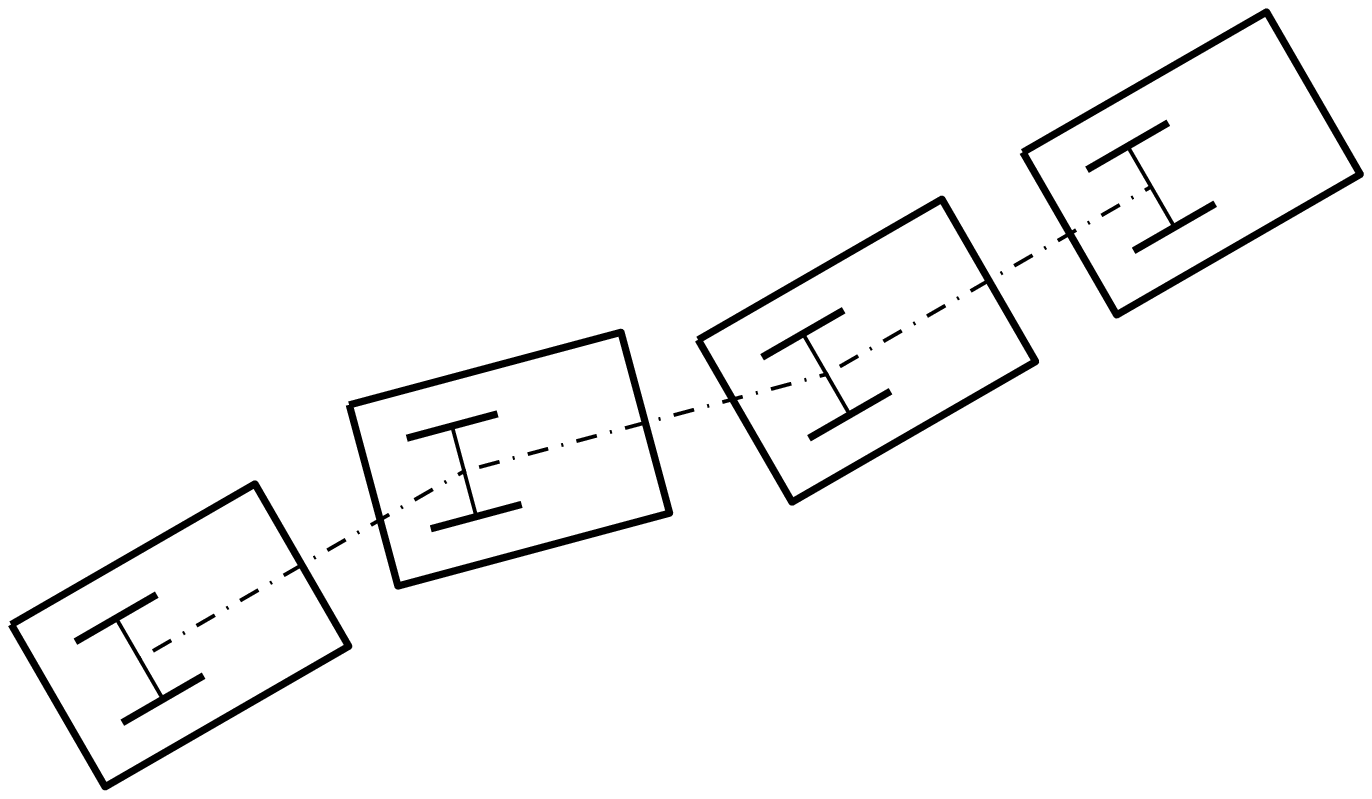}
\end{center}
\end{minipage}
\hspace{0.1\hsize} \begin{minipage}[h]{0.3\hsize}
\begin{center}
$$
\begin{array}{ccl}
\dot x_1 & = & u_1 \\
\dot x_2 & = & x_1\,u_2 \\
\dot x_3 & = & x_2\,u_2 \\
\dot x_4 & = & x_3\,u_2 \\
\dot x_5 & = & x_4\,u_2 \\
\dot x_6 & = & u_2
\end{array}
$$
\end{center}
\end{minipage}
\caption{A two-trailer and its normal form. Growth vector: $(2,3,4,5,6)$.
Singularity type: $a_{0}a_{0}a_{0}$.}%
\end{figure}

\bigskip

\begin{figure}[h]
\begin{minipage}[h]{0.5\hsize}
\begin{center}
\includegraphics[width=\hsize]{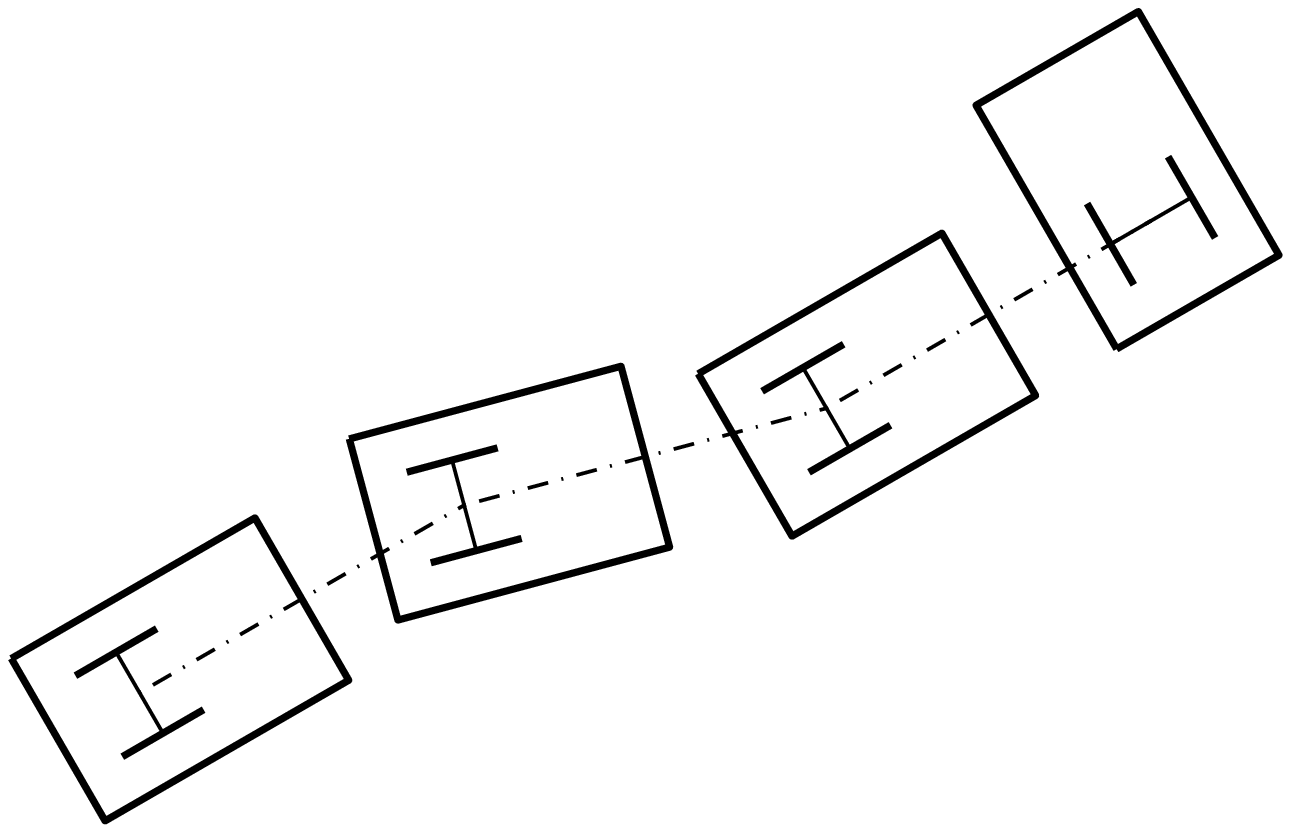}
\end{center}
\end{minipage}
\hspace{0.1\hsize} \begin{minipage}[h]{0.3\hsize}
\begin{center}
$$
\begin{array}{ccl}
\dot x_1 & = & u_1 \\
\dot x_2 & = & u_2 \\
\dot x_3 & = & x_1\,x_2\,u_2 \\
\dot x_4 & = & x_1\,x_3\,u_2 \\
\dot x_5 & = & x_1\,x_4\,u_2 \\
\dot x_6 & = & x_1\,u_2
\end{array}
$$
\end{center}
\end{minipage}
\caption{A two-trailer and its normal form. Growth vector: $(2,3,4,4,5,5,6)$.
Singularity type: $a_{0}a_{0}a_{1}$.}%
\end{figure}

\newpage

\bigskip

\begin{figure}[h]
\begin{minipage}[h]{0.45\hsize}
\begin{center}
\includegraphics[width=\hsize]{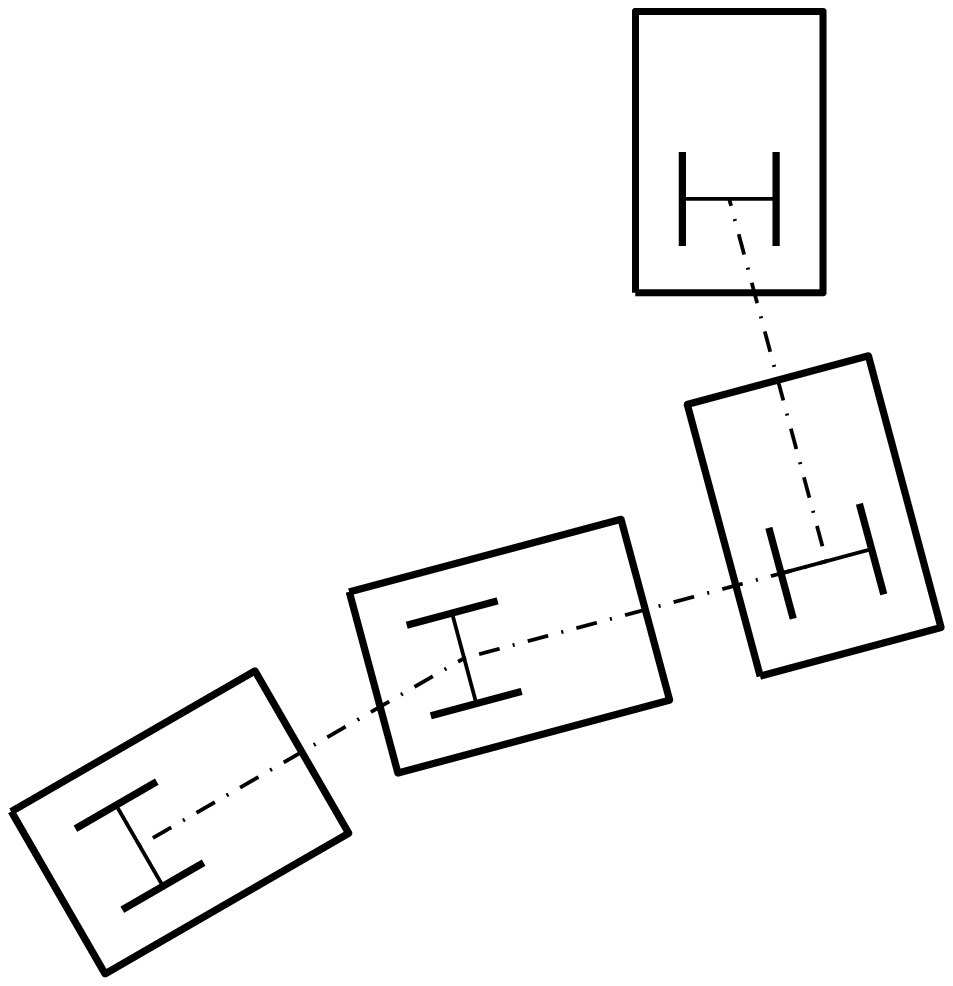}
\end{center}
\end{minipage}
\hspace{0.1\hsize} \begin{minipage}[h]{0.3\hsize}
\begin{center}
$$
\begin{array}{ccl}
\dot x_1 & = & u_1 \\
\dot x_2 & = & (x_1+1)\,u_2 \\
\dot x_3 & = & u_2 \\
\dot x_4 & = & x_2\,x_3\,u_2 \\
\dot x_5 & = & x_2\,x_4\,u_2 \\
\dot x_6 & = & x_2\,u_2
\end{array}
$$
\end{center}
\end{minipage}
\caption{A two-trailer and its normal form. Growth vector: $(2,3,4,5,5,6)$.
Singularity type: $a_{0}a_{1}a_{0}$.}%
\end{figure}

\bigskip

\begin{figure}[h]
\begin{minipage}[h]{0.45\hsize}
\begin{center}
\includegraphics[width=\hsize]{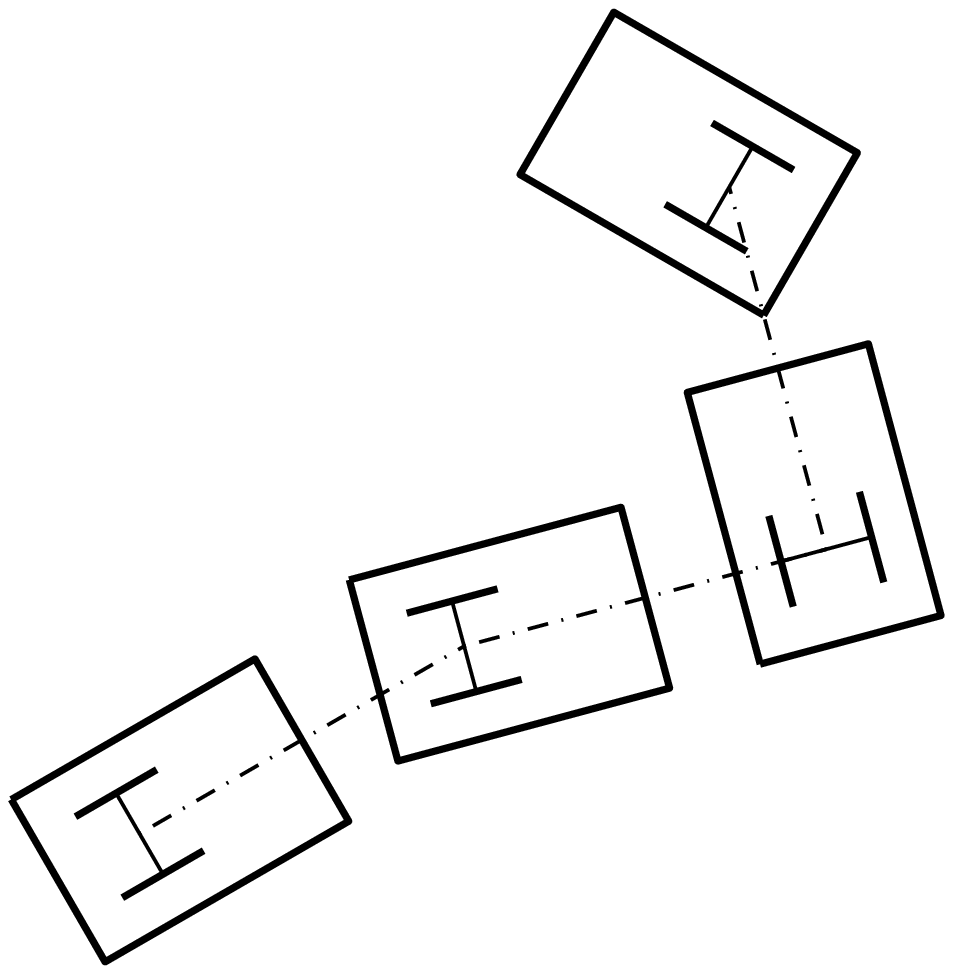}
\end{center}
\end{minipage}
\hspace{0.1\hsize} \begin{minipage}[h]{0.3\hsize}
\begin{center}
$$
\begin{array}{ccl}
\dot x_1 & = & u_1 \\
\dot x_2 & = & x_1\,u_2 \\
\dot x_3 & = & u_2 \\
\dot x_4 & = & x_2\,x_3\,u_2 \\
\dot x_5 & = & x_2\,x_4\,u_2 \\
\dot x_6 & = & x_2\,u_2
\end{array}
$$
\end{center}
\end{minipage}
\caption{A two-trailer and its normal form. Growth vector: $(2,3,4,5,5,5,6)$.
Singularity type: $a_{0}a_{1}a_{2}$.}%
\end{figure}

\newpage

\bigskip

\begin{figure}[h]
\begin{minipage}[h]{0.45\hsize}
\begin{center}
\includegraphics[width=\hsize]{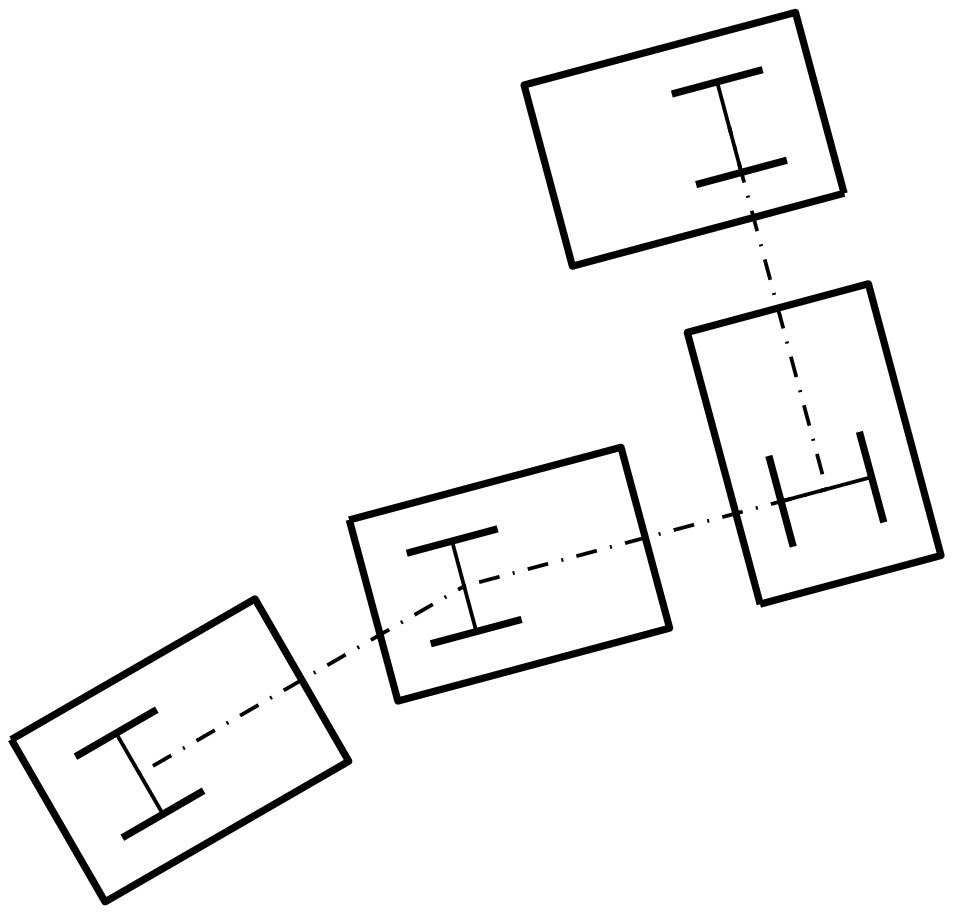}
\end{center}
\end{minipage}
\hspace{0.1\hsize} \begin{minipage}[h]{0.3\hsize}
\begin{center}
$$
\begin{array}{ccl}
\dot x_1 & = & u_1 \\
\dot x_2 & = & u_2 \\
\dot x_3 & = & x_1\,u_2 \\
\dot x_4 & = & x_1\,x_2\,x_3\,u_2 \\
\dot x_5 & = & x_1\,x_2\,x_4\,u_2 \\
\dot x_6 & = & x_1\,x_2\,u_2
\end{array}
$$
\end{center}
\end{minipage}
\caption{A two-trailer and its normal form. Growth vector: $(2,3,4,4,5,5,5,6)
$. Singularity type: $a_{0}a_{1}a_{1}$.}%
\end{figure}
\end{center}

\end{document}